\documentclass[sn-mathphys]{sn-jnl}

\usepackage{subcaption}
\usepackage{bm}

\begin{document}
\jyear{2022}%

\raggedbottom

\captionsetup[figure]{labelfont={bf},name={Fig.},labelsep=period}

\title{Multi-compartment human head modeling: generating adaptive tetrahedral mesh with GPU acceleration}

\author*[1]{\fnm{Fernando} \sur{Galaz Prieto}}\email{fernando.galazprieto@tuni.fi}
\author[1]{\fnm{Joonas} \sur{Lahtinen}}\email{joonas.lahtinen@tuni.fi}
\author[1]{\fnm{Maryam} \sur{Samavaki}}\email{maryamolsadat.samavaki@tuni.fi}
\author[1]{\fnm{Sampsa} \sur{Pursiainen}}\email{sampsa.pursiainen@tuni.fi}

\affil[1]{\orgdiv{Computing Sciences}, \orgname{Faculty of Information Technology and Communication Sciences}, \orgaddress{\street{Korkeakoulunkatu 3}, \city{Tampere}, \postcode{33014}, \state{Pirkanmaa}, \country{Finland}}}

\abstract{
\textbf{Purpose:} This paper introduces a highly adaptive and fast approach for generating a finite element (FE) discretization mesh for a given multi-compartment human head model obtained through a magnetic resonance imaging (MRI) dataset. The goal is to create accurate deep brain structures for electroencephalographic (EEG) source localization and applications.

\noindent \textbf{Methods:} We present a recursive application of solid angle labeling for surface segmentation with an adaptive scheme, i.e., a set of smoothing, inflation, and optimization tasks to enhance the mesh quality. 

\noindent \textbf{Results:} Our approach can produce a FE mesh with an accuracy greater than 1.0 millimeters, a significant milestone for 3D structure discretization and EEG source localization estimation.

\noindent \textbf{Conclusions:} Our method, implemented in the Matlab-based Zeffiro Interface, can manage the labeling aspect remarkably well to achieve human head FE meshes with complex deep brain structures using a time-effective parallel computing system.
}

\keywords{Mesh Generation, Finite Element Method (FEM), Inverse Problem, Boundary Conditions, Electroencephalography (EEG)}


\maketitle

\section{Introduction}
\label{sec:intro}
In this paper, we present a highly-adaptive volumetric Finite Element (FE) mesh \cite{de_munck_wolters_clerc_2012, rullmann2009eeg, braess2007finite} framework ---using magnetic resonance imaging (MRI) data as reference--- for generating realistic multi-compartment human head models. A critical aspect from this mesh generation process is to include forebrain structures, including the thalamus and the cerebrum, for studying accurate electromagnetic fields and their propagation across these compartments. In inverse modeling and imaging applications, solutions to the Maxwell's equations, given a source on the electromagnetic field, can be interpreted as a forward mapping. We evaluate the properties of the mesh generated using a simulated non-invasive Electroencephalography (EEG) \cite{niedermeyer2004} function wherein the brain activity is monitored through the potential in the electric field via a set of contact electrodes attached to the skin of the subject. EEG inverse problems are considered ill-posed \cite{kaipio2004}, i.e., solutions are not unique, are sensitive to noise and modeling errors, and strongly rely on the solution accuracy of the EEG forward problem \cite{hallez2007review, pursiainen2016electroencephalography, medani2015fem}. Stated in \cite{seeber2019subcortical, cuffin2001realistically, rezaei2020randomized}, sub-cortical components are known to be exceptionally sensitive and indistinguishable using only non-invasive data. Therefore, effective simulation of electric potentials throughout the head volume conductor model is essential.

The Finite Element Method (FEM) \cite{braess2007finite} is an advantageous method for modeling the head and its electromagnetic fields since it enables both surfaces- and volumetric-based mesh adaptations. FE mesh generation for the human head poses a challenging mathematical modeling task as the head includes various internal layers with complex geometrical properties. Well-known heuristic FE mesh generators, such as TetGen \cite{hang2015tetgen}, Netgen \cite{schoberl1997netgen}, and Gmsh \cite{geuzaine2009gmsh}, primarily aim at reconstructing non-intersecting surfaces. However, these generators might not deem applicable for reconstructing complex brain geometry in their standard form. Therefore, a tailored generator is required to produce high-quality brain and head FE mesh. Earlier solutions for optimized FE mesh generation include, e.g., SimBio-Vgrid \footnote{\url{http://vgrid.simbio.de}}, extensively applied in FEM-based EEG studies \cite{dannhauer2011modeling, antonakakis2019effect}, and coupling the TetGen mesh generator with other tools such as the implementation in Matlab's iso2mesh toolbox \cite{fang2009tetrahedral}. Through this study, we aim to provide an effective and automated alternative to perform a complete FE meshing process for an arbitrary number of cortical and sub-cortical compartments without the need for mesh pre-processing or combining different tools, such as in the recent study \cite{piastra2021comprehensive}. Our approach incorporates advanced features such as FE meshing of a complete head segmentation found by FreeSurfer Software Suite \cite{freesurfer}, robustness for intersecting surfaces, and a time-effective Graphic Processing Unit (GPU) accelerated recursive labeling process.

Our method to obtain an optimized FE mesh with an accuracy greater than 1.0 millimeter (mm) ---a crucial physical and physiological standpoint \cite{rullmann2009eeg}--- is approached through the following two stages: In the first stage (I), "mesh generation," we describe the design of a regular mesh, refinement of the tetrahedrons, and the application of recursive solid angle labeling for a multi-compartment surface segmentation \cite{lo2015finite}; as a result, this allows the mesh to match a set of closed segmentation boundaries, i.e., neurophysiological constraints in which the structure of the target brain model can be either (a) complex-shaped (e.g., cerebral cortex and sub-cortical nuclei), (b) thin-layered (skin), or (c) include material properties whose electrical conductivity can include a steep contrast (skin--skull--cerebrospinal fluid (CSF) boundaries) \cite{CHW:Dan2012}. The labeling process, in this regard, is the most predominant aspect in determining the computational cost of mesh generation and is effectively parallelizable either in a central or graphics processing unit (CPU or GPU). The second stage (II), "post-processing," consists of a set of tasks, i.e., surface smoothing, inflation, and optimization tasks.

We measure the accuracy of the EEG lead field matrices obtained with adaptive vs.\ regular meshes and perform a series of EEG source localization experiments to assess the quality of the matrices. As a reference, we used a three-layered spherical Ary model, and an MRI-based multi-compartment realistic head model \cite{piastra_maria_carla_2020_3888381} obtained using FreeSurfer. For source localization \cite{michel2019eeg}, we used the minimum norm estimate (MNE) \cite{CHW:Ham94}, the standardized low-resolution brain electromagnetic tomography (sLORETA) \cite{PascualMarqui2002}, and the dipole scan method, also known as the deviation scan \cite{fuchs1998improving} (the simplest form of beamformer). For result corroboration, we examine the 14.0 and 22.0 ms (milliseconds) median nerve stimulus responses reconstructed from experimental Somatosensory Evoked Potential (SEP) data \cite{piastra_maria_carla_2020_3888381, rezaei2021reconstructing} using a combination of sLORETA and gaussian mixture modeling (GMM) \cite{murphy2012machine}. 

Our results indicate that the highly-adaptive approach can produce robust FE meshes considering the 3D structure of the outcome of the discretization and the source localization estimates of EEG. Our algorithm, implemented in the open Matlab-based Zeffiro Interface forward and inverse toolbox\footnote{\url{https://github.com/sampsapursiainen/zeffiro_interface}} \cite{he2020zeffiro, rezaei2021reconstructing, miinalainen2019realistic}, whose forward solvers built upon high-resolution regular meshing approaches, demonstrates that mesh adaptation can considerably improve accuracy without significantly increasing computational time due to recursive re-labeling and GPU acceleration. 

\section{Materials and Methods}

\subsection{Mesh generation}
\label{sec:nested_labeling}
To obtain an appropriate FE mesh \cite{lo2015finite} for an arbitrary surface-based head tissue segmentation, we implement a hierarchical approach for filling the complete computation domain $\Omega$ where the tissue compartments are assumed to be nested. Layer labeling corresponding to the presented mesh takes precedence from the innermost to the outermost one. Whenever an intersecting segmentation boundary occurs, the one that comes first in the layer hierarchy will take priority. Since the number of compartments included in the model can be arbitrary, a temporal bounding box, behaving as the outermost compartment, ensure any adverse shrinkage effects that might take place on the external surface during the optimization steps.

\subsubsection{Regular meshing}
\begin{figure}[h!]
    \centering
    \begin{scriptsize}
        \begin{minipage}{3.0cm}
        \centering
            \begin{subfigure}{1.3cm}
            \centering
                \includegraphics[width = 1.2cm]{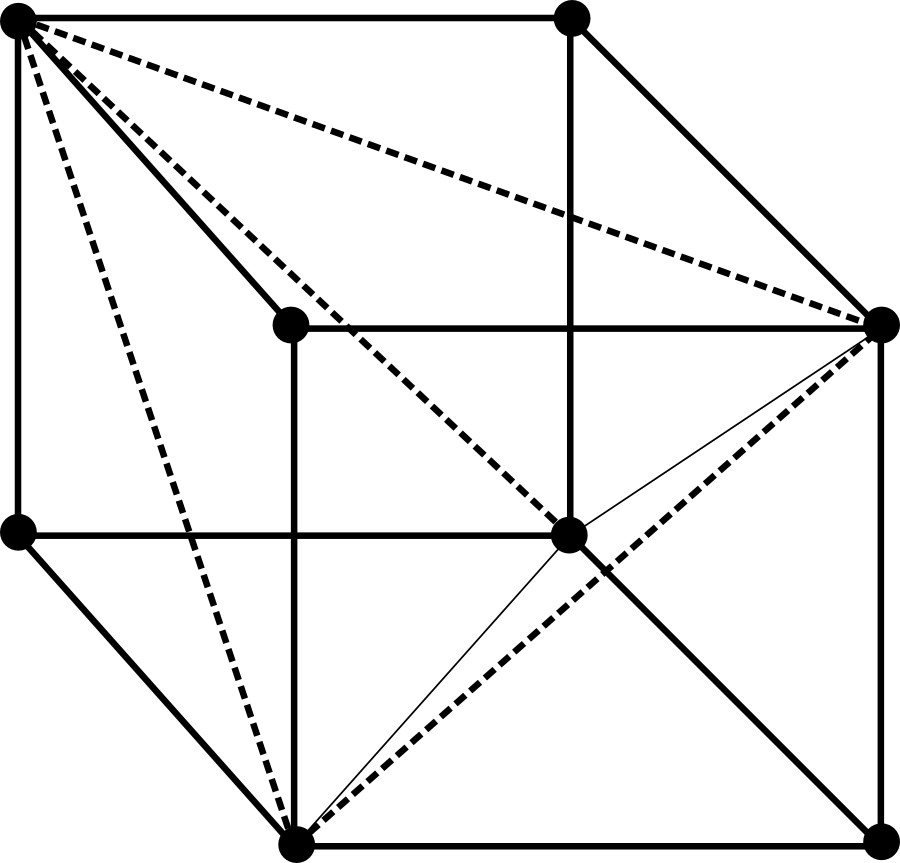}
                \caption{}
                \label{fig:subdivision_a}
            \end{subfigure}
            \begin{subfigure}{1.3cm}
            \centering
                \includegraphics[width = 1.2cm]{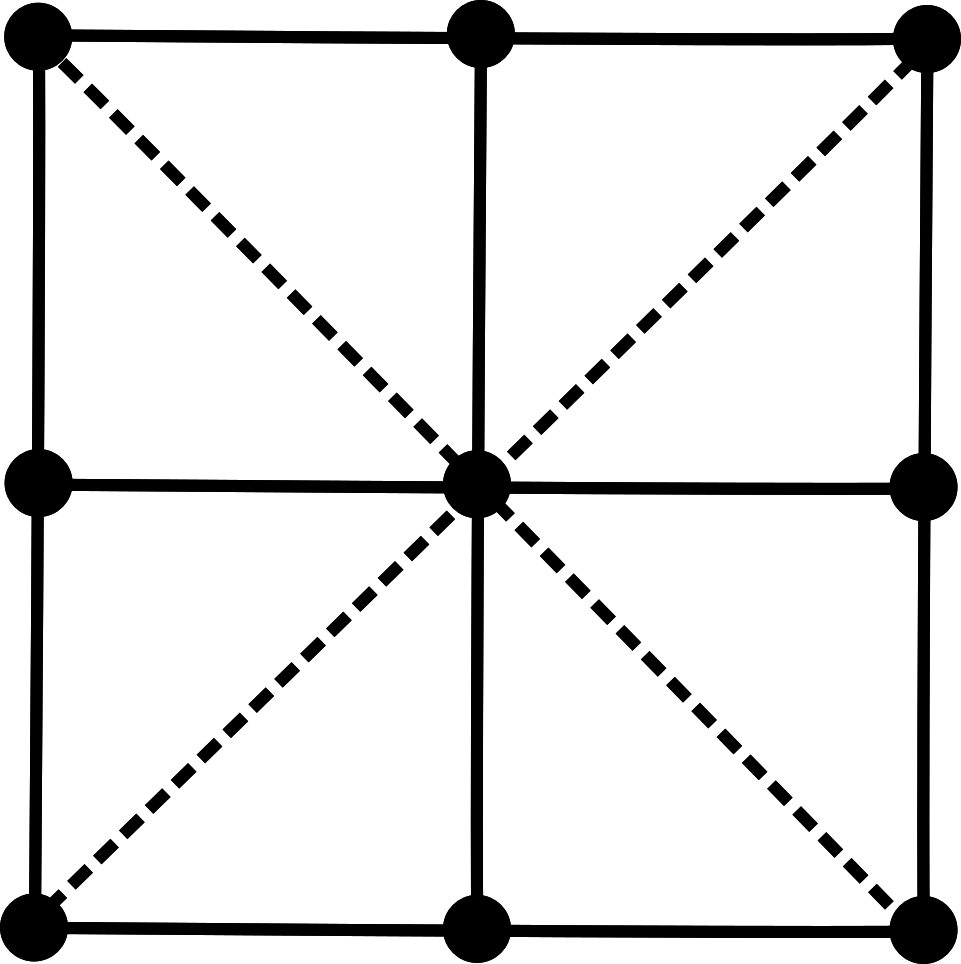}
                \caption{}
                \label{fig:subdivision_b}
            \end{subfigure}
            \\
            \textit{5 subdivisions}
        \end{minipage}
        \begin{minipage}{3.0cm}
        \centering
            \begin{subfigure}{1.3cm}
            \centering
                \includegraphics[width = 1.2cm]{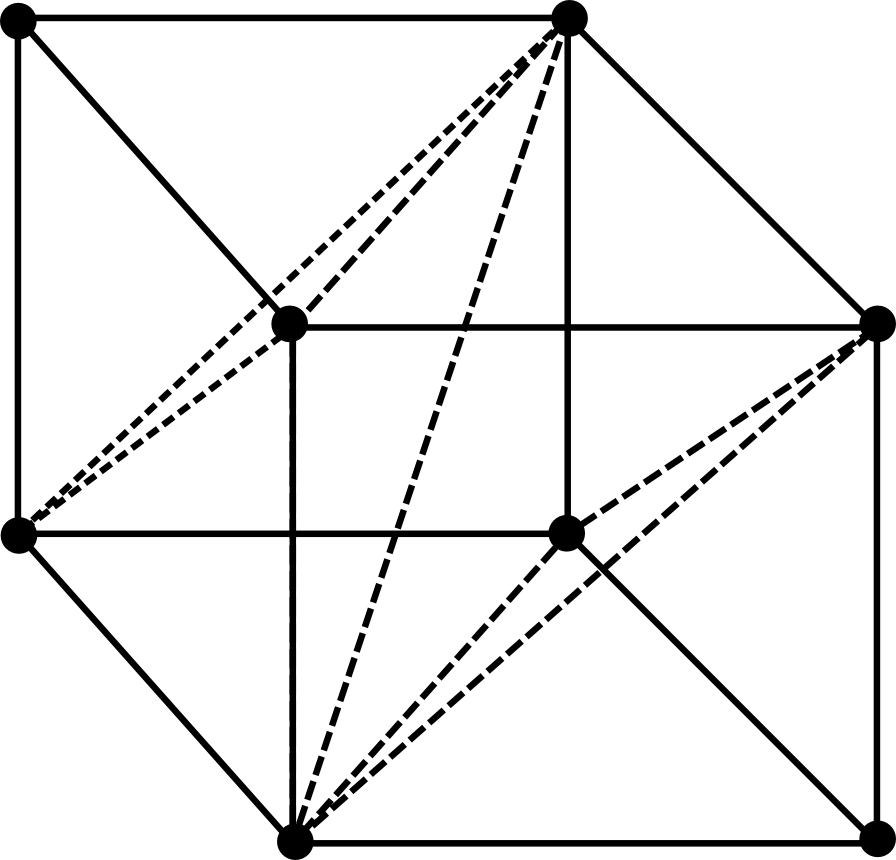}
                \caption{}
                \label{fig:subdivision_c}
            \end{subfigure} 
            \begin{subfigure}{1.3cm}
            \centering
                \includegraphics[width = 1.2cm]{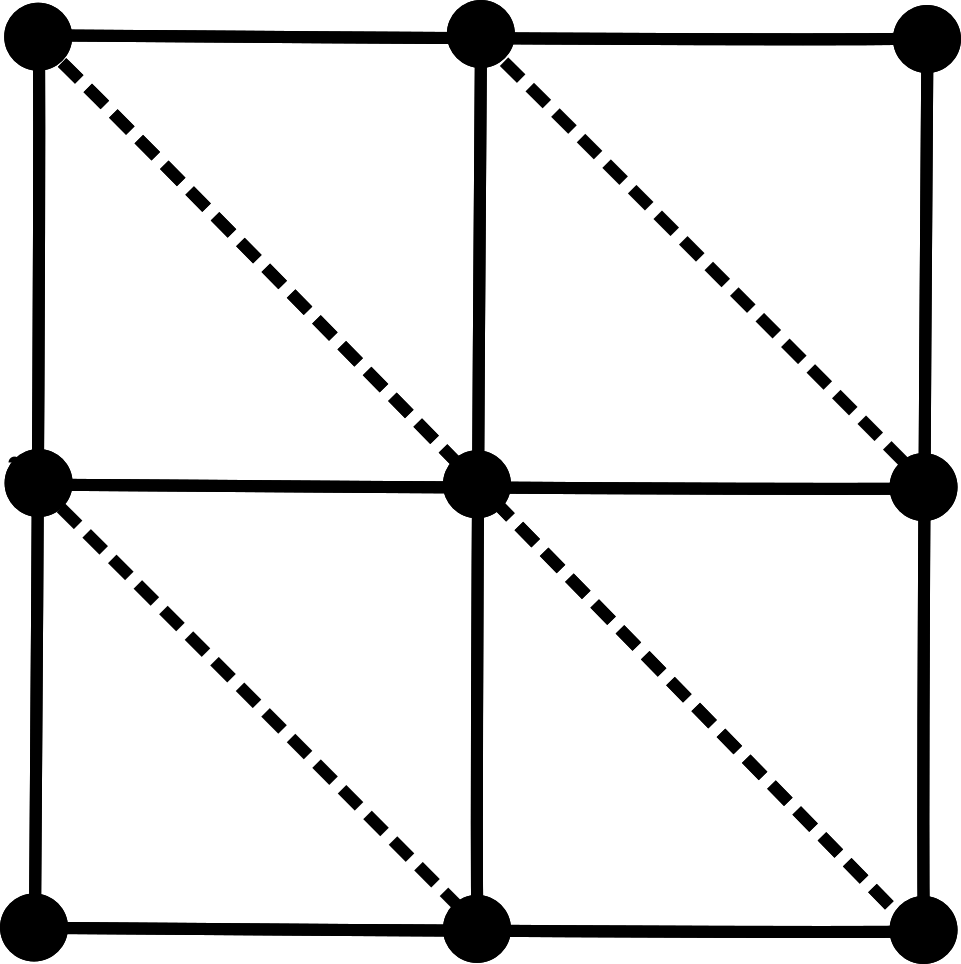}
                \caption{}
                \label{fig:subdivision_d}
            \end{subfigure}
           \\
           \textit{6 subdivisions}
        \end{minipage}
    \end{scriptsize}
\caption{Schematic illustration of an (\ref{fig:subdivision_a}) hexahedra subdivided into five tetrahedra \cite{pellerin2018there}. We chose this subdivision strategy as the resulting (planar) alternating pattern of faces and edges in adjacent hexahedra (\ref{fig:subdivision_b}) is advantageous to reduce lattice-dependent deficiencies in the mesh. In contrast, a subdivision into six tetrahedra (\ref{fig:subdivision_c}) results into a non-alternating pattern creating potential geometrical bias, e.g., directional bias of smoothing (\ref{fig:subdivision_d}).}
\label{fig:regular_meshing}
\end{figure}

The meshing process initiates by setting a regular hexahedral lattice in the domain. The hexahedral is subdivided into five tetrahedral elements (See Figure \ref{fig:regular_meshing}), the lowest subdivisions possible \cite{pellerin2018there}. This subdivision provides several benefits: it leads towards an alternating direction of the diagonal sub-dividing the edges over the adjacent hexahedra, restricts the generation of weighted directional bias affecting further mesh manipulation operations (e.g., smoothing, inflation, optimization), and minimizes the total amount of computer memory required to support the dataset.

\subsubsection{Mesh refinement}
\label{sec:mesh_refinement}
\begin{figure}[h!]
    \centering
    \begin{scriptsize}
        \begin{minipage}{3.0cm}
        \centering
            \begin{subfigure}{1.3cm}
            \centering
                \includegraphics[width = 1.2cm]{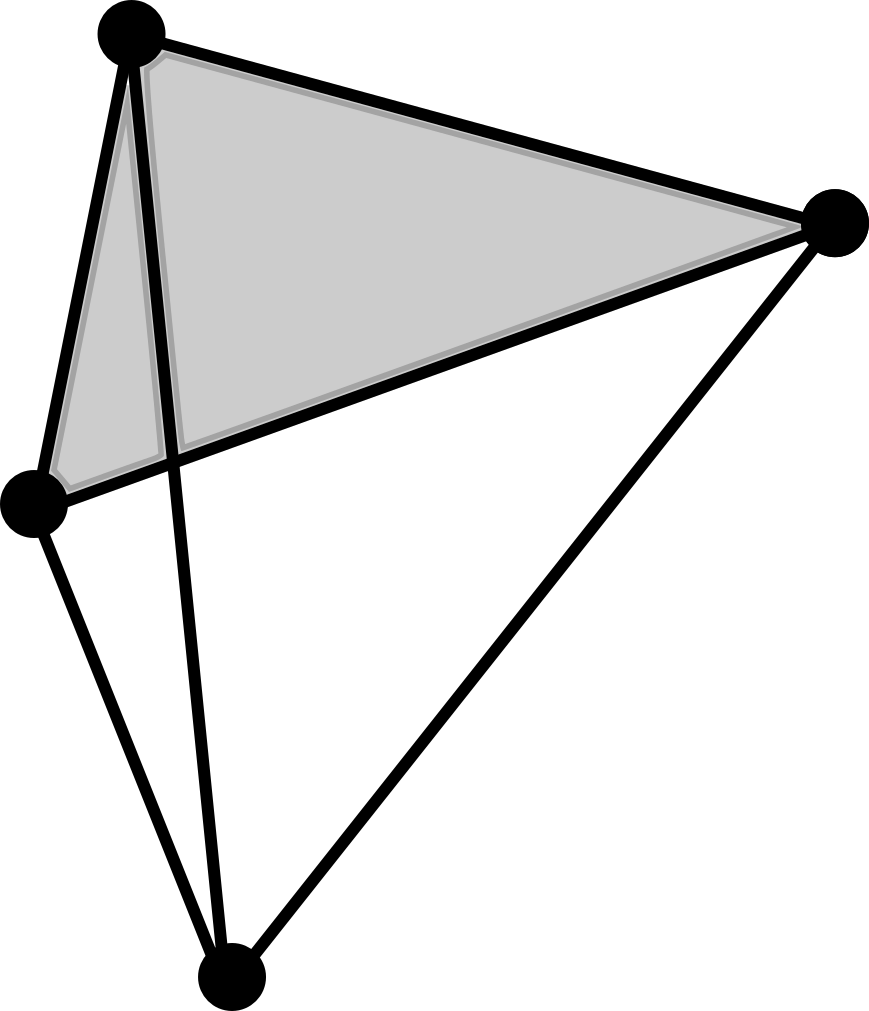}
                \caption{}
                \label{fig:mesh_refinement_a}
            \end{subfigure}
            \begin{subfigure}{1.3cm}
            \centering
                \includegraphics[width = 1.2cm]{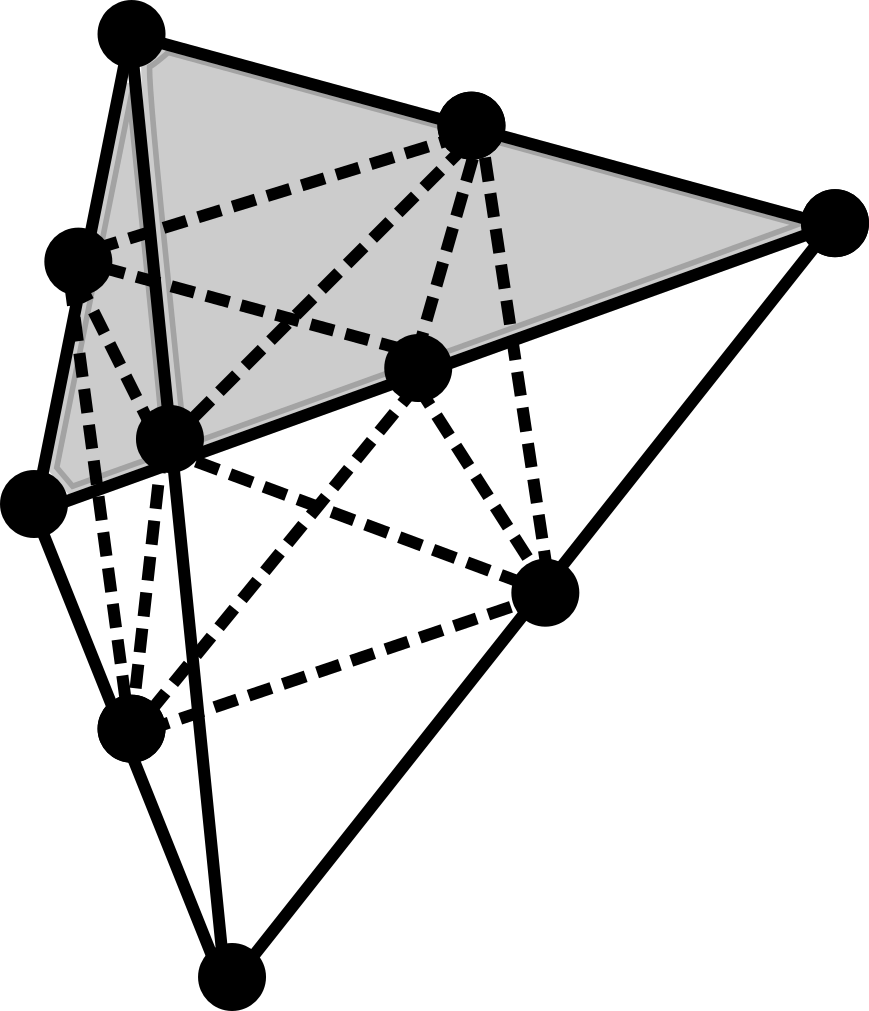}
                \caption{}
                \label{fig:mesh_refinement_b}
            \end{subfigure} 
            \textit{Volume refinement}
        \end{minipage}
        \begin{minipage}{4.3cm}
        \centering
            \begin{subfigure}{1.3cm}
            \centering
                \includegraphics[width = 1.2cm]{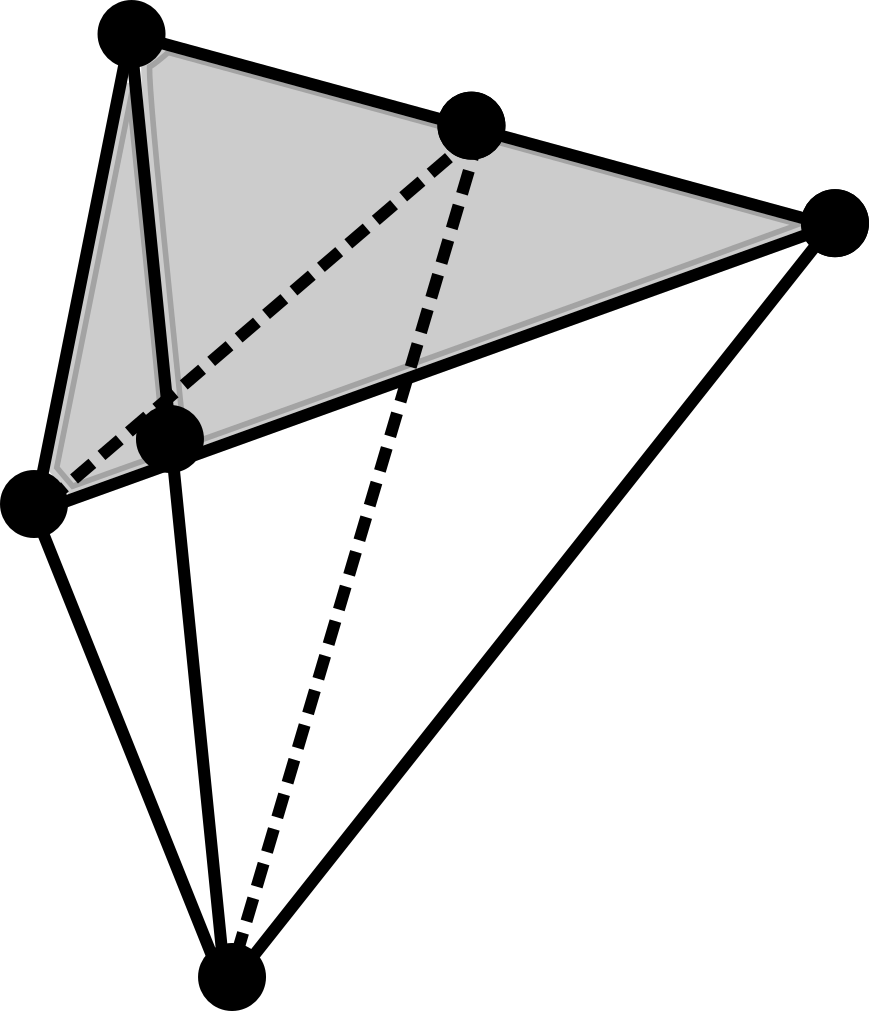}
                \caption{}
                \label{fig:mesh_refinement_c}
            \end{subfigure}
            \begin{subfigure}{1.3cm}
            \centering
                \includegraphics[width = 1.2cm]{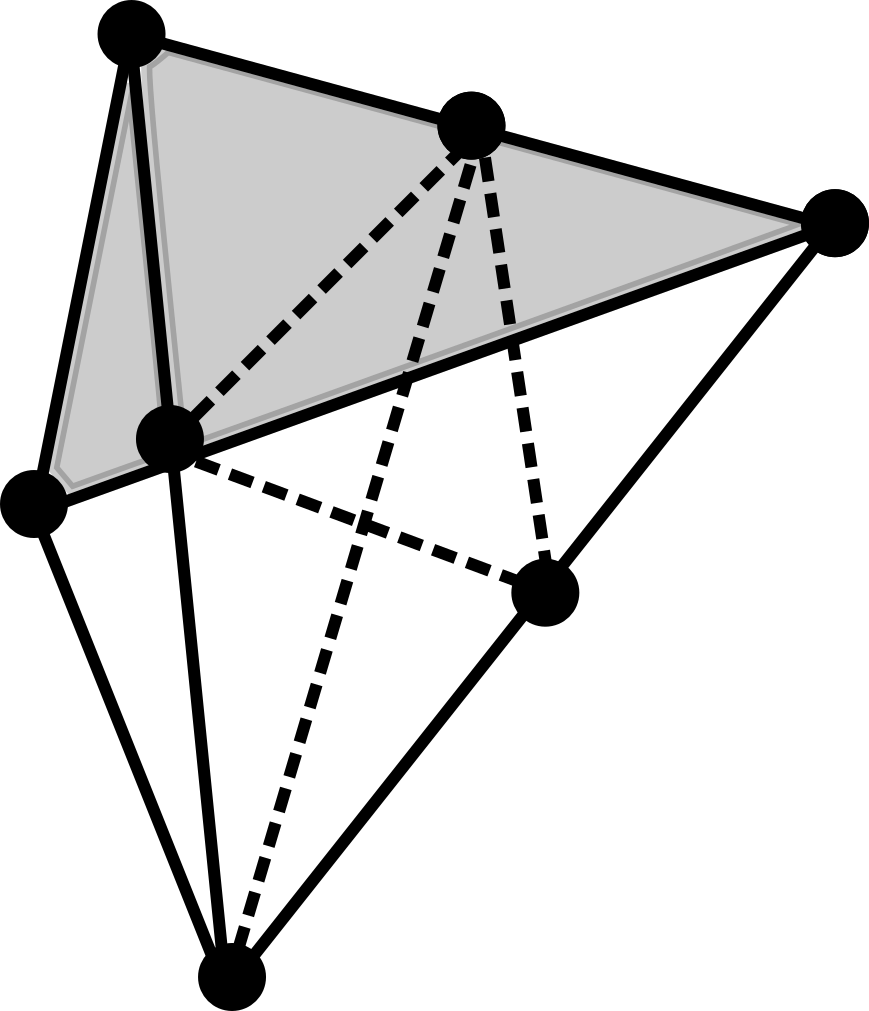}
                \caption{}
                \label{fig:mesh_refinement_d}
            \end{subfigure}
            \begin{subfigure}{1.3cm}
            \centering
                \includegraphics[width = 1.2cm]{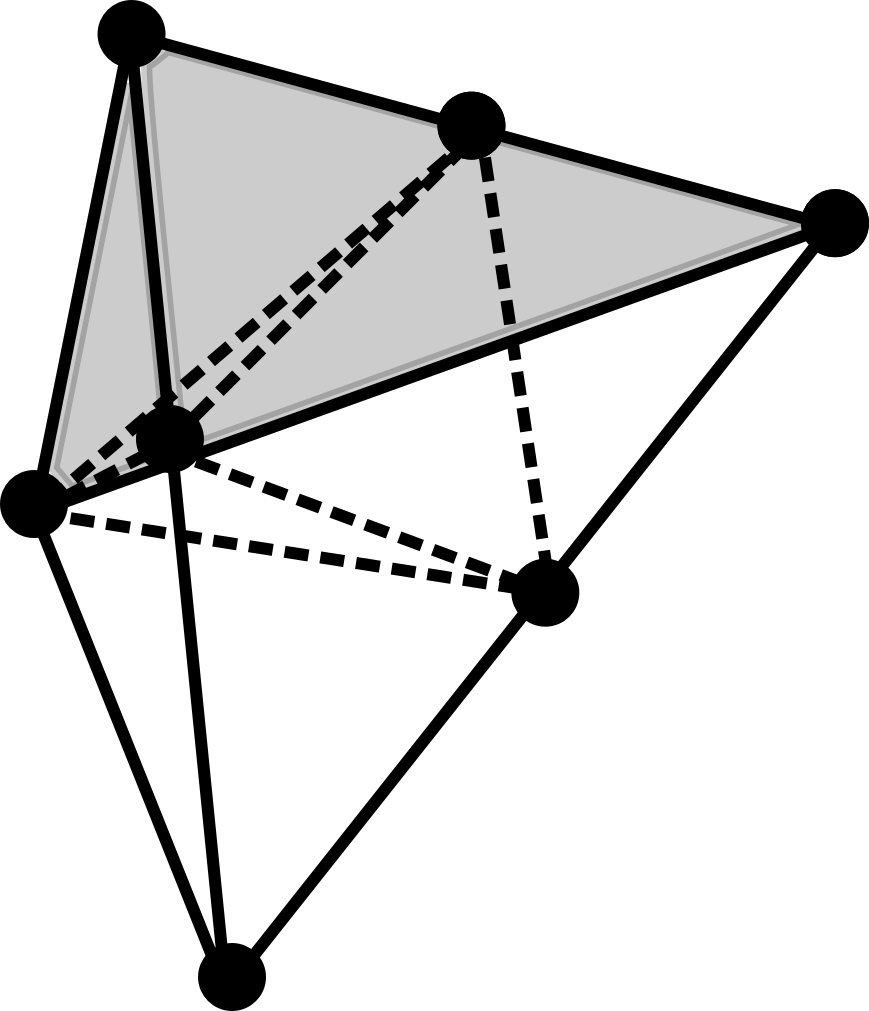}
                \caption{}
                \label{fig:mesh_refinement_e}
            \end{subfigure}
            \textit{Surface refinement}
        \end{minipage}
    \end{scriptsize}
\caption{Schematic illustration of the volume $V$, and surface $S$ refinement. In volume refinement, each tetrahedral element (\ref{fig:mesh_refinement_a}) in the given subdomain $V$ of $\Omega$ is subdivided into eight tetrahedra (\ref{fig:mesh_refinement_b}). The refined volume $V$ is connected to the remaining part of the domain $\Omega \setminus S$ via surface refinement where the elements can have one to three edges subdivided on the face adjacent to the boundary $\partial S$, leading to tetrahedral subdivisions (\ref{fig:mesh_refinement_c}, \ref{fig:mesh_refinement_d}, \ref{fig:mesh_refinement_e}), respectively. A refinement with respect to a boundary $\partial S$ can be defined as two similar refinements on both sides of $\partial S$.}
\label{fig:refinement}
\end{figure}

The mesh can be refined to a given subdomain $V$ or its boundary $\partial V$ (see Figure \ref{fig:refinement}). The elements in $V$ are subdivided into eight tetrahedra faces, while in $\Omega \setminus V$ ---those having nodes adjacent to the boundary--- are subdivided to their corresponding number of edges shared with $\partial V$.

\subsubsection{Labeling}
\label{sec:labeling}
\begin{figure}[h!]
    \centering
    \begin{scriptsize}
        \begin{subfigure}{2.5cm}
            \centering
            \includegraphics[width = 2.2cm]{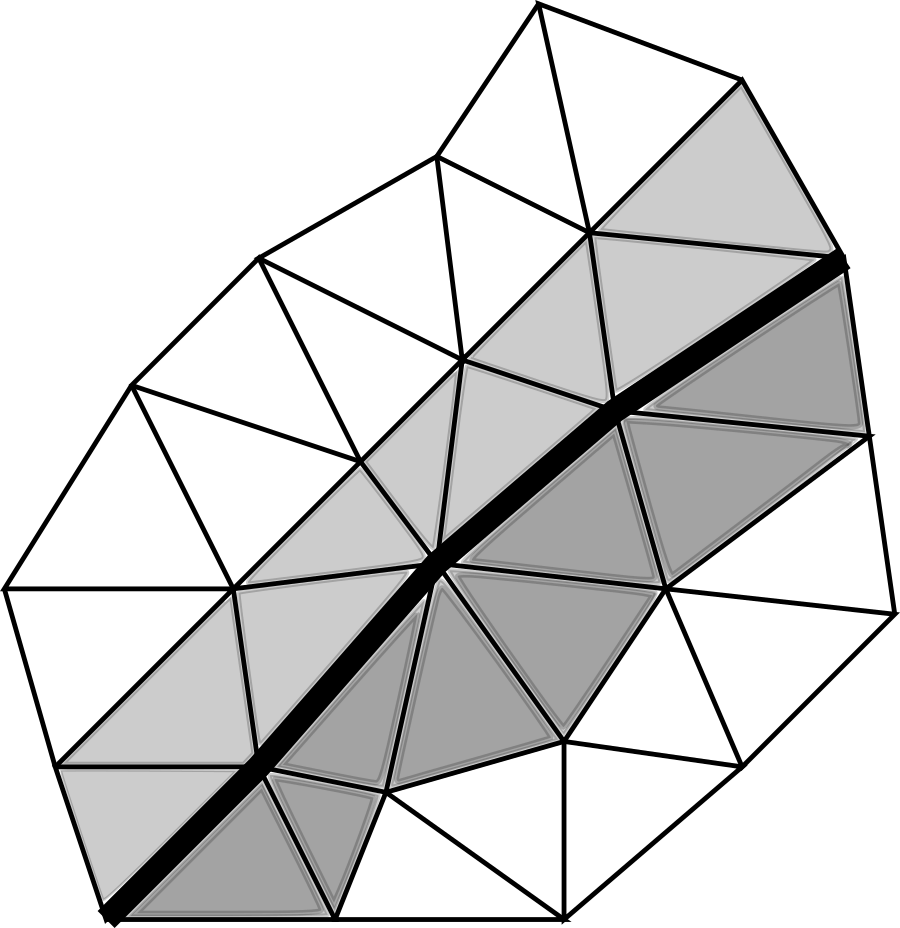}
            \caption{}
            \label{fig:labeling_a}
        \end{subfigure}
        \begin{subfigure}{2.5cm}
            \centering
            \includegraphics[width = 2.2cm]{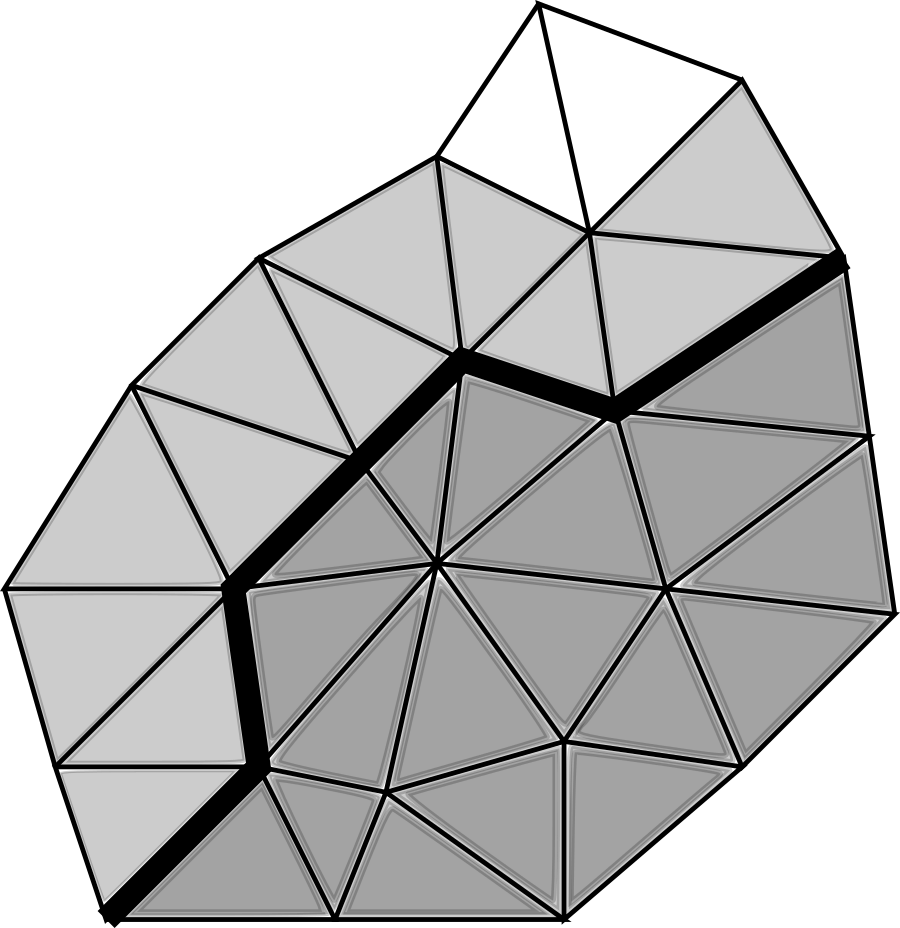}
            \caption{}
            \label{fig:labeling_b}
        \end{subfigure}
    \end{scriptsize}
    \caption{Recursive solid angle labeling process re-labels the tetrahedra (triangles) adjacent to a given boundary $\partial V$ (black curve) between two different tissue compartments (light and dark grey). The re-labeling process is performed until $\partial V$ does not change between two consecutive recursion steps, as it does between (\ref{fig:labeling_a}) and (\ref{fig:labeling_b}).}
    \label{fig:labeling}
\end{figure}

To obtain the compartment of a given element, we apply the solid angle labeling method \cite{lo2015finite}, i.e., a mesh node at position $\vec{r}_i$ associated with the integral
\begin{equation}
s_i = \frac{1}{4 \pi} \int_{\partial S} (\vec{r} - \vec{r}_i) \cdot \vec{n} \,  \hbox{d} A, \end{equation}
where $\mathcal{S}$ is a closed segmentation boundary, $\hbox{d} A$ is a surface area differential, and $\vec{n}$ is the normal vector at the point $\vec{r}$. The integral $s_i$ defines the ratio between the solid angle and the angle covered by $\partial S$ with respect to $\vec{r}_i$, a point is enclosed by $\partial S$ if $s_i >= T$ with $T$ denoting a given threshold value between 0 and 1. Tetrahedrons with four nodes inside $\partial S$ are labeled as the elements contained within the surface $S$. The labeling aspect is applied as follows:
\begin{enumerate}
    \item Initial labeling: all nodes within a tetrahedral mesh are labeled using the solid angle method.
    \item Recursive re-labeling (see Figure \ref{fig:labeling}): the tissue boundaries in the existing mesh are re-labeled after changing their nodes, i.e., for each labeling-based tissue boundary, the solid angle integral $s$ is re-evaluated for the nodes shared by the tetrahedra adjacent to the tissue boundary. Thus, the labels of the tetrahedra are updated, potentially changing the tissue boundary. The process is repeated recursively until achieving convergence, i.e., when the status of one or more nodes remains unchanged. 
\end{enumerate}

Once the regular mesh has been generated, labeled, refined, and re-labeled, the resulting mesh is post-processed via smoothing, inflation, and optimization techniques to fit any tetrahedral compartment surface $\partial V$ included in the mesh with the corresponding closed segmentation boundary $\mathcal{S}$.

\subsection{Post processing}

\subsubsection{Smoothing}
\label{sec:smoothing}
As a smoothing technique, we use the Taubin's method \cite{pursiainen2012raviart} which performs alternating Laplacian forward and backward smoothing steps to smooth the mesh, and to reduce shrinkage, respectively. Both volume and tissue boundaries are smoothed using the following formula: 
\begin{eqnarray}
\vec{x}^{(k+1/2)}_i & = & \vec{x}^{(k)}_i + \lambda \sum_{j \in \mathcal{N}_i}  x^{(k)}_j, \\
\vec{x}^{(k + 1)}_i & = & \vec{x}^{(k)}_i - \mu \sum_{j \in \mathcal{N}_i}  x^{(k+1/2)}_j, 
\end{eqnarray}
where $\lambda < 1$ and $\mu < 1$ are the smoothing parameters, and $\mathcal{N}_i$ is the index set containing the $i$-th nodes together with their neighbors.

In volumetric smoothing, all the neighbors connected by an edge with $\vec{x}^{(k)}_i$ are included in $\mathcal{N}_i$. The stopping criterion for smoothing is set to be 
\begin{equation}
    \frac{\| \vec{x}^{(k+1/2)} - \vec{x}^{(k)}  \|_2}{\| \vec{x}^{(0)}  \|_2} \leq \xi, 
\end{equation}
where $\xi$ is a user-defined value.

\subsubsection{Inflation}
\label{sec:inflation}
After the labeling and re-labeling routines, the surface $\partial V$ is then enclosed by the closed segmentation boundary $\mathcal{S}$ as each tetrahedra adjacent to $\partial V$ has now its four nodes inside $\partial S$. Once four nodes have been enclosed, they are then lifted by inflation. This method is accomplished by initially detecting the intersecting edges with $\mathcal{S}$, followed by determining the intersection point, and finally sliding the nodes of $\partial V$ towards those points along their associated edges. The amount of inflation for a given intersecting edge is calculated via the formula: 
\begin{equation}
    \vec{x}^{(1)}_i = \vec{x}^{(0)}_i + \zeta d \vec{e}, 
\end{equation}
where $\vec{e}$ is a unit vector parallel to the edge, $d$ is the distance between $\vec{x}_i$ and $\partial S$, and $\zeta < 1$ is an inflation parameter. 

\subsubsection{Optimization}
\label{ref:Optimization}
After smoothing and inflation, the optimization step finalizes the post-processing in two phases, ensuring that the elements generated are valid and have an appropriate condition which we measure by the volume $\mathcal{V}_T$ of a given tetrahedron divided by the length of its longest edge $\ell_T^{\hbox{\scriptsize (max)}}$, i.e., 
\begin{equation}
    \kappa = \frac{\mathcal{V}_T}{\ell_T^{\hbox{\scriptsize (max)}}}.
\end{equation}
Here, $\kappa$ is the condition number. In the first phase, the optimization searches for inverted tetrahedrons containing elements with negative determinant. These elements are fixed by bringing any node $\vec{x}$ outside of the outer surface of its associated supernode back into the surface set containing all the neighbors sharing an edge with $\vec{x}$. The process repeats recursively until no negative determinants are found. The second phase executes Delaunay turns \cite{delaunay1934sphere, braess2007finite} for each pair of adjacent elements in which the condition number is smaller than a given threshold $\tau$, ($\kappa < \tau$), at least, in one of the two tetrahedrons. The tetrahedron with minimal condition number determines the orientation of the shared face.

\subsection{EEG modelling}
\label{ref:EEG_modelling}
In EEG, an unknown primary source current distribution ${\bf \vec{J}_p}$, generated by the brain activity, is to be localized given a set of point measurements ${\bf y}$ of the electric potential distribution $u$.

Given an electric potential $\vec{J}_p$ a solution to the second-order partial differential equation can be obtained as
\begin{equation} 
\label{eeg_equation}
\nabla \cdot (\sigma \nabla u) = \nabla \cdot \vec{J}_p,
\end{equation}
where ${\bf \sigma}$ is the electric conductivity distribution, affecting the solution significantly. 

Equation (\ref{eeg_equation}) can be expressed in a discretized form ${\bf A} {\bf z} = {\bf G} {\bf x}$, where the vectors ${\bf z}$ and ${\bf x}$ discretize the volumetric distributions of $u$ and $\vec{J}_p$, and the matrices ${\bf A}$ and ${\bf G}$ the operators acting on both left- and right-side, respectively. Solving the discretized form with respect to ${\bf z}$ and taking into account that ${\bf y}$ can be obtained by multiplying ${\bf z}$ with a restriction matrix ${\bf R}$, i.e., ${\bf y} = {\bf R} {\bf z}$, implies the linear forward model
\begin{equation} 
\label{forward_model}
{\bf y} = {\bf L} {\bf x},
\end{equation} 
where ${\bf L} ={\bf R} {\bf A}^{-1} {\bf G} $ is the lead field matrix. The equation (\ref{forward_model}) poses an ill-posed inverse problem with respect to the unknown activity ${\bf x}$. The source currents corresponding to the entries of ${\bf x}$ are assumed to be dipolar \cite{CHW:Ham93}.  

To obtain the lead field matrix, we apply the FEM along with the divergence conforming H(div) source model \cite{pursiainen2016electroencephalography} to approximate dipolar sources. The position-based optimization (PBO) technique \cite{bauer2015, pursiainen2016electroencephalography} provides us with the dipolar source's position and orientation. As for solving the inverse problem, we apply and compare three well-known alternative source localization methods: the minimum norm estimation (MNE), the Standardized low-resolution brain electromagnetic tomography (sLORETA), and the dipole scan methods.

\subsubsection{Minimum norm estimation}
\label{sec:MNE}
The Minimum Norm Estimate (MNE) \cite{CHW:Ham94} is evaluated as the solution to the following Tikhonov regularized minimization problem:
\begin{equation}
    \min_{{\bf x}}\left\lbrace \left \| {\bf C}^{-1/2}\left({\bf L}{\bf x}-{\bf y}\right) \right \|^2+\left \| \lambda \bm{\Sigma}^{-1/2}{\bf x} \right \|^2\right\rbrace,
    \label{tikhonov}
\end{equation}
where ${\bf C}$ is the noise covariance matrix, $\bm{\Sigma}$ is the prior covariance matrix, and ${\bf \lambda}$ is the penalty term. This problem can be associated with a Bayesian setting in which both the prior and likelihood have gaussian densities \cite{rezaei2021reconstructing}. The applied regularization parameter in this regard is understood as a prior-over-measurement signal-to-noise ratio (PM-SNR), i.e., the relative weight of the prior compared to the likelihood. 
 
\subsubsection{Standardized low-resolution brain electromagnetic tomography}
\label{sec:sLORETA}
The standardized low-resolution brain electromagnetic tomography (sLORETA) \cite{PascualMarqui2002} is a re-weighted MNE reconstruction, where the weighting coefficients are inversely proportional to the noise variances reconstructed via MNE  (\ref{tikhonov}). Consequently, the energy of the reconstruction is not affected by the variable amplitude of the lead field in different parts of the domain where the reconstruction is found.

\subsubsection{Dipole scan}
\label{sec:dipolescan}
The dipole scan \cite{fuchs1998improving} method evaluates the subsequent goodness-of-fit (GOF) $g$ between the data and a dipole source at the $i$-th position of the space as 
\begin{equation}
g  = 1 -\frac{ \| {\bf y}-{\bf L}_i {\bf L}_i^\dag {\bf y}\|^2}{\| {\bf y} \|^2}.
\end{equation}
Here, ${\bf L}_i$ is the submatrix of ${\bf L}$ containing the associated source(s) at the $i$-th position. The reconstruction is a GOF distribution covering the whole source space. 

\subsection{Computational domains}

\subsubsection{Spherical three-layer Ary model}
\label{sec:ary_sphere}
The three-layer Ary model \cite{ary1981location} is composed of three concentric spheres with radii of 87, 92, and 100 mm and conductivities of 0.33, 0.0042, and 0.33 S/m (Siemens per meter), respectively. The steep contrast between the adjacent brain and skull layers, analyzed using a universal heuristic mesh generator in \cite{pursiainen2011forward}, poses a challenge for forward simulation and dipole localization. Therefore, one ought to find the Berg parameters of the analytical solution. The numerical experiments used 180-point electrodes as reference evenly spread over the external surface.

\subsubsection{Realistic head model and segmentation}
\label{sec:head_segmentation}
The realistic head model was obtained from an openly available dataset\footnote{\url{https://doi.org/10.5281/zenodo.3888381}}  \cite{piastra_maria_carla_2020_3888381} including MRI data of a healthy 48 years old adult subject. The standard FreeSurfer Software Suite{\footnote{\url{https://surfer.nmr.mgh.harvard.edu/}}} reconstruction procedure segmented the data, providing 18 volumetric compartments in which the segmentation of the subcortical regions corresponded to FreeSurfer's Aseg Atlas. The tissue conductivities were set according to literature reference values \cite{CHW:Dan2012}: 0.33 S/m for scalp, 0.0064 S/m for skull, 1.79 S/m for CSF, 0.33 S/m for gray matter, 0.14 S/m for white matter. As for the subcortical compartments, 0.33 S/m was chosen \cite{shahid2014value}.

\subsection{Analytical solution}
\label{sec:analytical_solution}

\subsubsection{Earth mover's distance}
\label{sec:EMD}
To measure the accuracy of a source localization estimate for a given dipole, we use the earth mover's distance (EMD) \cite{RubnerY1998Amfd} which measures the minimum amount of work required to transfer a given mass distribution to another. The measure was first proposed as an analytical distance function between probability distributions in metric spaces in \cite{Kantorovich, Vaserstein}. We use EMD to measure the similarity between a distributional reconstruction and a dipole associated with Dirac's delta function (or unit impulse).

\subsubsection{Gaussian mixture modeling}
\label{sec:GMM}
For the realistic head model, we estimate source localization through clustering since there is no 'ground truth' for such a case. We use the gaussian mixture modeling (GMM) \cite{murphy2012machine} clustering technique, where the aim is to identify the main concentration areas of the primary current density as a mixture, i.e., a superposition of gaussian distributions. We assume that the current density estimate corresponds to a single or multimodal distribution of source activity, where the individual mean location of a distribution component indicates the possible location of the corresponding neural activity. By assigning the weight
\begin{equation}
    w_i = \left\|\vec{J}(\vec{ r}_i)\right\|_2^2\left(\sum_{j=1}^N\left\|\vec{J}(\vec{ r}_j)\right\|_2^2\right)^{-1},
\end{equation}
for a given reconstruction $\vec{J}$ and for each source point $\vec{r}_i$, $i = 1,2, \ldots, N$, the GMM can detect a controlled number of neural activity clusters contained within $\vec{J}$. The number of possible clusters are recursively added in the configuration through a deterministic procedure: by selecting the best fitting set for a probability-thresholded Mahalanobis distance \cite{mahalanobis1936generalized} not yet added in the configuration. The most fitting number of Gaussian components is decided through the GMM algorithm by selecting the model with the lowest Bayes Information Criterion (BIC).
 
\subsection{Numerical experiments}
\label{sec:NumExp}
This section describes the numerical experiments performed with the adaptive FEM mesh generation routine for source localization accuracy and EEG forward simulation in the three-layer spherical Ary model and realistic head model as computational domains. We create and compare mesh generation and post-processing effectiveness over these domains using a regular mesh created by dividing the hexahedral lattice size of 1.0 mm edge length into a tetrahedral lattice, and an adaptive mesh using hexahedral lattice sizes of 3.0, 2.0, and 1.3 mm.

\subsubsection{Finite element meshes}
To avoid redundant computations during the labeling phase, the FE mesh generation process first evaluates the density of the surface segmentation. If the calculated density is greater than the initial FE mesh resolution, it downsamples the segmentation to match the resolution of said FE mesh. The tissue boundaries of the volumetric mesh are then refined once for the compartments, including the cerebrum, cerebellum, and deep active nuclei. The smoothing parameters were set to be equal $\lambda = \mu = 0.4$, and the stopping criterion was set $\xi = 0.9$ for the volumetric smoothing and $\xi = 0.1$ for the surface smoothing. The inflation parameter was chosen to be $\zeta = 0.05$.

\subsubsection{Source space}
\label{sec:source_space}
Divergence conforming \cite{pursiainen2016electroencephalography} source model was used to generate each EEG lead field matrix with a set of 10,000 source positions distributed uniformly throughout the domain $\Omega$. Each source position in the spherical Ary model has a triplet source, each one corresponding to the three Cartesian orientations. The sources in the head model were normally constrained concerning the compartment surface following the orientation of the neurons and neural activity \cite{purves2019neuroscience}.

\subsubsection{Somatosensory evoked potential data}
\label{sec:seps}
The data package used in this study includes a 74-channel EEG measurement dataset from a somatosensory evoked potential (SEP) experiment. This experiment consisted of stimulating the median nerve of the right wrist using monophasic square-wave electrical pulses with 0.5-millisecond duration \cite{piastra_maria_carla_2020_3888381}. The analyzed SEPs in these measurements indicated that, through non-invasive stimulation methods, one can approximately localize the originators of both sub-cortical and cortical early responses \cite{rezaei2021reconstructing}. Using this data package, we approximate the originators of the P14/N14 and P22/N22 (14 and 22-millisecond post-stimulus, respectively) peaks using sLORETA/GMM reconstructions and lead field matrices obtained with different FE meshes. The first one of these peaks originated in the medial lemniscus of the brain stem above the cuneate nucleus \citep{noel1996origin, buchner1995somatotopy, mauguiere1983neural, urasaki1990origin, passmore2014origin}. The second one originated in two locations: a cortical originator has been observed in the crown of either the first pre-or post-central gyrus \citep{buchner1994source, allison1991cortical, fuchs1998improving}, and a sub-cortical one in the thalamus \cite{papadelis2011ba3b}.

\subsubsection{Analysis of FE mesh quality}
\label{sec:error_measurements}
We define mesh accuracy by measuring the distance between the surface of the grey matter compartment and the actual tissue boundary of the segmentation. We also measure the element condition distribution and provide a visual interpretation of the meshes described.

For the spherical model, we evaluate the lead field matrix accuracy at different eccentricities values, i.e., relative distances from the origin in the brain compartment using the relative difference (RDM) and magnitude (MAG) measures \cite{pursiainen2016electroencephalography}:
{\setlength\arraycolsep{2pt} \begin{eqnarray}
\hbox{RDM}(\vec{j}_1,\vec{j}_2)& = & \,\left\|\frac{\vec{j}_1}{\left\|\vec{j}_1\right\|_2} - \frac{\vec{j}_2}{\left\|\vec{j}_2\right\|_2}\right\|_2  \\
\hbox{MAG}(\vec{j}_1,\vec{j}_2)& = & 1-\,\left(\frac{\left\|\vec{j}_1\right\|_2}{\left\|\vec{j}_2\right\|_2}\right).
\end{eqnarray}}
Here, $\vec{j}_1$ and $\vec{j}_2$ denote FEM-based and analytical forward models, respectively. The analytical sources have random orientations but share the same set of positions as their FEM-based approximations. The RDM reflects the topographical forward modeling error respecting location and orientation, while MAG reflects the variations in potential amplitude. We calculate distributional source localization estimates for 4 different sets of 18 dipoles placed at eccentricities of 0.06, 0.29, 0.63, and 0.98, evaluating the EMD between the actual dipole source and the estimate.

\subsubsection{Computing platform}
The computation of the numerical results in this study used a Dell Precision 5820 Workstation with 256 GigaByte (GB) Random Access Memory (RAM), a 10-core Intel i9-10900X Central Processing Unit (CPU), and NVidia Quadro RTX 4000 Graphics Processing Unit (GPU). Both CPU and GPU were applied as two alternative solutions for parallelizing the labeling process. In the CPU parallelization, ten serial execution threads were run simultaneously (one per CPU core), while on the GPU, the process was decomposed into individually handled blocks of ten vectorized labeling operations.

\section{Results}
\label{sec:results}

\begin{figure}[h!]
\centering
    \begin{scriptsize}
        \begin{subfigure}[t]{6.0cm}
        \centering
            \includegraphics[width=5.8cm]{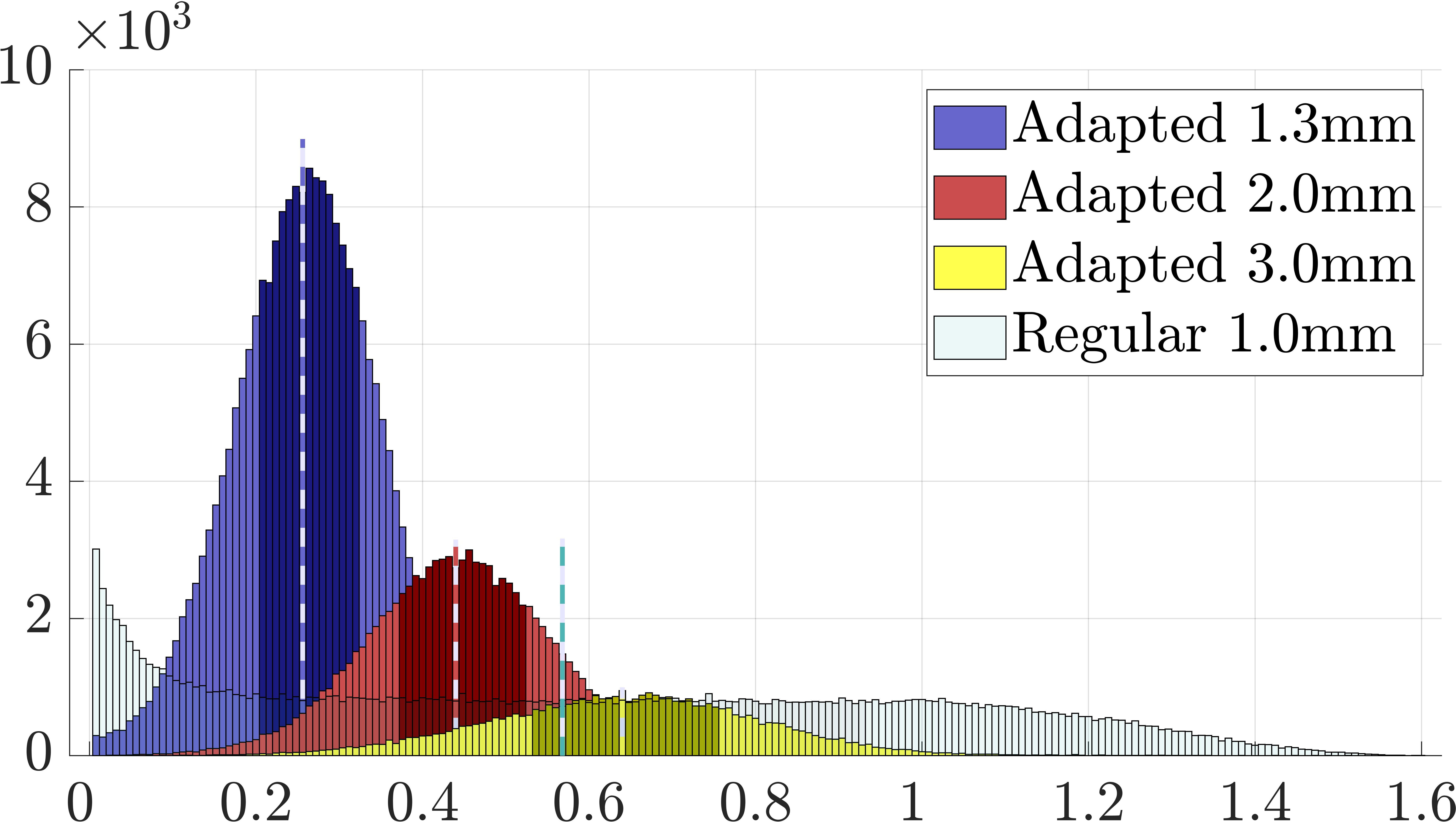}
            \caption{Spherical model}
            \label{fig:sphere_hist}
        \end{subfigure}
        \begin{subfigure}[t]{6.0cm}
        \centering
            \includegraphics[width=5.8cm]{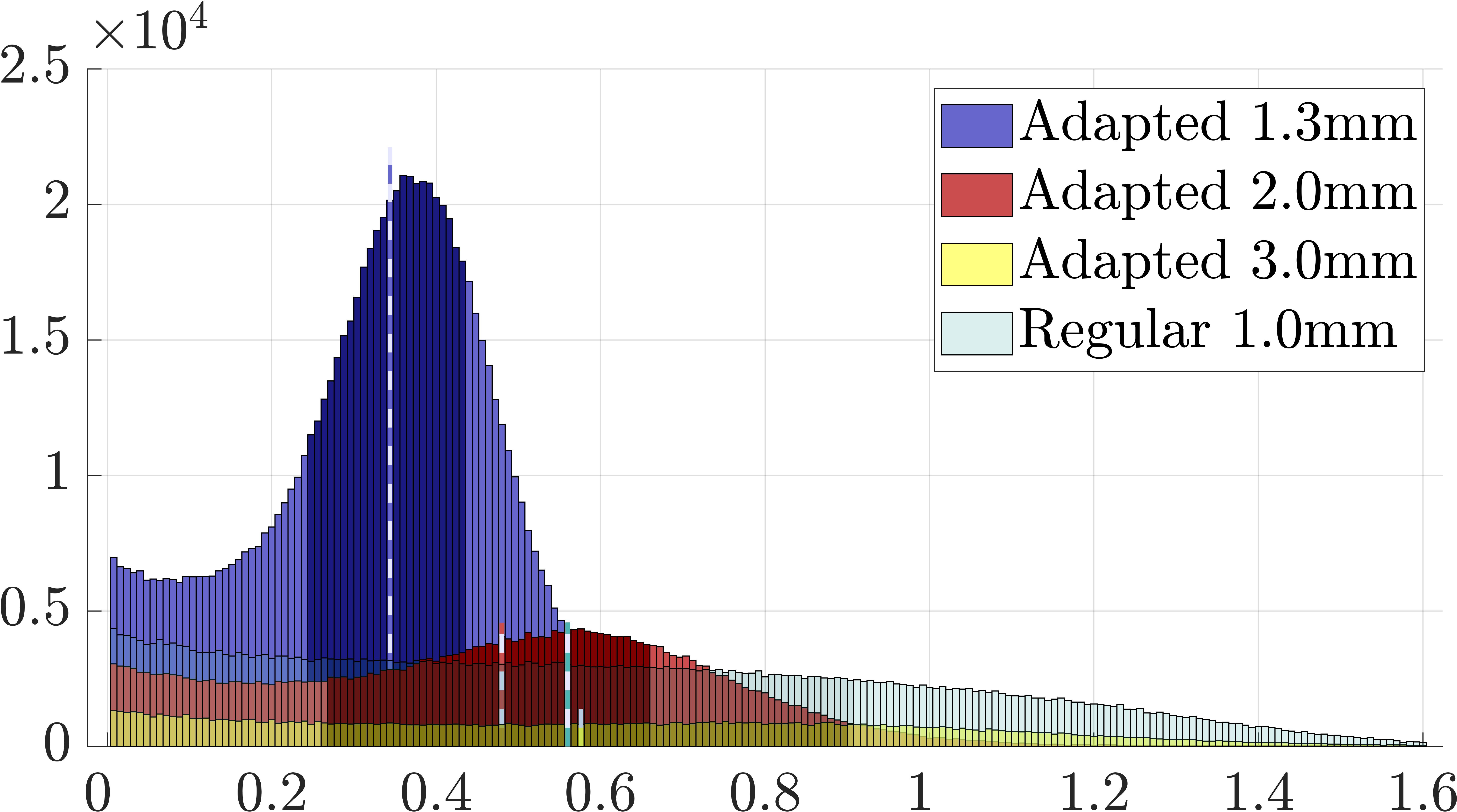}
            \caption{Realistic model}
            \label{fig:head_hist}
        \end{subfigure}
        \caption{Histograms for the spherical three-layer Ary model (\ref{fig:sphere_hist}) and realistic head model (\ref{fig:head_hist}) depicting the distance between the grey matter boundary in the FE mesh and in the segmentation.}
    \end{scriptsize}
\end{figure}

\begin{figure*}[h!]
    \centering
    \begin{scriptsize}
        \begin{subfigure}{4.2cm}
            \centering
            \includegraphics[width=4.1cm]{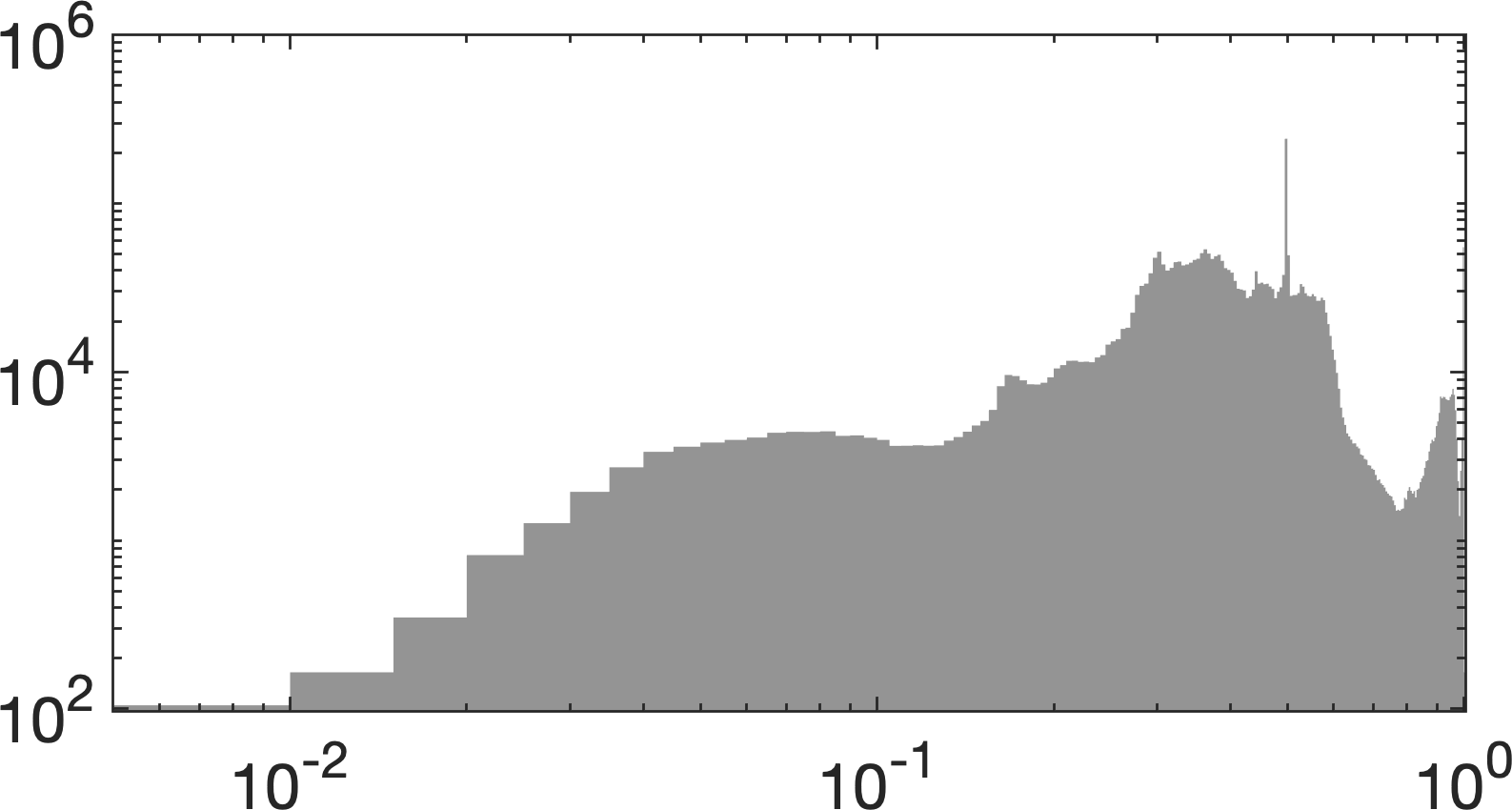}
            \caption{Sphere 3.0 mm}
            \label{fig:ele_cond_histogram_a}
        \end{subfigure}
        \begin{subfigure}{4.2cm}
            \centering
            \includegraphics[width=4.1cm]{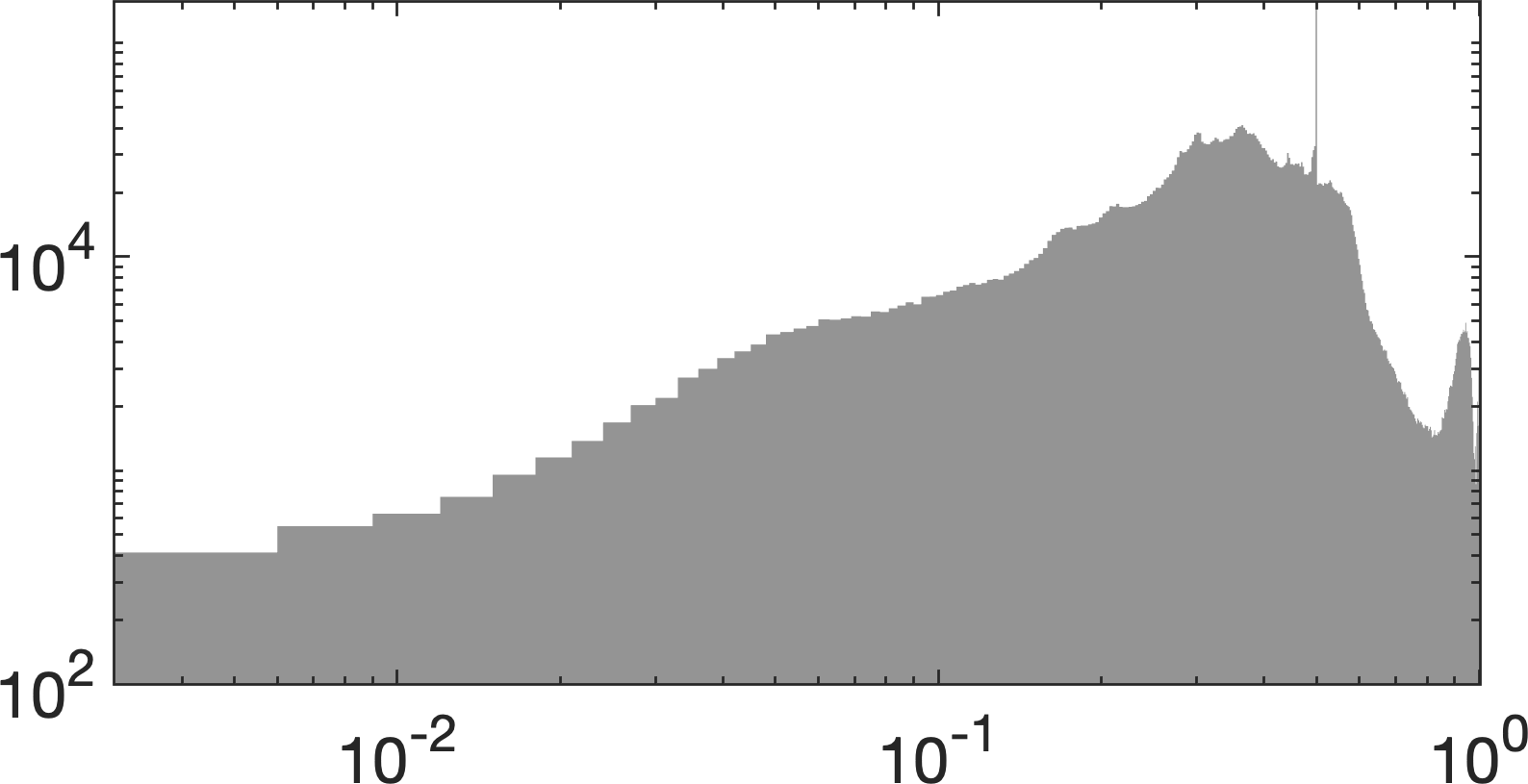}
            \caption{Realistic 3.0 mm}
            \label{fig:ele_cond_histogram_b}
        \end{subfigure}
        \begin{subfigure}{3cm}
            \centering
            \includegraphics[width=2.6cm]{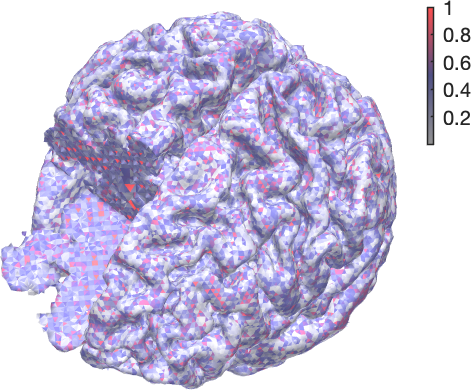}
            \caption{Realistic 3.0 mm}
            \label{fig:ele_cond1}
        \end{subfigure} 
        \vskip0.2cm
        \begin{subfigure}{4.2cm}
            \centering
            \includegraphics[width=4.1cm]{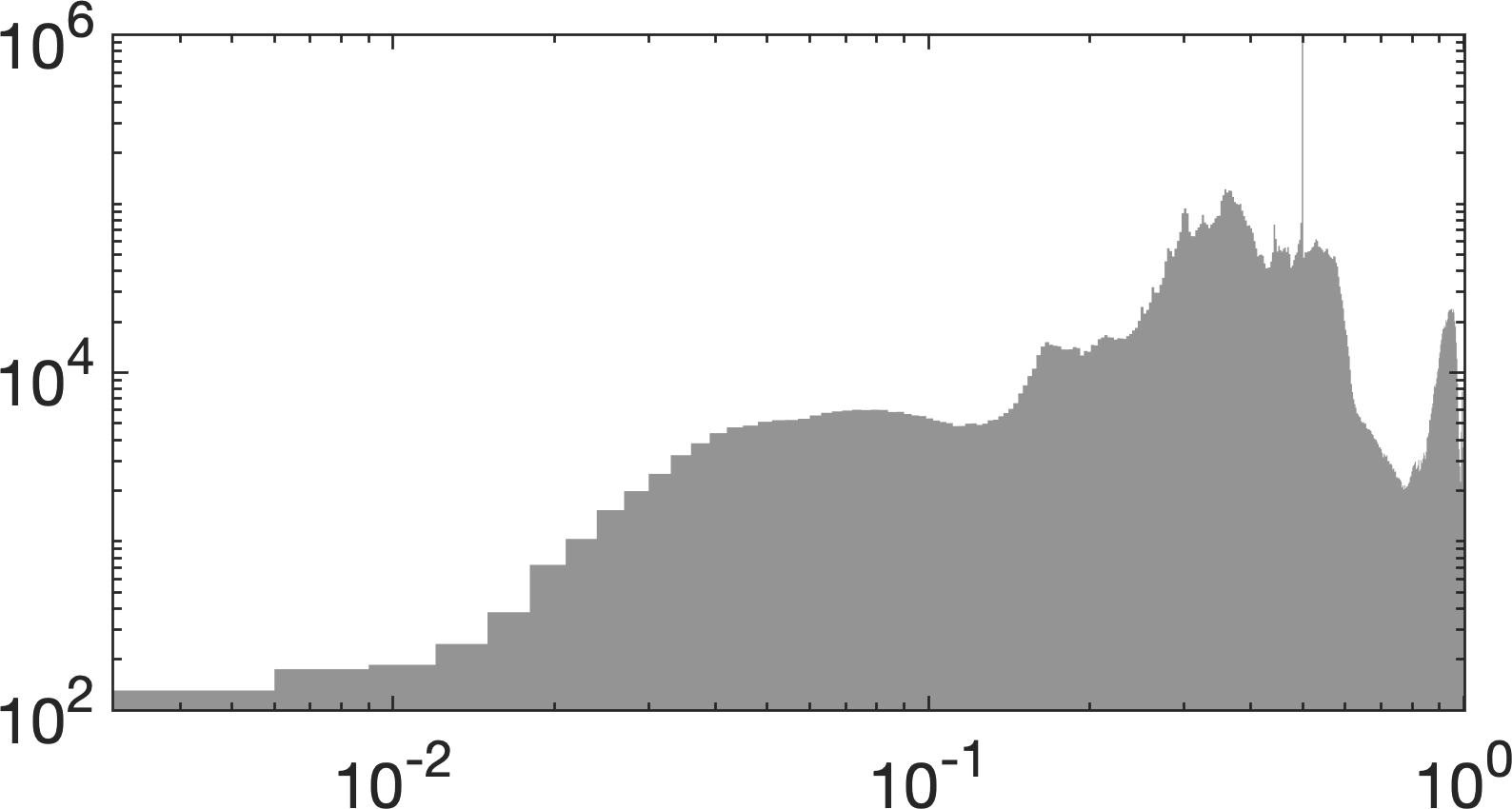}
            \caption{Sphere 2.0mm}
            \label{fig:ele_cond_histogram_c}
        \end{subfigure}
        \begin{subfigure}{4.2cm}
            \centering
            \includegraphics[width=4.1cm]{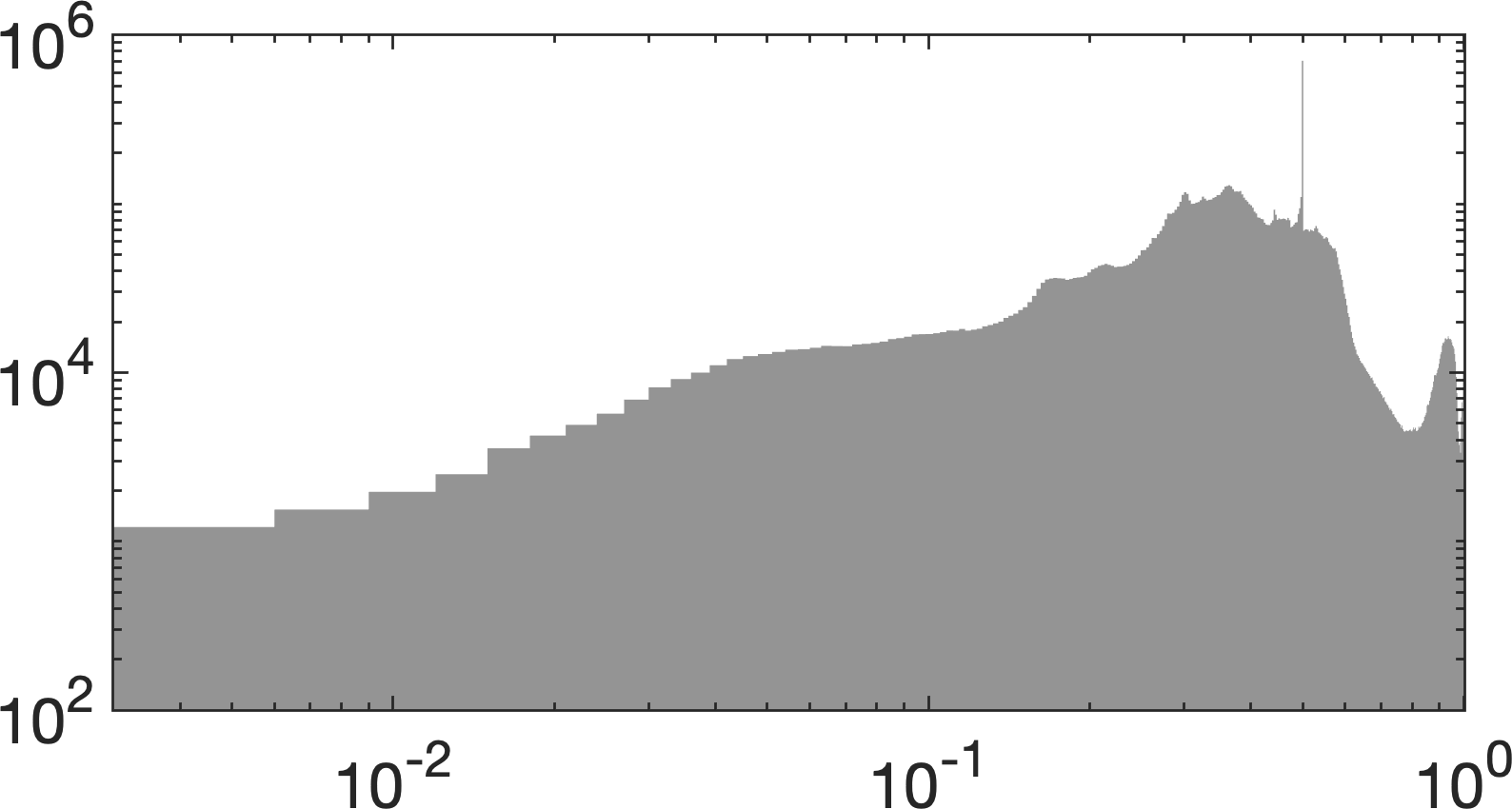}
            \caption{Realistic 2.0mm}
            \label{fig:ele_cond_histogram_d}
        \end{subfigure}
        \begin{subfigure}{3cm}
            \centering
            \includegraphics[width=2.6cm]{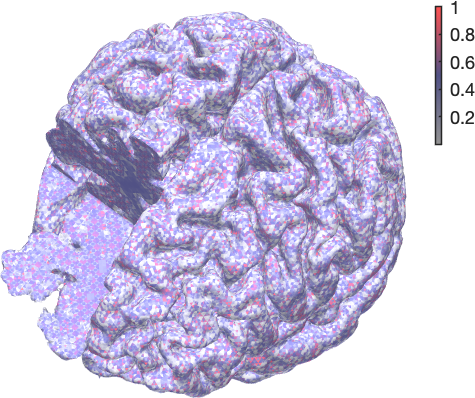}
            \caption{Realistic 2.0 mm}
            \label{fig:ele_cond2}
        \end{subfigure} 
        \vskip0.2cm
        \begin{subfigure}{4.2cm}
            \centering
            \includegraphics[width=4.1cm]{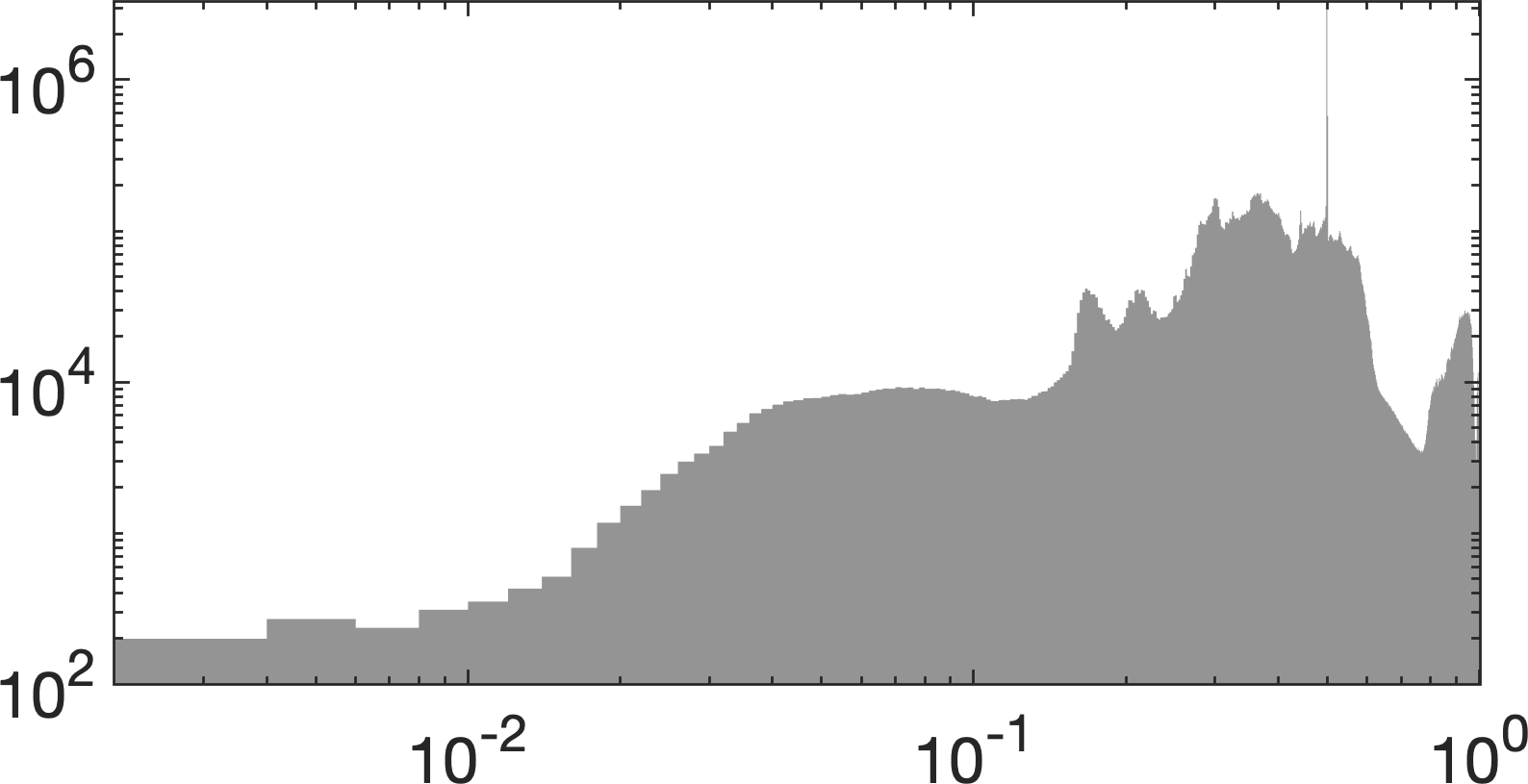}
            \caption{Sphere 1.3mm}
            \label{fig:ele_cond_histogram_e}
        \end{subfigure}
        \begin{subfigure}{4.2cm}
            \centering
            \includegraphics[width=4.1cm]{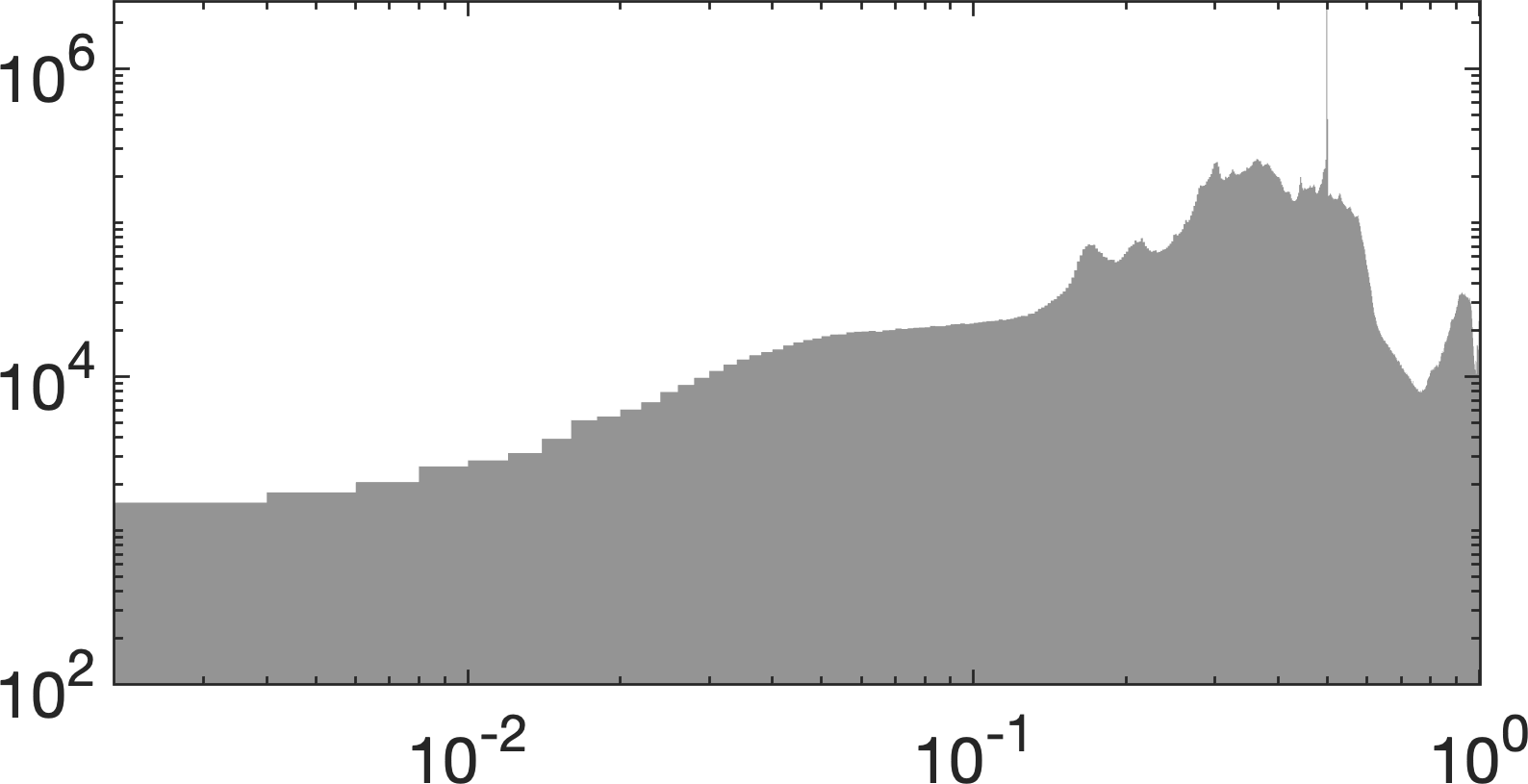}
            \caption{Realistic 1.3mm}
            \label{fig:ele_cond_histogram_f}
        \end{subfigure}
        \begin{subfigure}{3cm}
            \centering
            \includegraphics[width=2.6cm]{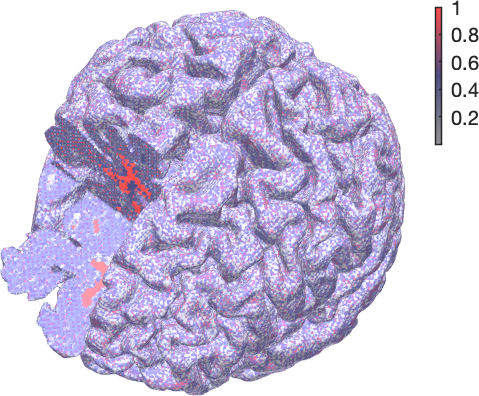}
            \caption{Realistic 1.3 mm}
            \label{fig:ele_cond3}
        \end{subfigure} 
    \end{scriptsize}
\caption{Element condition distributions of the adapted 3.0, 2.0 and 1.3 mm  FE mesh shown as a histogram for  the spherical three-layer Ary model (\ref{fig:ele_cond_histogram_a}, \ref{fig:ele_cond_histogram_c}, \ref{fig:ele_cond_histogram_e}) and realistic head model (\ref{fig:ele_cond_histogram_b}, \ref{fig:ele_cond_histogram_d}, \ref{fig:ele_cond_histogram_f}). The spatial condition number mappings (\ref{fig:ele_cond1}, \ref{fig:ele_cond2}, \ref{fig:ele_cond3}) for the realistic head model show as expected that the non-refined internal parts of the mesh, in particular, the interior of the white matter, correspond to an elevated condition compared to the thin refined layers closer to the surface.}
\label{fig:element_condition}
\end{figure*}

\begin{figure}[h!]
    \centering
    \begin{scriptsize}
        \begin{minipage}[T]{7.5cm}
            \centering

            \begin{minipage}[T]{1.60cm}
                \centering
                \textbf{Adapted \\ 3.0mm}
            \end{minipage}
            \begin{minipage}[T]{1.60cm}
                \centering
                \textbf{Adapted \\ 2.0mm}
            \end{minipage}
            \begin{minipage}[T]{1.60cm}
                \centering
                \textbf{Adapted \\ 1.3mm}
            \end{minipage}
            \begin{minipage}[T]{1.60cm}
                \centering
                \textbf{Regular \\ 1.0mm}
            \end{minipage}
            
        \end{minipage}
        
        \vskip0.1cm
        
        \begin{minipage}[T]{7.5cm}
        \centering
            \begin{subfigure}[T]{1.60cm}
                \centering
                \includegraphics[width=1.55cm]{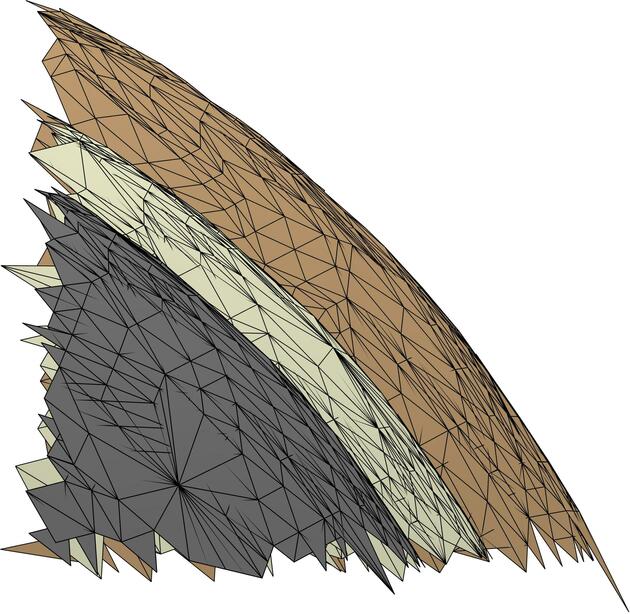}
                \caption{}
                \label{fig:sphere_comparison_a}
            \end{subfigure}
            \begin{subfigure}[T]{1.60cm}
                \centering
                \includegraphics[width=1.55cm]{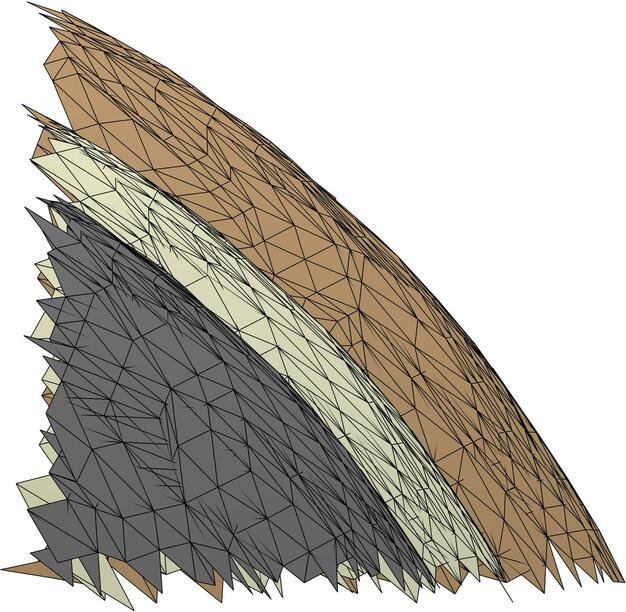}
                \caption{}
                \label{fig:sphere_comparison_b}
            \end{subfigure}
            \begin{subfigure}[T]{1.60cm}
                \centering
                \includegraphics[width=1.55cm]{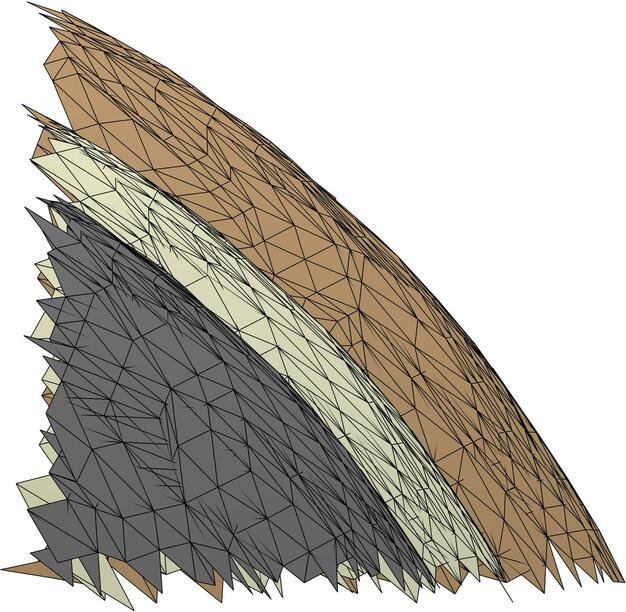}
                \caption{}
                \label{fig:sphere_comparison_c}
            \end{subfigure}
            \begin{subfigure}[T]{1.60cm}
                \centering
                \includegraphics[width=1.55cm]{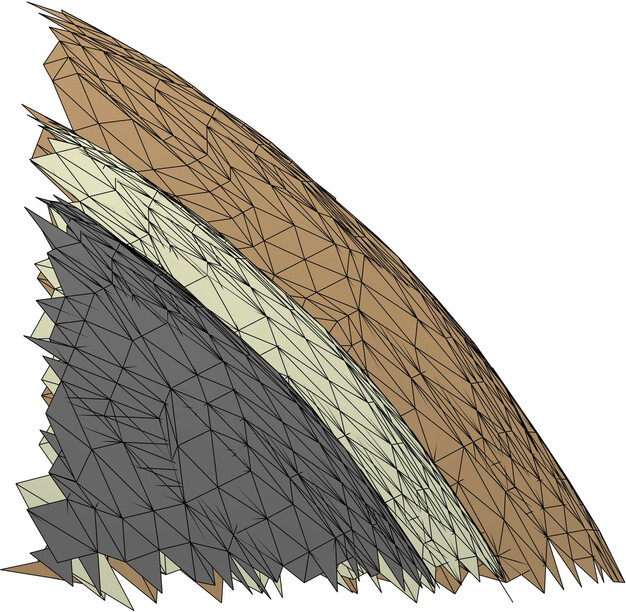}
                \caption{}
                \label{fig:sphere_comparison_d}
            \end{subfigure}
        
        \vskip0.1cm
        
            \begin{subfigure}[T]{1.60cm}
                \centering
                \includegraphics[width=1.55cm]{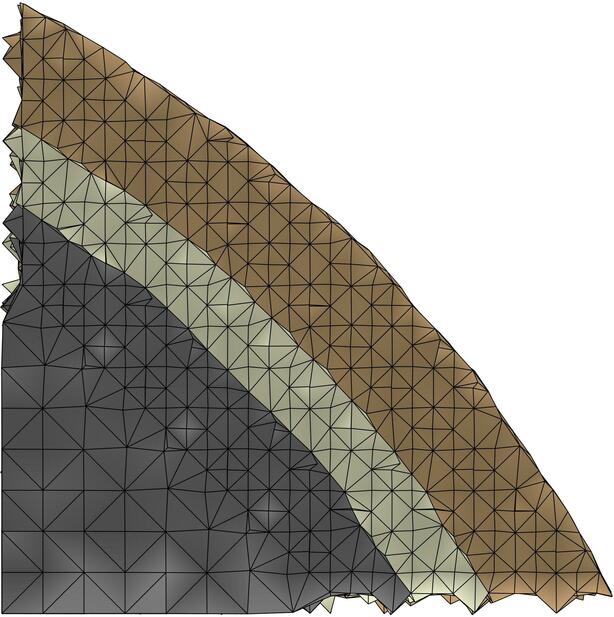}
                \caption{}
                \label{fig:sphere_comparison_e}
            \end{subfigure}
            \begin{subfigure}[T]{1.60cm}
                \centering
                \includegraphics[width=1.55cm]{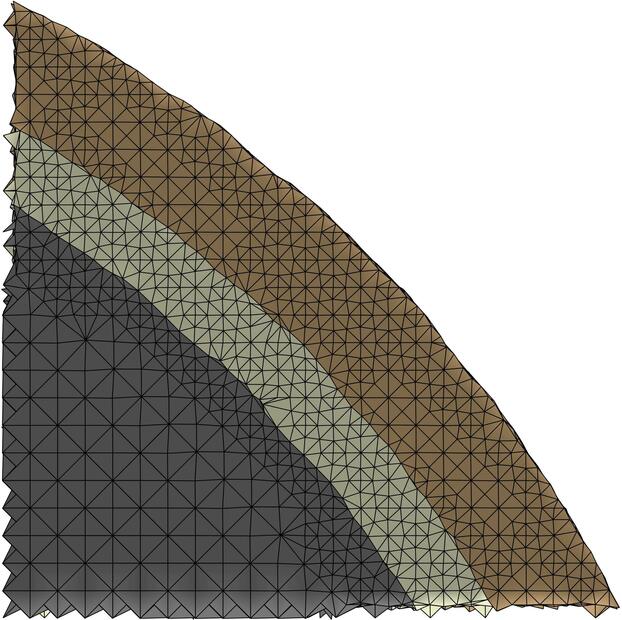}
                \caption{}
                \label{fig:sphere_comparison_f}
            \end{subfigure}
            \begin{subfigure}[T]{1.60cm}
                \centering
                \includegraphics[width=1.55cm]{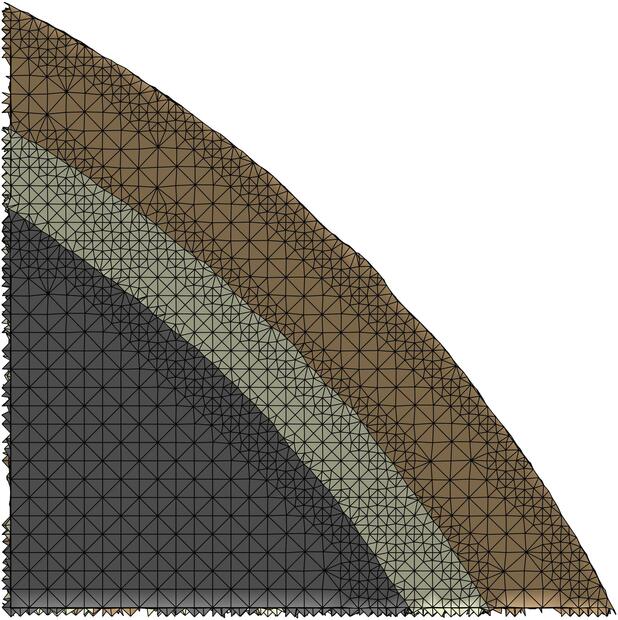}
                \caption{}
                \label{fig:sphere_comparison_g}
            \end{subfigure}
            \begin{subfigure}[T]{1.60cm}
                \centering
                \includegraphics[width=1.55cm]{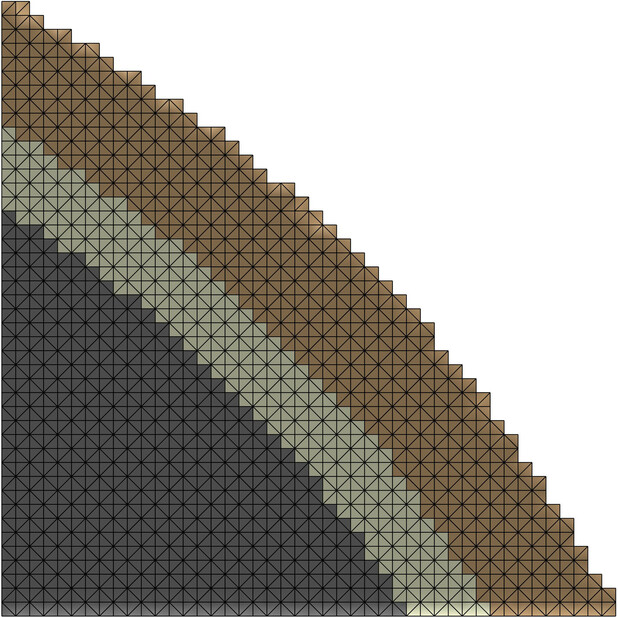}
                \caption{}
                \label{fig:sphere_comparison_h}
            \end{subfigure}
        \end{minipage}
    
\caption{Quadrants of the adapted surface- (\ref{fig:sphere_comparison_a}-\ref{fig:sphere_comparison_c}) and volume- (\ref{fig:sphere_comparison_e}-\ref{fig:sphere_comparison_g}) based reconstructions of the spherical three-layer Ary model depicting proposed mesh generation and post-processing methods using 3.0, 2.0 and 1.3 millimeter mesh sizes, respectively. Regular surface- (\ref{fig:sphere_comparison_d}) and volume- \ref{fig:sphere_comparison_h}) with 1.0 mm mesh size without post-processing effect are included for comparison. The presented layers are (top-bot): scalp (brown), skull (white) and gray matter.
\label{fig:adapted_slice_sphere}}
\end{scriptsize}
\end{figure}

\begin{figure}[h!]
\centering
    \begin{scriptsize}
        \begin{minipage}[T]{7.5cm}
        \centering
            \begin{minipage}[T]{1.60cm}
                \centering
                \textbf{Adapted \\ 3.0mm}
            \end{minipage}
            \begin{minipage}[T]{1.60cm}
                \centering
                \textbf{Adapted \\ 2.0mm}
            \end{minipage}
            \begin{minipage}[T]{1.60cm}
                \centering
                \textbf{Adapted \\ 1.3mm}
            \end{minipage}
            \begin{minipage}[T]{1.60cm}
                \centering
                \textbf{Regular \\ 1.0mm}
            \end{minipage}
        \end{minipage}
        
        \vskip0.1cm
        
        \begin{minipage}[T]{7.5cm}
        \centering
            \begin{subfigure}[T]{1.60cm}
                \centering
                \includegraphics[width=1.55cm]{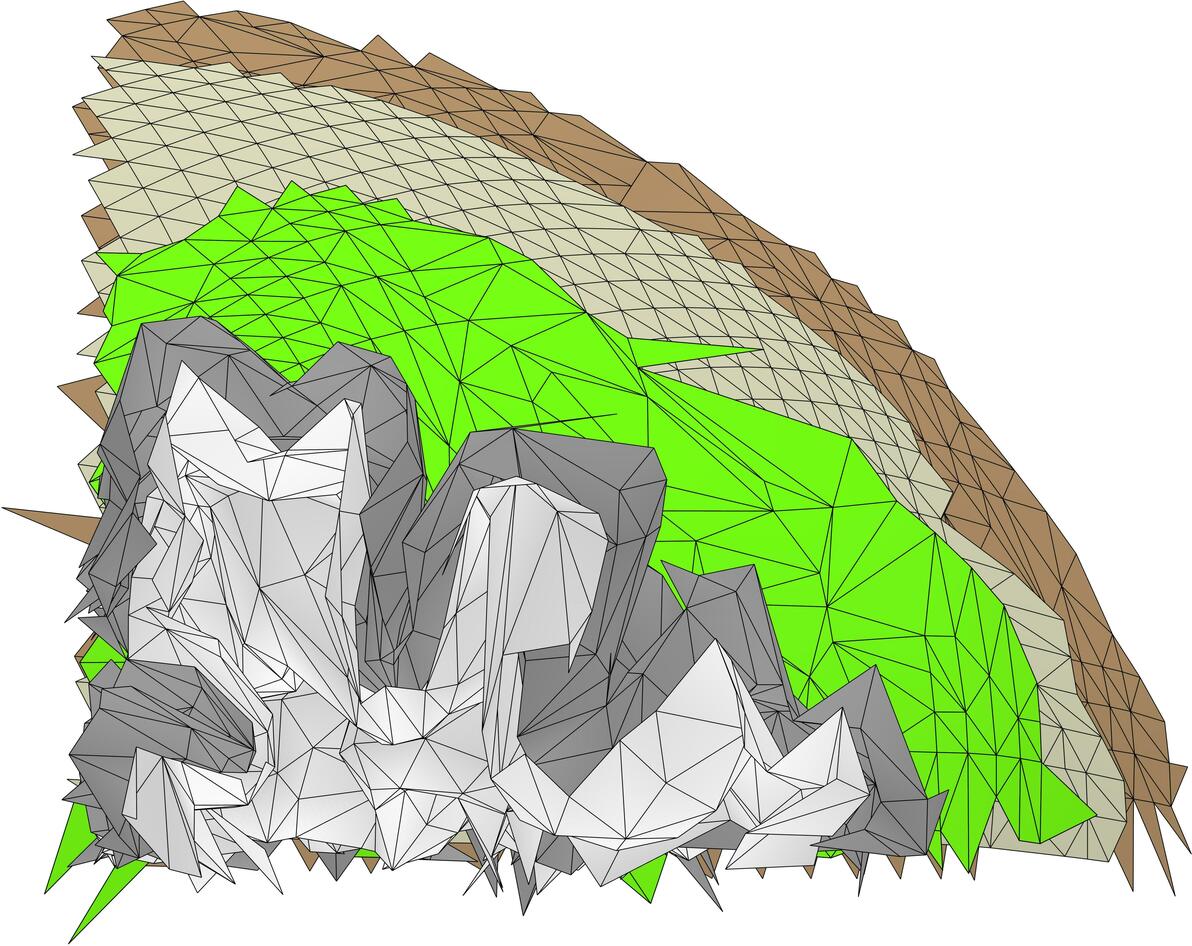}
                \caption{}
                \label{fig:carsten_comparision_a}
            \end{subfigure}
            \begin{subfigure}[T]{1.60cm}
                \centering
                \includegraphics[width=1.55cm]{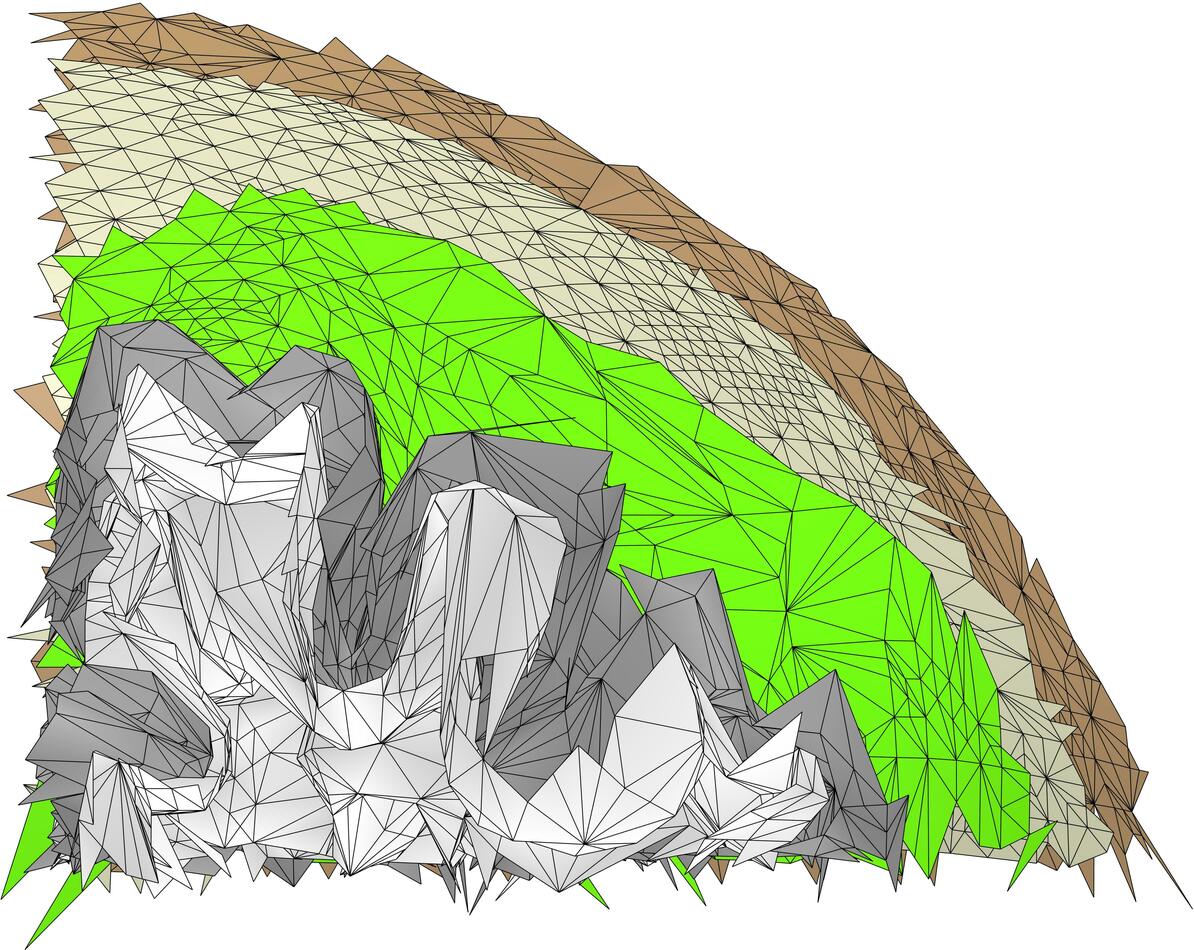}
                \caption{}
                \label{fig:carsten_comparision_b}
            \end{subfigure}
            \begin{subfigure}[T]{1.60cm}
                \centering
                \includegraphics[width=1.55cm]{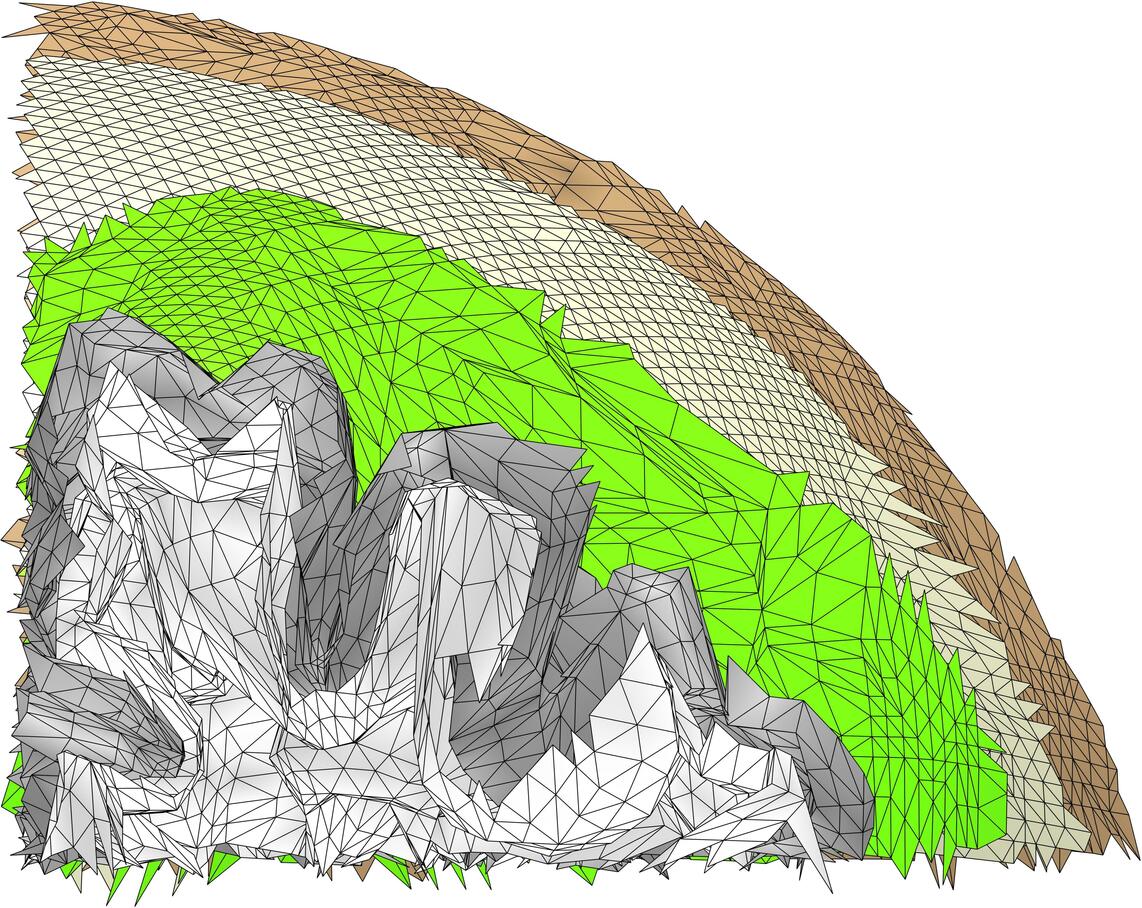}
                \caption{}
                \label{fig:carsten_comparision_c}
            \end{subfigure}
            \begin{subfigure}[T]{1.60cm}
                \centering
                \includegraphics[width=1.55cm]{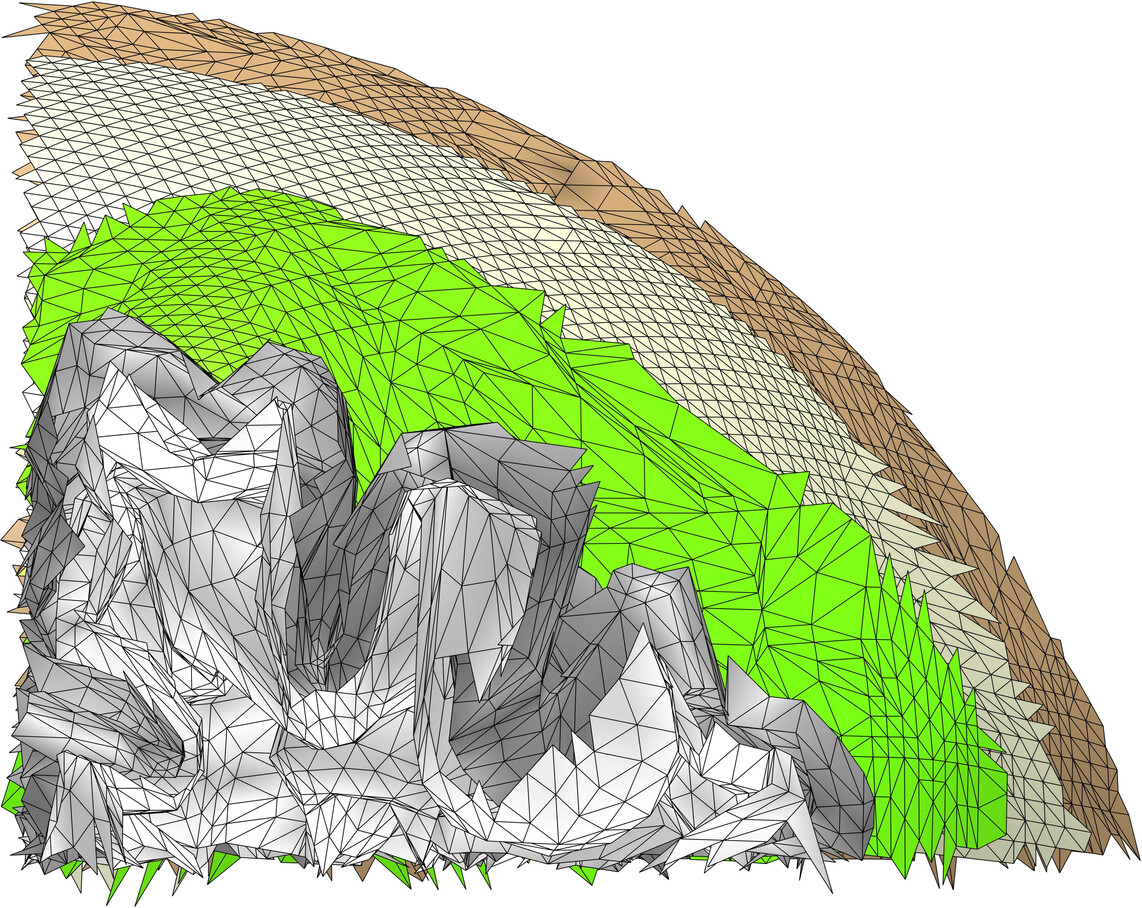}
                \caption{}
                \label{fig:carsten_comparision_d}
            \end{subfigure}
        \end{minipage}
        
        \vskip0.1cm
        
        \begin{minipage}[T]{7.5cm}
        \centering
            \begin{subfigure}[T]{1.60cm}
            \centering
                \includegraphics[width=1.55cm]{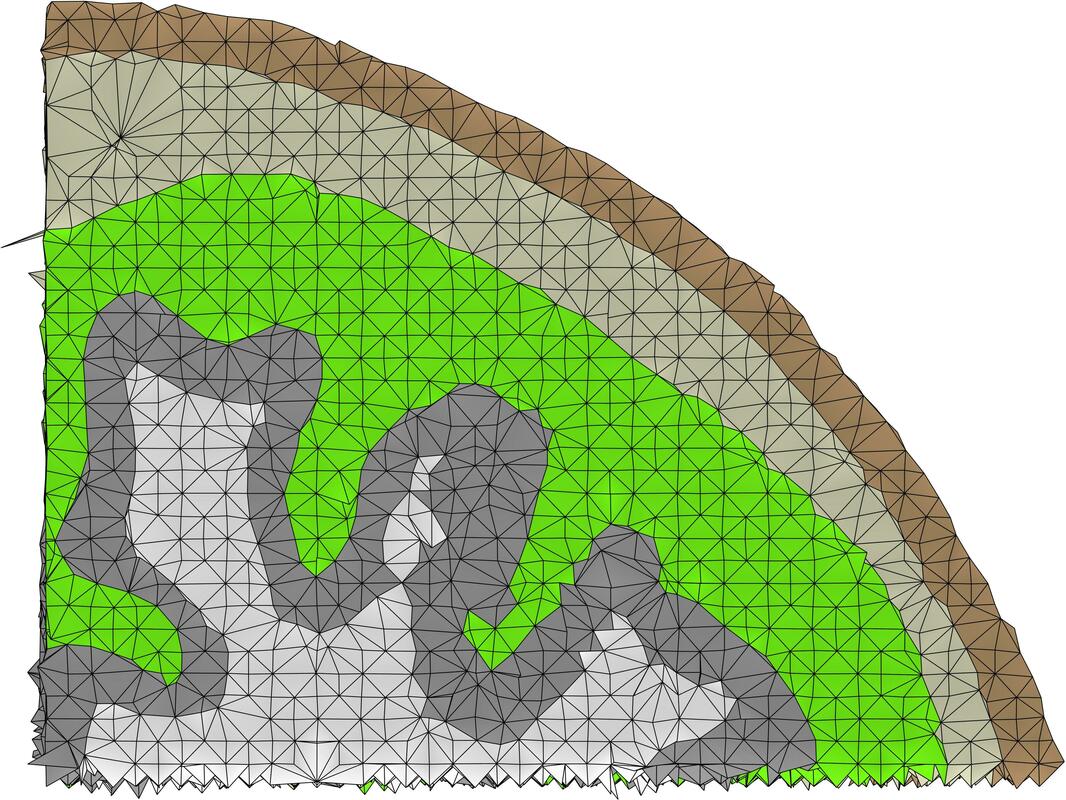}
                \caption{}
                \label{fig:carsten_comparision_e}
            \end{subfigure}
            \begin{subfigure}[T]{1.60cm}
                \centering
                \includegraphics[width=1.55cm]{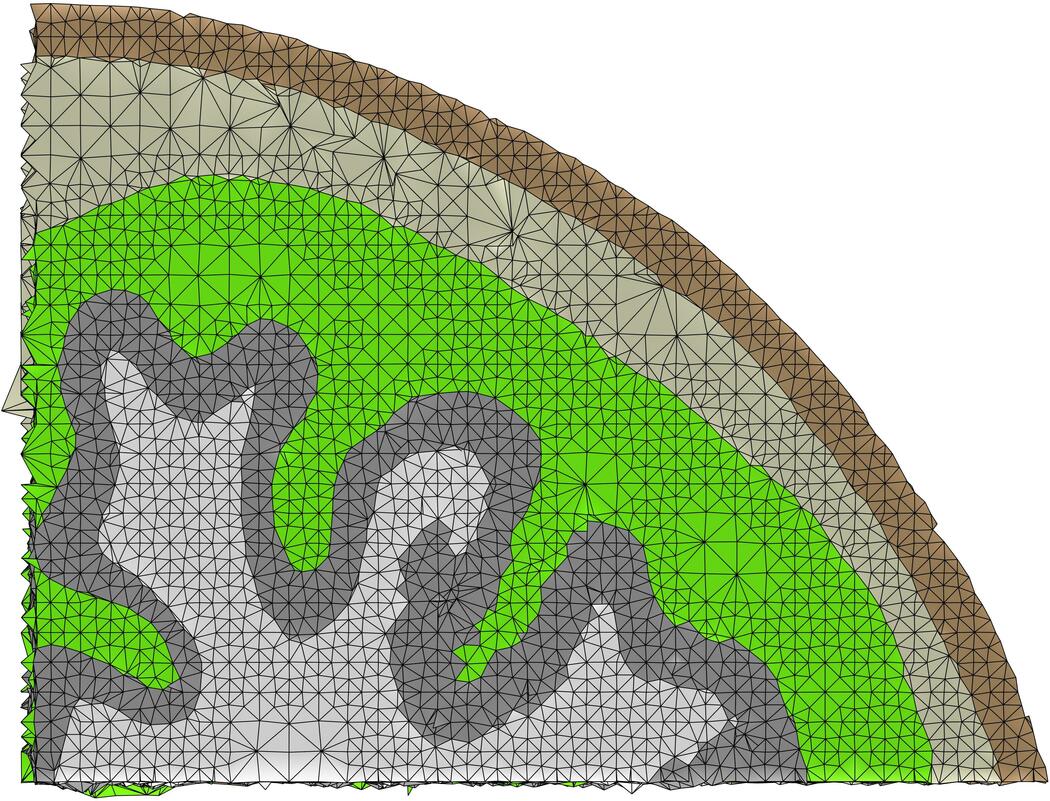}
                \caption{}
                \label{fig:carsten_comparision_f}
            \end{subfigure}
            \begin{subfigure}[T]{1.60cm}
                \centering
                \includegraphics[width=1.55cm]{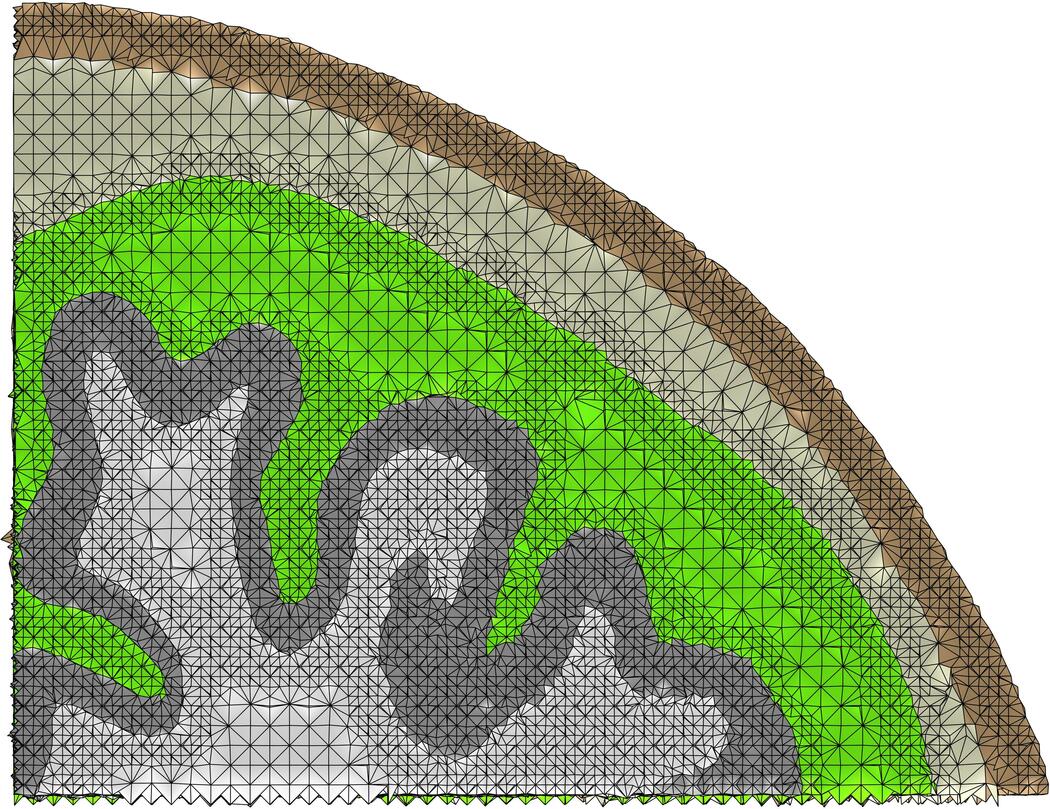}
                \caption{}
                \label{fig:carsten_comparision_g}
            \end{subfigure}
            \begin{subfigure}[T]{1.60cm}
                \centering
                \includegraphics[width=1.55cm]{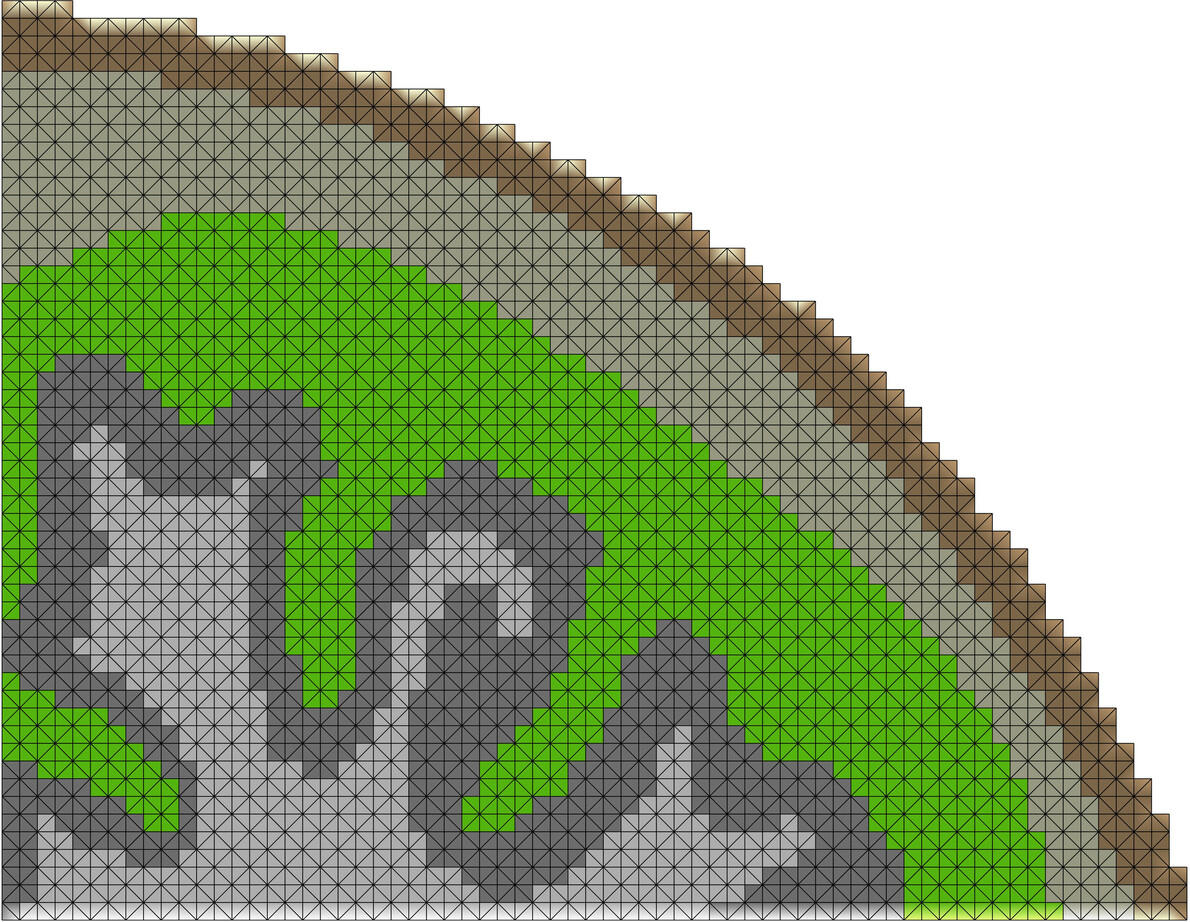}
                \caption{}
                \label{fig:carsten_comparision_h}
            \end{subfigure}
        \end{minipage}
    
\caption{Quadrants of the adapted surface- (\ref{fig:carsten_comparision_a}-\ref{fig:carsten_comparision_c}) and volume- (\ref{fig:carsten_comparision_e}-\ref{fig:carsten_comparision_g}) based reconstructions of the realistic head model depicting proposed mesh generation and post-processing methods using 3.0, 2.0 and 1.3 millimeter mesh sizes, respectively. Regular surface- (\ref{fig:carsten_comparision_d}) and volume- \ref{fig:carsten_comparision_h}) with 1.0 mm mesh size without post-processing effect are included for comparison. The presented layers (top-bot) are: scalp (dark brown), skull (light brown), cerebrospinal fluid (green), gray matter, and white matter.
\label{fig:adapted_slice_realistic}}
\end{scriptsize}
\end{figure}%

\begin{figure}[h!]
\centering
    \begin{scriptsize}
        \begin{minipage}[T]{7.5cm}
        \centering
            \begin{minipage}[T]{1.60cm}
                \centering
                \textbf{Adapted \\ 3.0mm}
            \end{minipage}
            \begin{minipage}[T]{1.60cm}
                \centering
                \textbf{Adapted \\ 2.0mm}
            \end{minipage}
            \begin{minipage}[T]{1.60cm}
                \centering
                \textbf{Adapted \\ 1.3mm}
            \end{minipage}
            \begin{minipage}[T]{1.60cm}
                \centering
                \textbf{Regular \\ 1.0mm}
            \end{minipage}
        \end{minipage}
        
        \vskip0.1cm
        
        \begin{minipage}[T]{7.5cm}
            \centering
            \begin{subfigure}[T]{1.60cm}
                \centering
                \includegraphics[width=1.55cm]{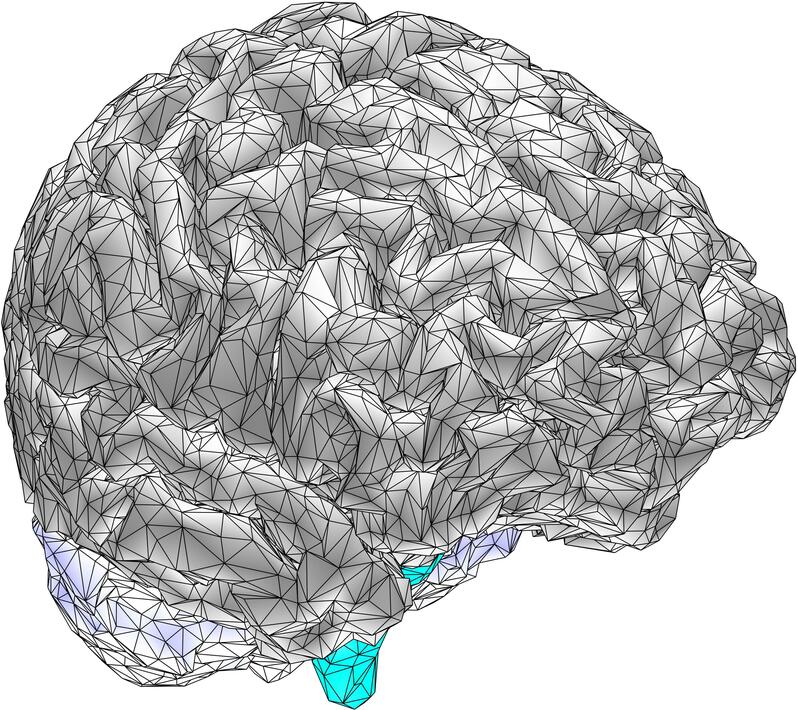}
            \end{subfigure}
            \begin{subfigure}[T]{1.60cm}
                \centering
                \includegraphics[width=1.55cm]{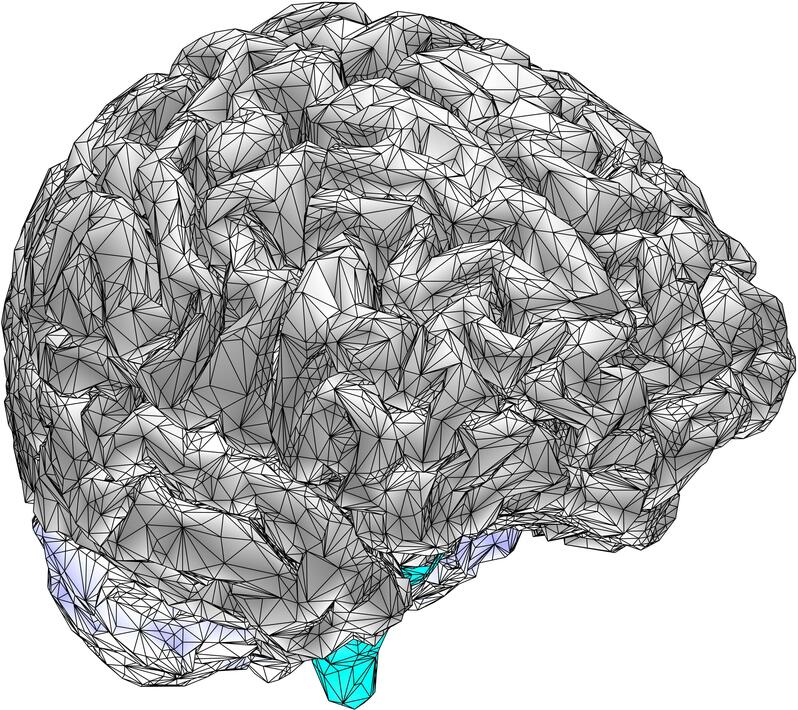}
            \end{subfigure}
            \begin{subfigure}[T]{1.60cm}
                \centering
                \includegraphics[width=1.55cm]{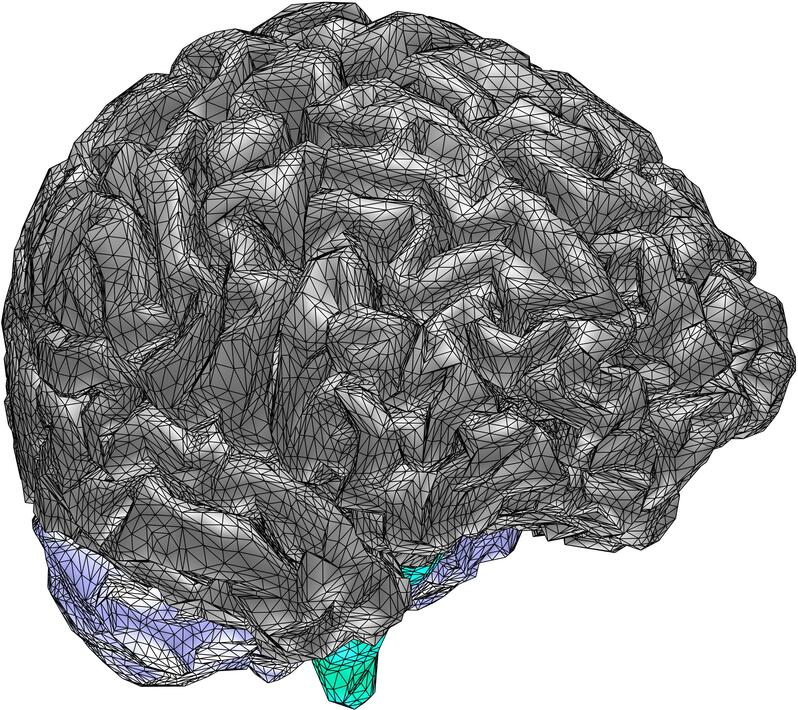}
            \end{subfigure}
            \begin{subfigure}[T]{1.60cm}
                \centering
                \includegraphics[width=1.55cm]{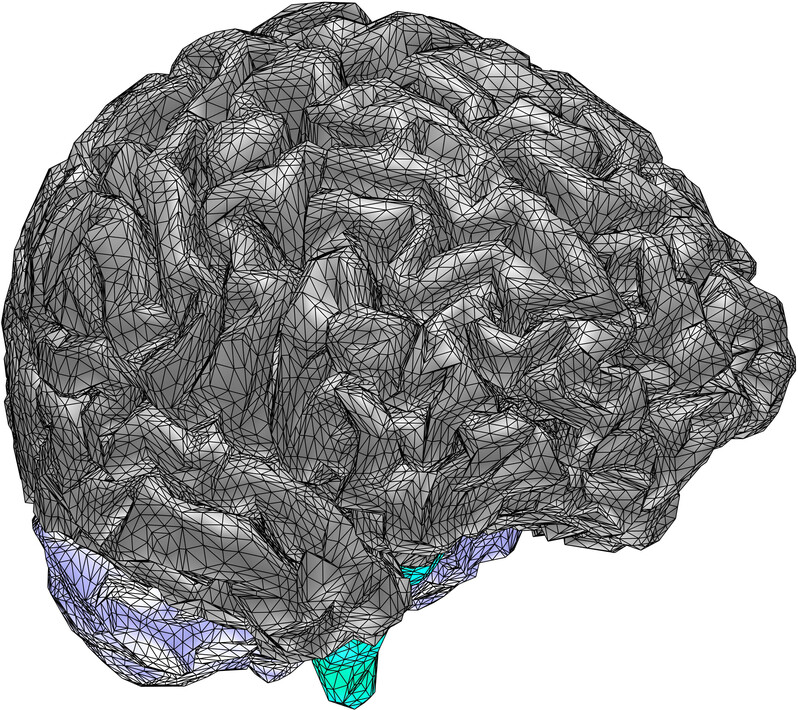}
            \end{subfigure}
            
            \vskip0.1cm
            
            \begin{subfigure}[T]{1.60cm}
                \centering
                \includegraphics[width=1.55cm]{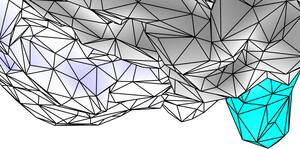}
                \caption{}
                \label{fig:cortex_a}
            \end{subfigure}
            \begin{subfigure}[T]{1.60cm}
                \centering
                \includegraphics[width=1.55cm]{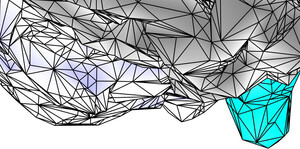}
                \caption{}
                \label{fig:cortex_b}
            \end{subfigure}
            \begin{subfigure}[T]{1.60cm}
                \centering
                \includegraphics[width=1.55cm]{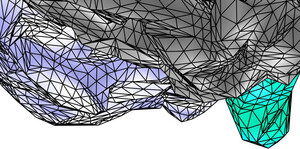}
                \caption{}
                \label{fig:cortex_c}
            \end{subfigure}
            \begin{subfigure}[T]{1.60cm}
                \centering
                \includegraphics[width=1.55cm]{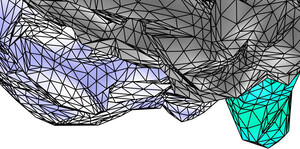}
                \caption{}
                \label{fig:cortex_d}
            \end{subfigure}
        \end{minipage}

        \hrule
        \vskip0.2cm
        
        \begin{minipage}[T]{7.5cm}
        \centering
            \begin{subfigure}[T]{1.60cm}
                \centering
                \includegraphics[width=1.55cm]{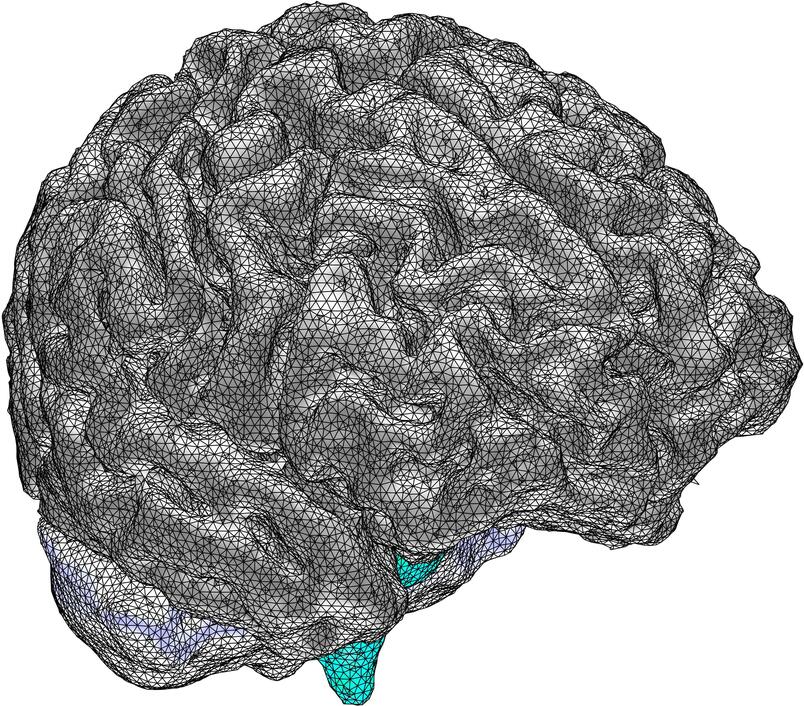}
            \end{subfigure}
            \begin{subfigure}[T]{1.60cm}
                \centering
                \includegraphics[width=1.55cm]{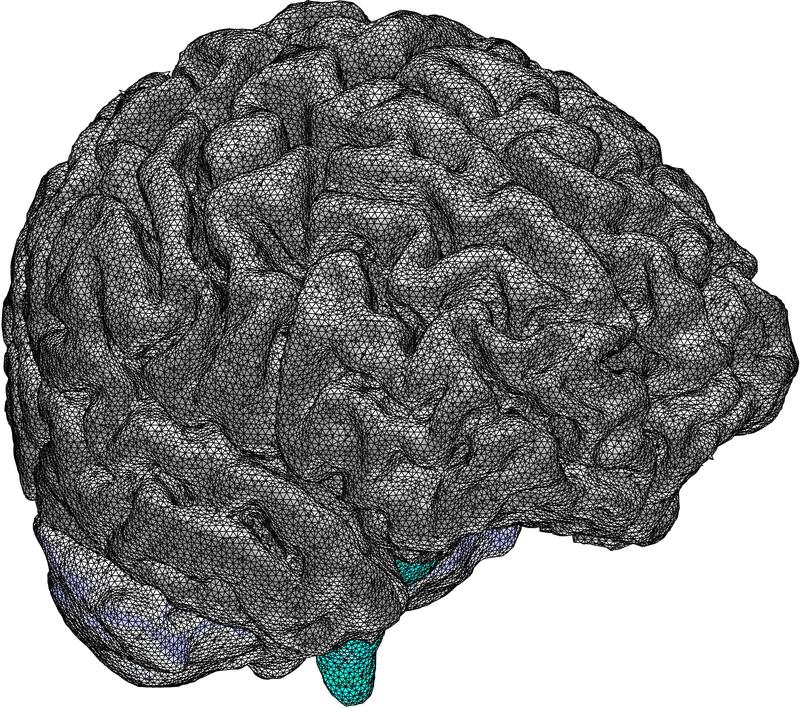}
            \end{subfigure}
            \begin{subfigure}[T]{1.60cm}
                \centering
                \includegraphics[width=1.55cm]{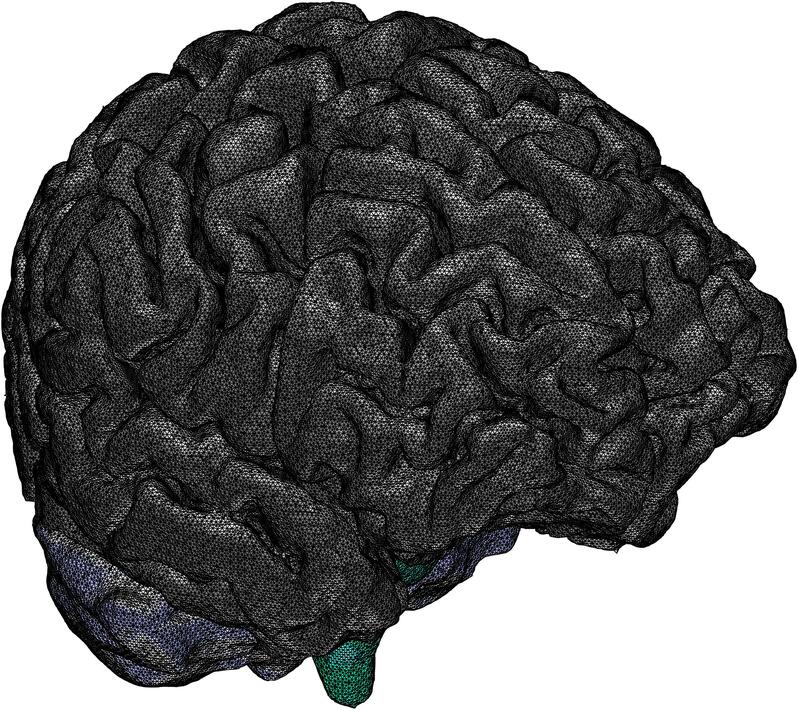}
            \end{subfigure}
            \begin{subfigure}[T]{1.60cm}
                \centering
                \includegraphics[width=1.55cm]{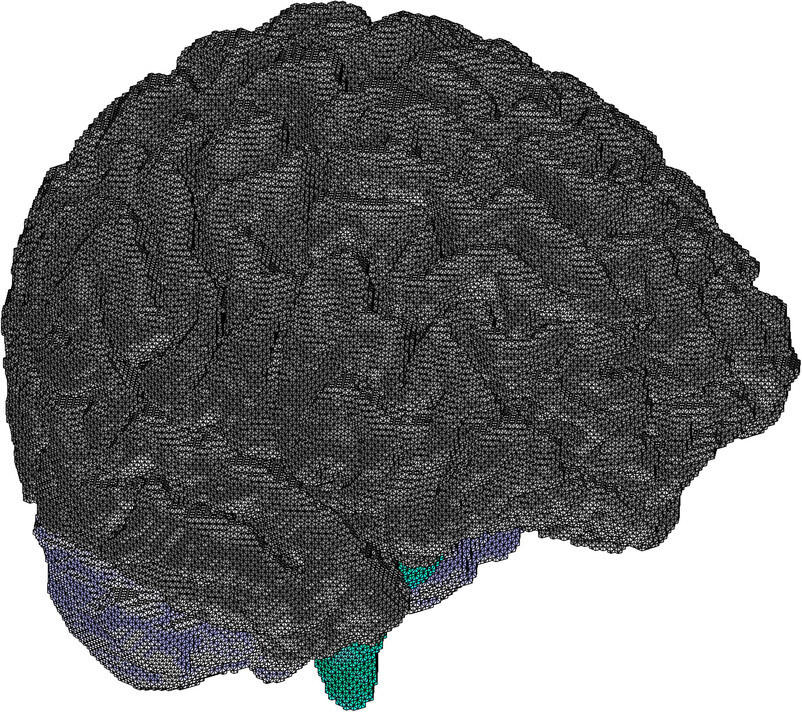}
            \end{subfigure}
            \\ \vskip0.1cm
            \begin{subfigure}[T]{1.60cm}
                \centering
                \includegraphics[width=1.55cm]{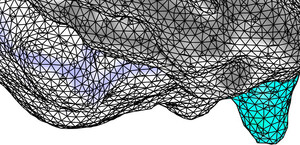}
                \caption{}
                \label{fig:cortex_e}
            \end{subfigure}
            \begin{subfigure}[T]{1.60cm}
                \centering
                \includegraphics[width=1.55cm]{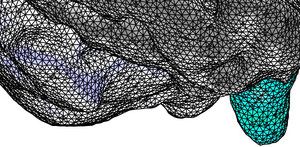}
                \caption{}
                \label{fig:cortex_f}
            \end{subfigure}
            \begin{subfigure}[T]{1.60cm}
                \centering
                \includegraphics[width=1.55cm]{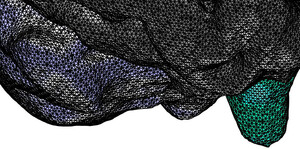}
                \caption{}
                \label{fig:cortex_g}
            \end{subfigure}
            \begin{subfigure}[T]{1.60cm}
                \centering
                \includegraphics[width=1.55cm]{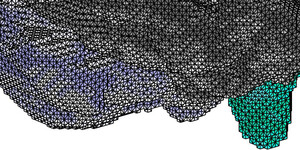}
                \caption{}
                \label{fig:cortex_h}
            \end{subfigure}
        \end{minipage} 
\caption{Reconstruction of a realistic head model depicting the performance of adaptive FE mesh generation using 3.0, 2.0, and 1.3 mm mesh sizes (Surface: (\ref{fig:cortex_a}-\ref{fig:cortex_c}); Volume: (\ref{fig:cortex_e}-\ref{fig:cortex_g}), respectively) after applying recursive solid angle labeling, smoothing, inflation and optimization procedures in contrast with a regular 1.0 mm mesh generation without post-process procedures (\ref{fig:cortex_d} and \ref{fig:cortex_h}). Regular mesh with a 1.0 mm mesh size without a post-processing effect is included for comparison. Cerebrum (gray), cerebellum (light blue) and brain stem (cyan).}
\label{fig:carsten_surface}
\end{scriptsize}
\end{figure}

\begin{figure}[h!]
\centering
    \begin{scriptsize}
        \begin{minipage}[T]{7.5cm}
        \centering
            \begin{minipage}[T]{1.60cm}
                \centering
                \textbf{Adapted 3.0mm}
            \end{minipage}
            \begin{minipage}[T]{1.60cm}
                \centering
                \textbf{Adapted 2.0mm}
            \end{minipage}
            \begin{minipage}[T]{1.60cm}
                \centering
                \textbf{Adapted 1.3mm}
            \end{minipage}
            \begin{minipage}[T]{1.60cm}
                \centering
                \textbf{Regular 1.0mm}
            \end{minipage}
        \end{minipage}
        
        \vskip0.1cm
        
        \begin{minipage}[T]{7.5cm}
        \centering
            \begin{subfigure}[T]{1.60cm}
                \centering
                \includegraphics[width=1.55cm]{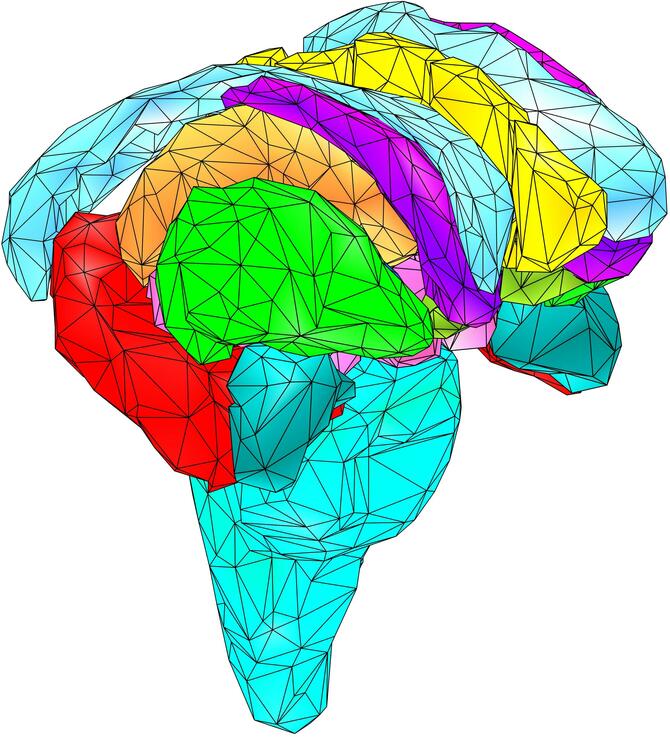}
            \end{subfigure}
            \begin{subfigure}[T]{1.60cm}
                \centering
                \includegraphics[width=1.55cm]{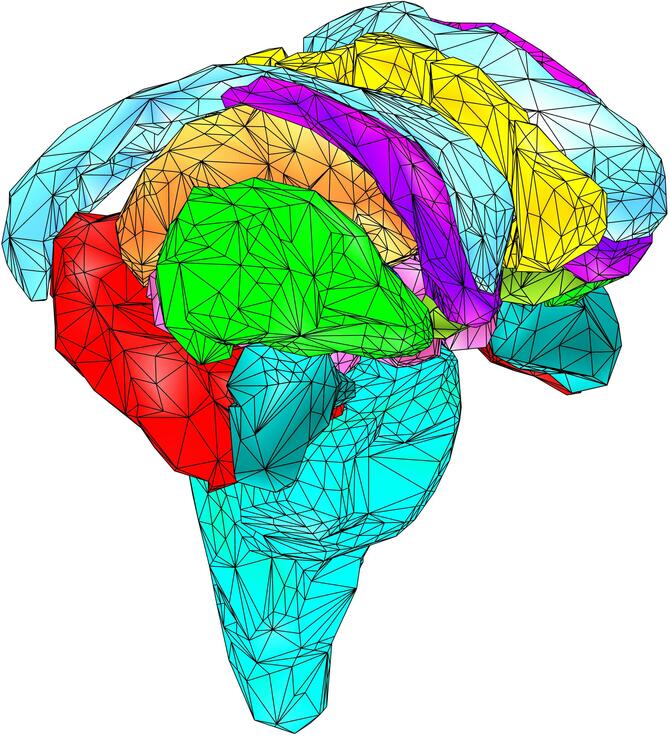}
            \end{subfigure}
            \begin{subfigure}[T]{1.60cm}
                \centering
                \includegraphics[width=1.55cm]{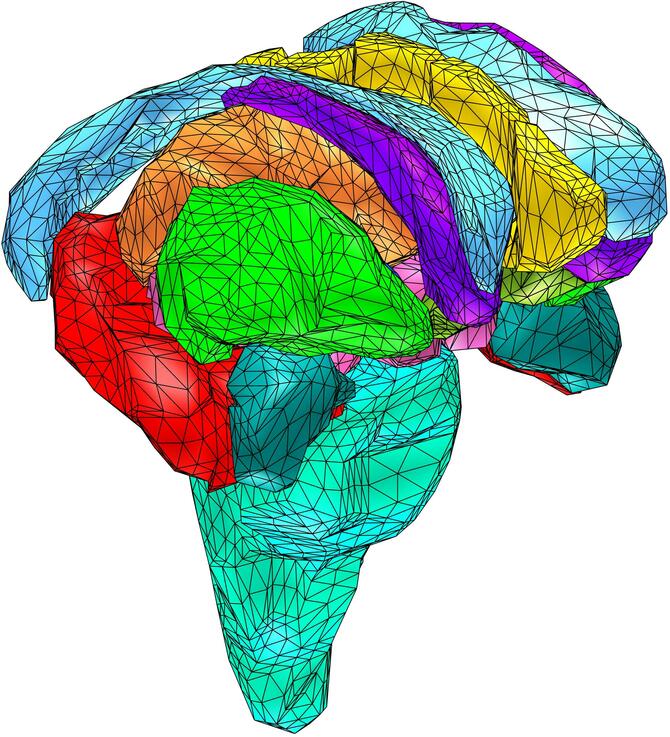}
            \end{subfigure}
            \begin{subfigure}[T]{1.60cm}
                \centering
                \includegraphics[width=1.55cm]{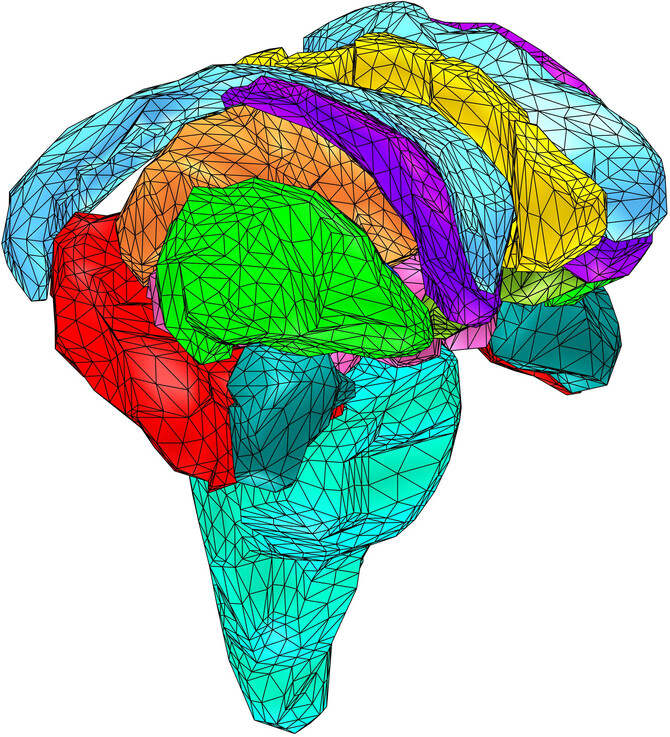}
            \end{subfigure}
            \\ \vskip0.1cm
            \begin{subfigure}[T]{1.60cm}
                \centering
                \includegraphics[width=1.55cm]{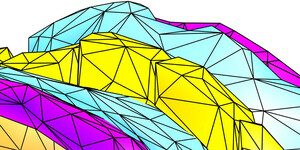}
                \caption{}
                \label{fig:deep_a}
            \end{subfigure}
            \begin{subfigure}[T]{1.60cm}
                \centering
                \includegraphics[width=1.55cm]{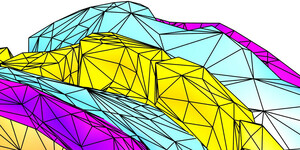}
                \caption{}
                \label{fig:deep_b}
            \end{subfigure}
            \begin{subfigure}[T]{1.60cm}
                \centering
                \includegraphics[width=1.55cm]{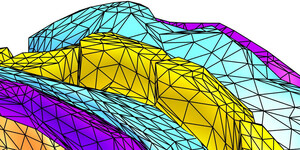}
                \caption{}
                \label{fig:deep_c}
            \end{subfigure}
            \begin{subfigure}[T]{1.60cm}
                \centering
                \includegraphics[width=1.55cm]{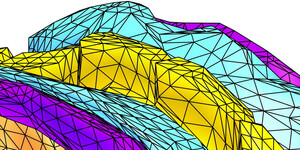}
                \caption{}
                \label{fig:deep_d}
            \end{subfigure}
        \end{minipage}
        
        \hrule
        \vskip0.2cm
        
        \begin{minipage}[T]{7.5cm}
        \centering
            \begin{subfigure}[T]{1.60cm}
            \centering
                \includegraphics[width=1.55cm]{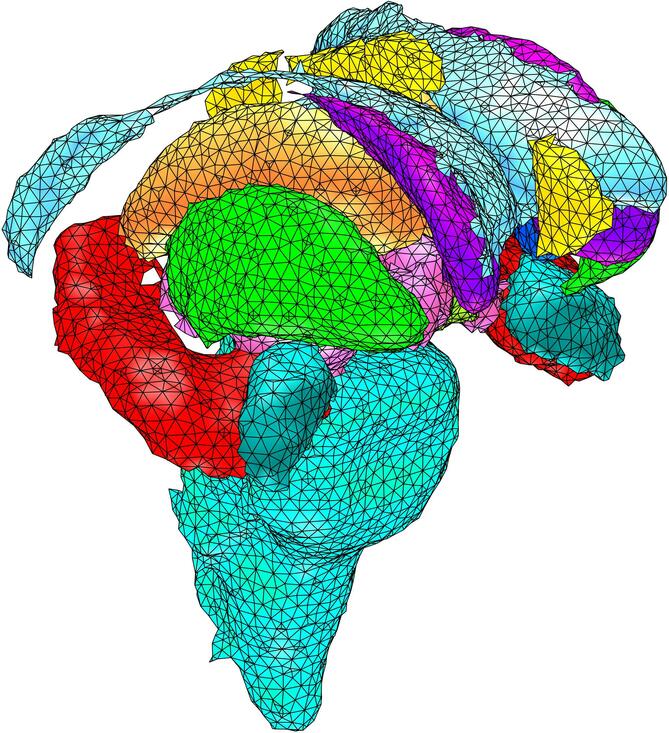}
            \end{subfigure}
            \begin{subfigure}[T]{1.60cm}
                \centering
                \includegraphics[width=1.55cm]{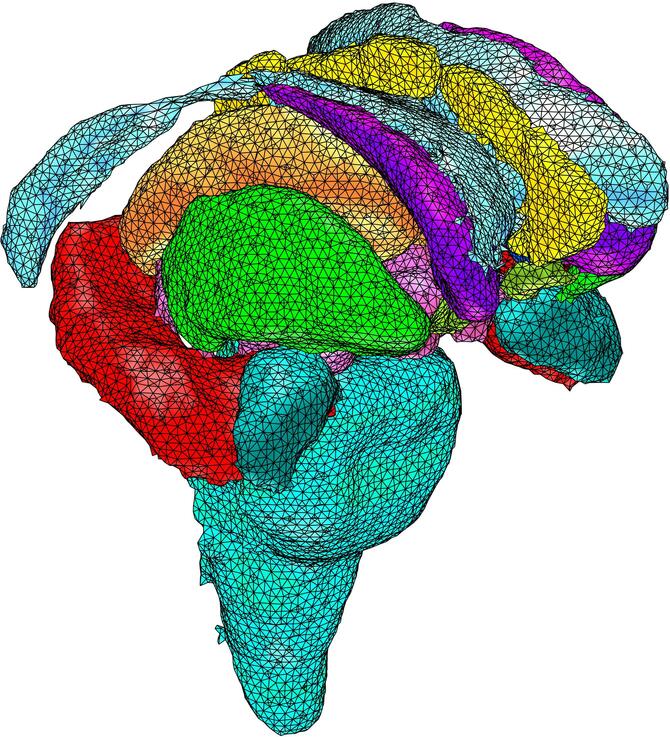}
            \end{subfigure}
            \begin{subfigure}[T]{1.60cm}
                \centering
                \includegraphics[width=1.55cm]{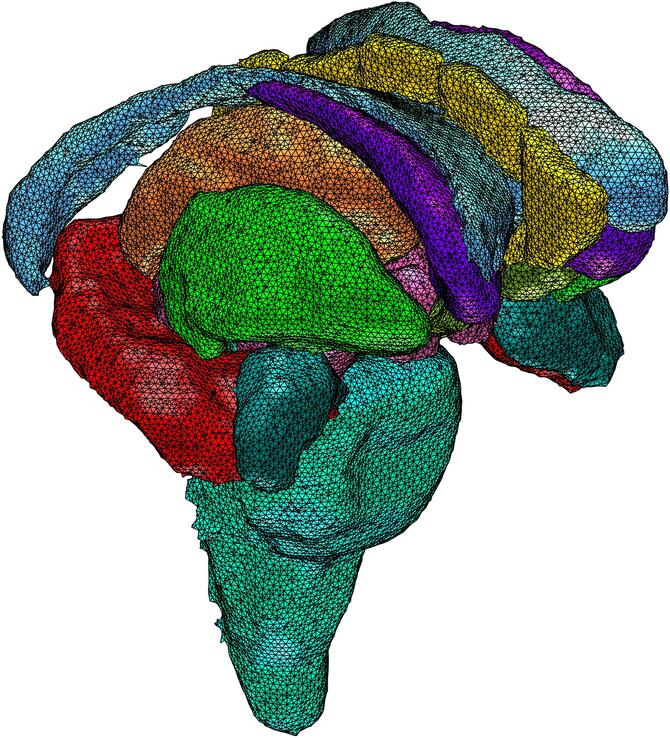}
            \end{subfigure}
            \begin{subfigure}[T]{1.60cm}
                \centering
                \includegraphics[width=1.55cm]{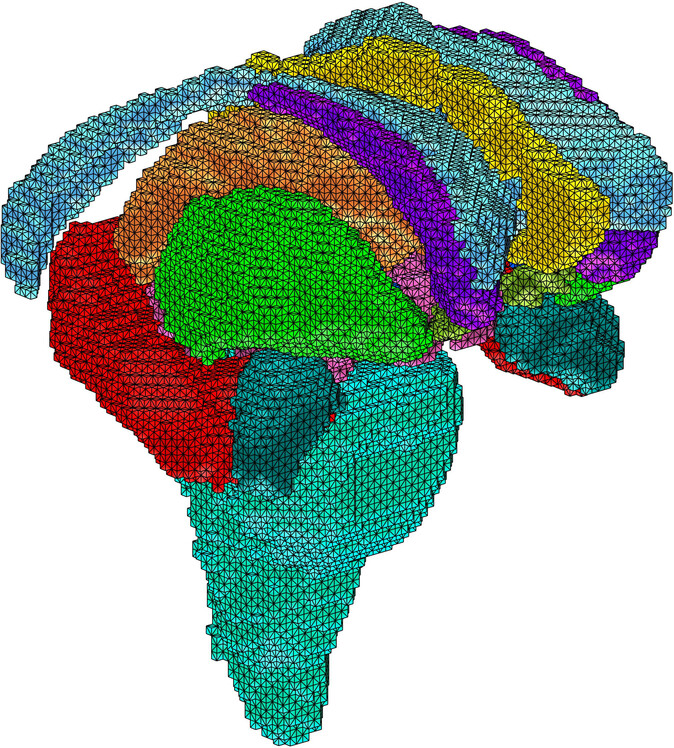}
            \end{subfigure}
            \\ \vskip0.1cm
            \begin{subfigure}[T]{1.60cm}
            \centering
                \includegraphics[width=1.55cm]{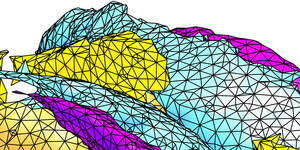}
                \caption{}
                \label{fig:deep_e}
            \end{subfigure}
            \begin{subfigure}[T]{1.60cm}
                \centering
                \includegraphics[width=1.55cm]{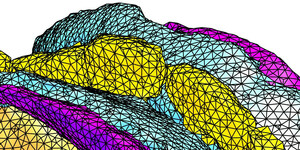}
                \caption{}
                \label{fig:deep_f}
            \end{subfigure}
            \begin{subfigure}[T]{1.60cm}
                \centering
                \includegraphics[width=1.55cm]{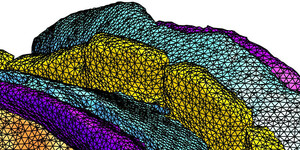}
                \caption{}
                \label{fig:deep_g}
            \end{subfigure}
            \begin{subfigure}[T]{1.60cm}
                \centering
                \includegraphics[width=1.55cm]{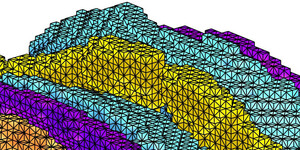}
                \caption{}
                \label{fig:deep_h}
            \end{subfigure}
        \end{minipage}
\caption{Reconstruction of the deep structures for a realistic head model featuring the brain stem (cyan), hippocampus (red), amygdala (dark green), putamen (light green), thalamus (orange), caudate (purple), ventricles (blue), and ciangulate cortex (yellow) depicting the performance of adaptive FE mesh generation using 3.0 mm, 2.0 mm and 1.3 mm mesh sizes (Surface: \ref{fig:deep_a}-\ref{fig:deep_c}; Volume: \ref{fig:deep_e}-\ref{fig:deep_g}, respectively) after applying recursive solid angle labeling, smoothing, inflation and optimization procedures in contrast to a regular 1.0 mm mesh generation without post-process procedures (\ref{fig:deep_d} and \ref{fig:deep_h}). The focused image beneath illustrates that with a coarser mesh size over a complex domain the solver can yield inaccurate results, e.g., an incomplete cingulate cortex, while using a finer mesh size, e.g., 2.0 mm or lower, can model it with accuracy.}
\label{fig:carsten_deep}
\end{scriptsize}
\end{figure}

The downsampled surface segmentation (triangular) and volumetric (tetrahedral) FE mesh parameters for each model are presented in Table \ref{table:nodes_tetras}. The computing times (measured in seconds) for each case have been included in Tables \ref{table:cpu_gpu_sphere} and \ref{table:cpu_gpu_head}. We obtained the finest FE mesh with a 1.3 mm initial lattice resolution. This mesh resulted in 5.0 M (million) nodes and 27 M tetrahedrons for the sphere, and 7.5 M nodes and 40 M tetrahedrons for the realistic head model.

\begin{table*}
\centering
\begin{scriptsize}
    \caption{The number of points and triangles for (downsampled) surface segmentations together with the node and tetra count  for the FE meshes obtained using the spherical three-layer Ary model and a realistic head model.}
    \label{table:nodes_tetras}
        \centering
        \begin{tabular}{lrrrr}
            \toprule
            \multicolumn{1}{c}{} & \multicolumn{3}{c}{\textbf{Adapted}} & \multicolumn{1}{r}{\textbf{Regular}} \\
            \multicolumn{1}{l}{\textbf{Property}} & \textbf{3.0mm} & \textbf{2.0mm} & \textbf{1.3mm} & \textbf{1.0mm} \\
            \midrule
            Points      & 34,773        & 42,284        & 42,284        & 42,284     \\
            Triangles   & 72,556        & 88,636        & 88,636        & 88,636     \\
            Nodes       & 607,184       & 1,833,945     & 5,003,931     & 4,269,175  \\
            Tetras      & 3,245,733     & 9,903,208     & 26,801,840    & 20,871,950 \\
            Size (MB)   & 707           & 1,860         & 5,092         & 6,612      \\
            \bottomrule
            \multicolumn{5}{l}{\em{(a)} Spherical overview. }
        \end{tabular}
        \vskip0.2cm
        \begin{tabular}{lrrrr}
            \toprule
            \multicolumn{1}{c}{} & \multicolumn{3}{c}{\textbf{Adapted}} & \multicolumn{1}{r}{\textbf{Regular}} \\
            \multicolumn{1}{l}{\textbf{Property}} & \textbf{3.0mm} & \textbf{2.0mm} & \textbf{1.3mm} & \textbf{1.0mm} \\
            \midrule
            Points      & 46,555        & 104,684       & 186,057       & 186,057    \\
            Triangles   & 93,016        & 209,314       & 372,088       & 372,088    \\
            Nodes       & 837,029       & 2,517,262     & 7,526,427     & 4,798,880  \\
            Tetras      & 4,540,361     & 13,699,905    & 40,895,474    & 23,460,830 \\
            Size (MB)   & 1,297         & 2,812         & 7,254         & 6,463      \\
            \bottomrule
            \multicolumn{5}{l}{\em{(b)} Head model overview.}
        \end{tabular}
    \end{scriptsize}
\end{table*}

The results suggest that the recursive re-labeling method applied in adaptive mesh generation allows faster meshing for the mesh size as compared to the regular approach. Most computing effort is spent on labeling, which can be accelerated using GPU parallelization; GPU labeling took, in most cases, less than $1/10$ of the time consumed by the CPU. The labeling effect was most influential in the case of regular mesh, where the recursive surface extraction (Section \ref{sec:labeling}) was not applied. In adaptive meshing, the GPU accelerated total meshing time was close to  1/3--1/2 compared to the case of no acceleration. The post-processing phase (smoothing, inflation, and optimization) took roughly 30-40 \% of the total required time by the GPU accelerated meshing.

The distance between the given and the reconstructed grey matter boundary is shown in Figures (\ref{fig:sphere_hist}) and (\ref{fig:sphere_hist}). In the spherical case, the median of this distance is close to one-fourth of the initial mesh resolution, i.e., one-half of the refined mesh size in the vicinity of the boundary. For the realistic case, the errors were slightly larger and more dispersed, the median and spread (inter-quartile range between the 25 and 75 \% quantiles) being smaller than the refined mesh size. The regular mesh has a comparably large spread following the absence of the smoothing and inflation steps. 

The element condition distributions are shown in Figure \ref{fig:element_condition}. While the meshes distributions for the spherical and realistic meshes share similarities, we assume that the proportion of the elements with conditions below 0.01 for the realistic mesh compared to the spherical model is a result of having a greater number of smoothed surfaces. Spatial condition number mappings for the realistic head model showed, as expected, that non-refined internal parts of the mesh, particularly the white matter interior, can reach an elevated condition compared to thin, refined layers closer to the surface.

The effect of the mesh adaptation is illustrated in Figures \ref{fig:adapted_slice_sphere} and \ref{fig:adapted_slice_realistic} for the spherical and realistic cases, respectively. The effect of surface refinement is visible in the case of an adapted 1.3 mm mesh; the resolution of an adapted FE mesh decreases towards the interior in the compartments that share the refined boundaries, while in the case of regular mesh is the same anywhere. 

In Figure \ref{fig:carsten_surface}, we show reconstructions of the cerebrum, cerebellum, and brain stem. Figure \ref{fig:carsten_deep} displays the deep brain components, including the hippocampus, amygdala, putamen, thalamus caudate, ventricles, and cingulate cortex, from the same volumetric MRI data set described in \ref{sec:head_segmentation}. An increase in the mesh resolution improves the shape and structure of the modeled tissues, which is particularly prominent for the caudal anterior area of the cingulate cortex, not fully reconstructed in the case of the 3.0 mm resolution, attributed to a comparably coarse initial resolution. 

\subsection{EEG forward modelling}
\begin{figure}
\centering  
\begin{scriptsize}
    \begin{subfigure}[T]{7.5cm} \centering
        \includegraphics[width=7.5cm]{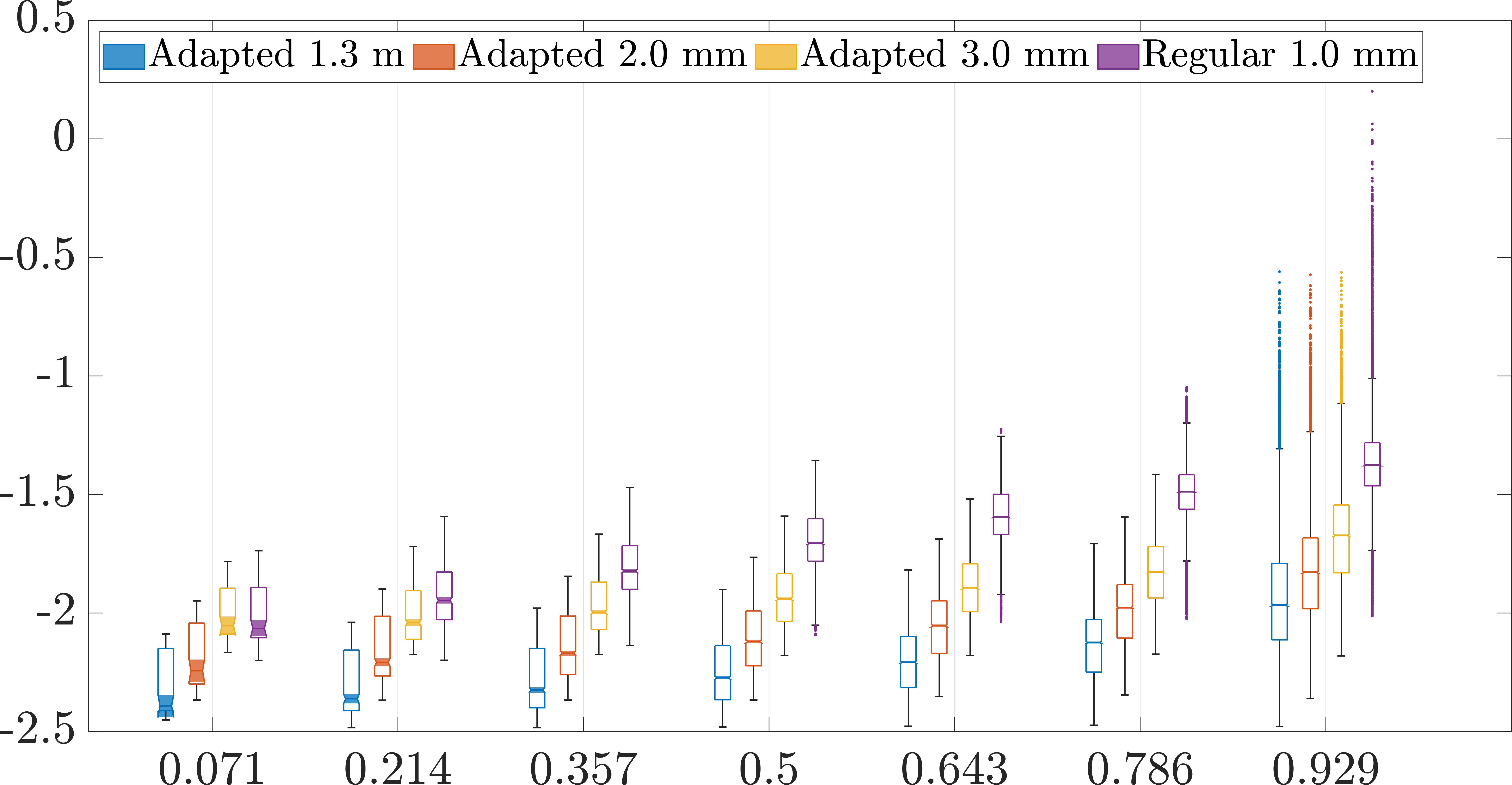}
        \caption{Full-range Relative Difference (RDM)}
        \label{fig:rdm_all}
    \end{subfigure}
    \begin{subfigure}[T]{7.5cm} \centering
        \includegraphics[width=7.5cm]{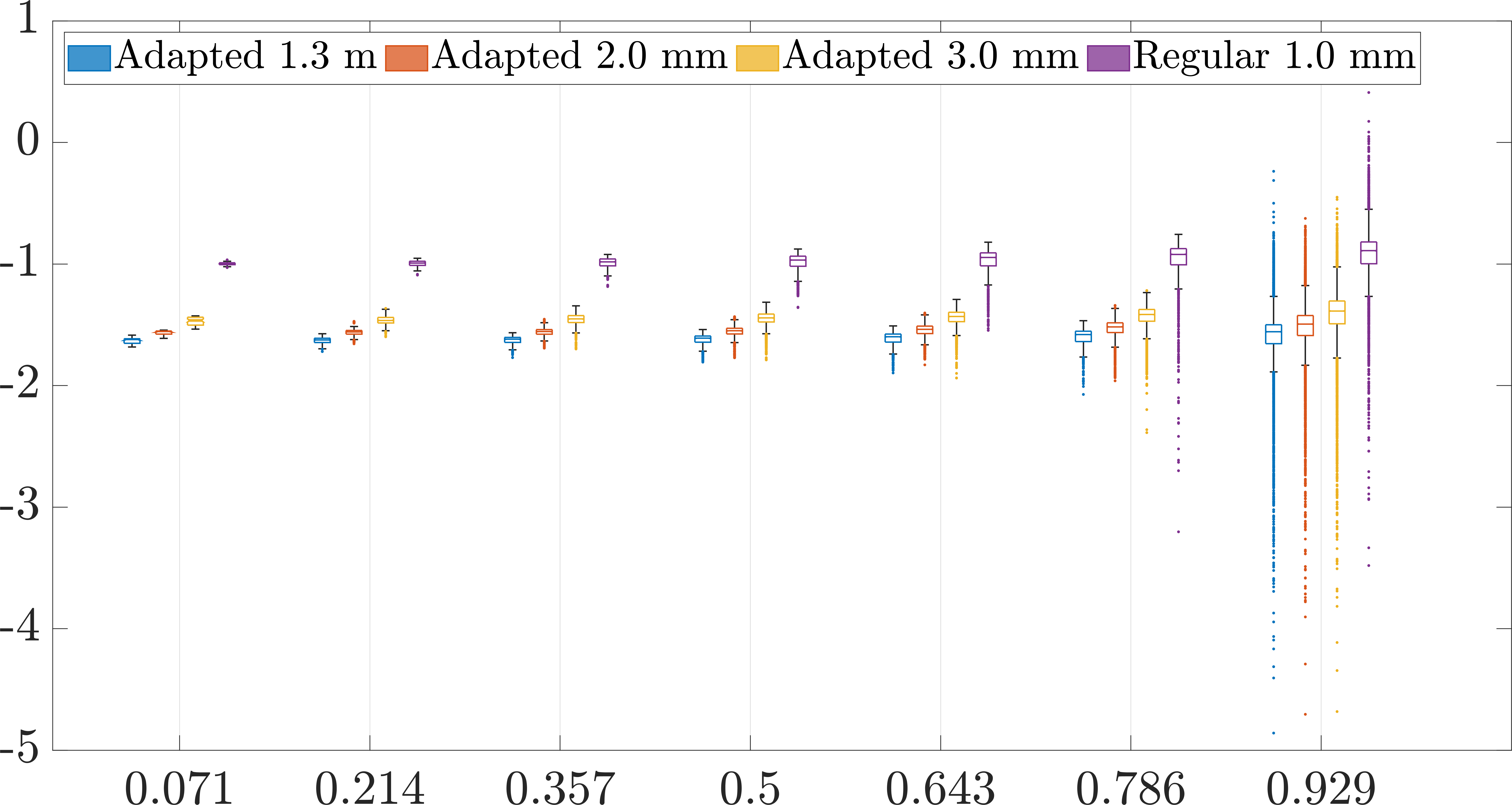}
        \caption{Full-range Magnitude (MAG)}
        \label{fig:mag_all}
    \end{subfigure} 
    \begin{subfigure}[T]{7.5cm} \centering
        \includegraphics[width=7.5cm]{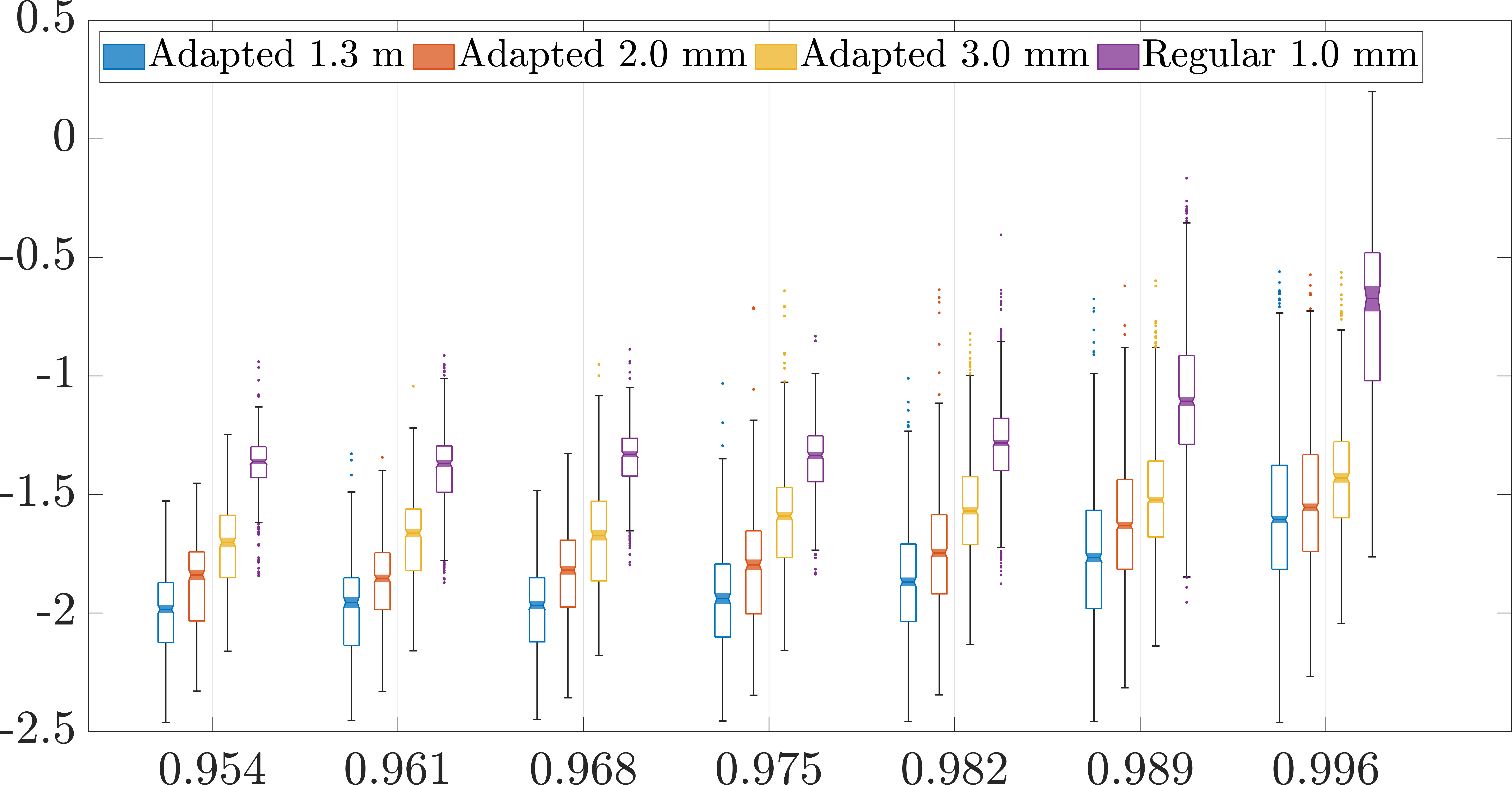}
        \caption{95\% - 99\% Relative Difference (RDM)}
        \label{fig:rdm_high}
    \end{subfigure}
    \begin{subfigure}[T]{7.5cm} \centering
        \includegraphics[width=7.5cm]{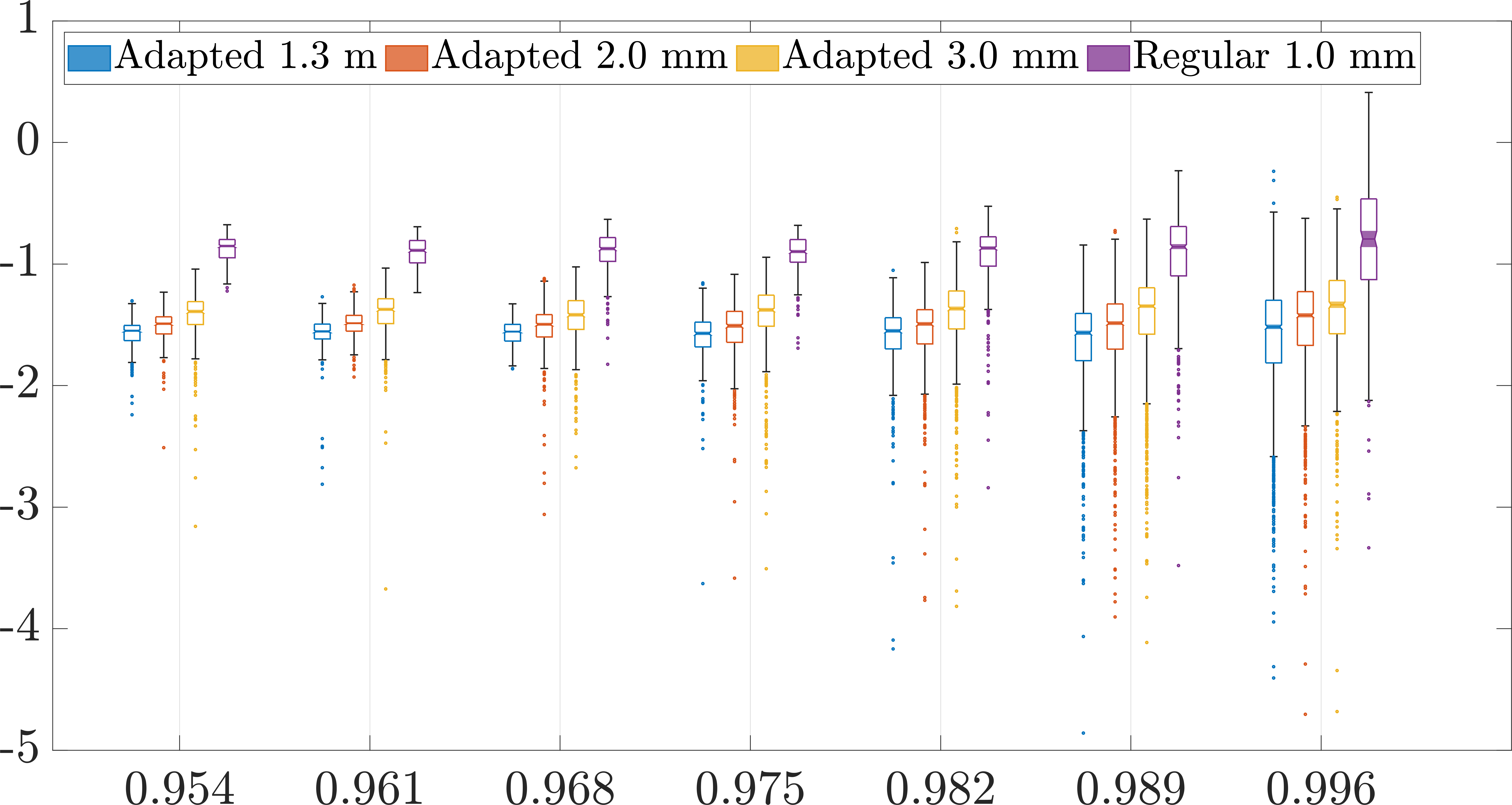}
        \caption{95\% - 95\% Magnitude (MAG)}
        \label{fig:mag_high}
    \end{subfigure}
\caption{Accuracy of the spherical three-layer Ary model in EEG forward modelling evaluated at different eccentricities levels, i.e. relative distances from the origin in the brain compartment using Relative Difference (RDM) and Magnitude (MAG).}
\label{fig:rdm_mag}
\end{scriptsize}
\end{figure}

The RDM and MAG obtained with the spherical geometry are visualized in Figure \ref{fig:rdm_mag}. The adaptive meshes resulted in systematically more accurate lead field matrices, where the smallest differences were obtained in the adapted 1.3 mm case. For each lead field matrix, the RDM had a median of below 5 \% and a MAG of below 20 \% up to 0.95. For the adapted meshes, these medians did not exceed 3 and 5 \%. With the finest adapted mesh, the upper limits of 3 and 4 \% were maintained for eccentricities up to 0.998. The effect of adapting was pronounced towards the high eccentricities, as is shown by the growing spread of the RDM and MAG distributions.

\subsection{Source localization}
\begin{figure}[t]
\centering
    \begin{scriptsize}
        \begin{subfigure}[T]{7.5cm}
        \centering
            \includegraphics[width=7.0cm]{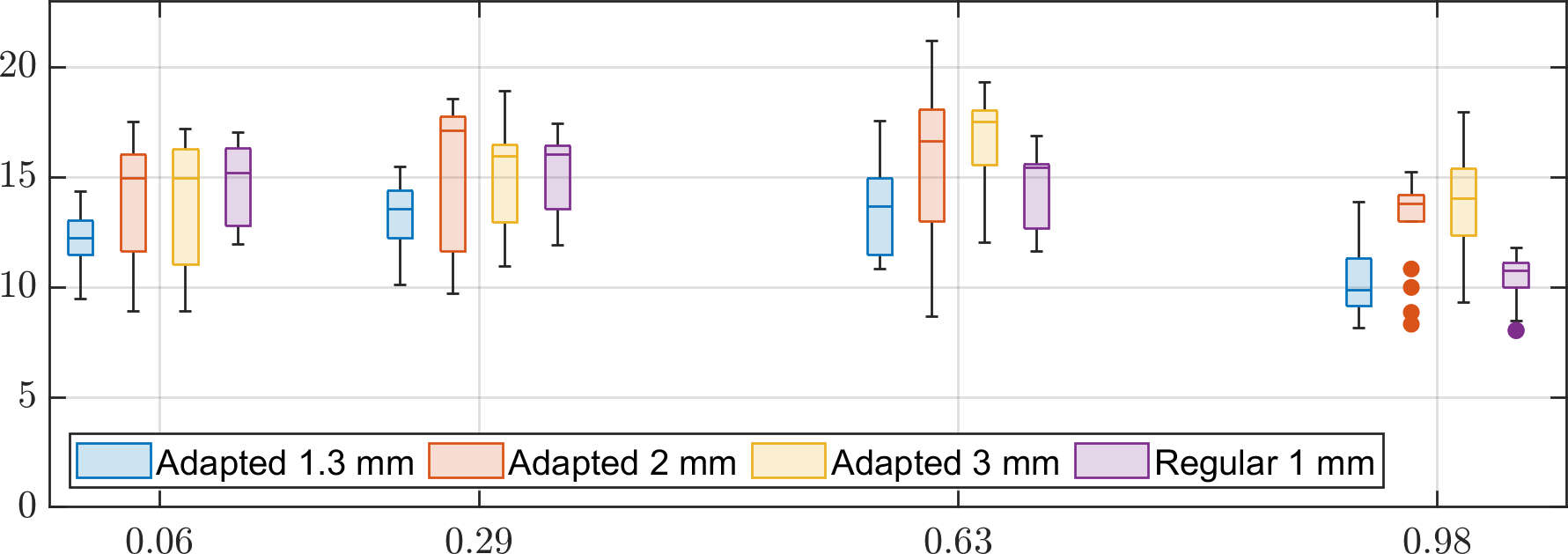}
            \caption{MNE}
            \label{fig:ecc_MNE}
        \end{subfigure}
        \vskip0.1cm
        \begin{subfigure}[T]{7.5cm}
        \centering
            \includegraphics[width=7.0cm]{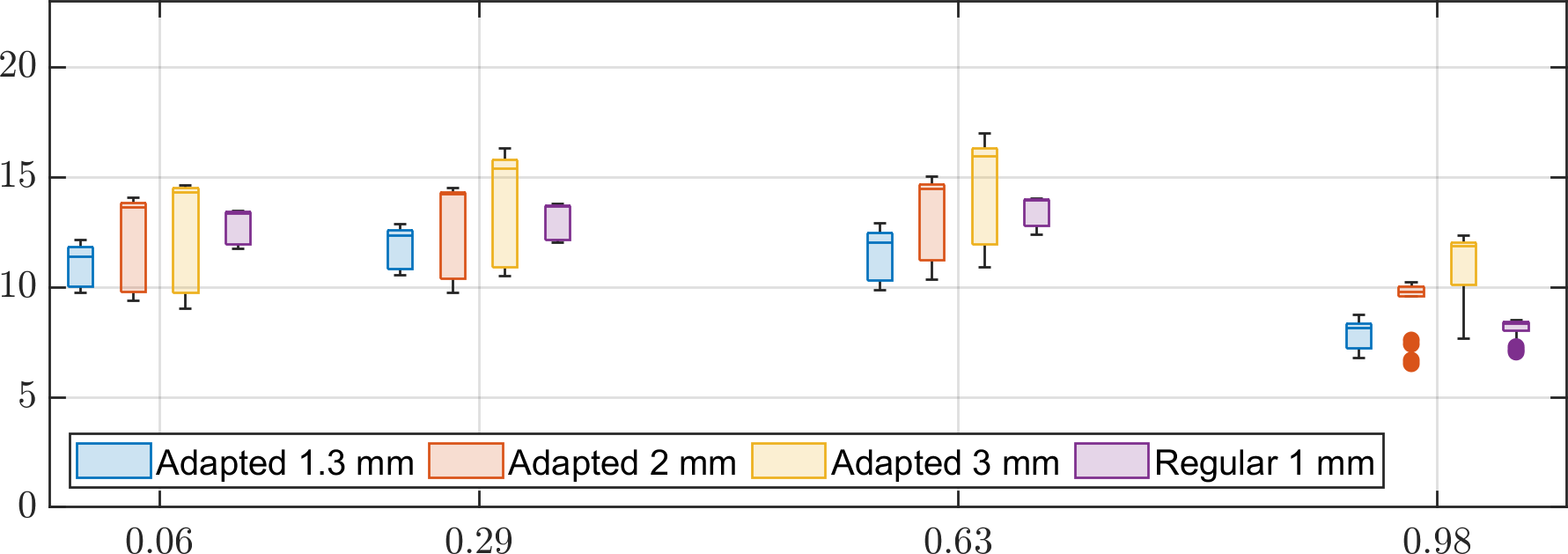}
            \caption{sLORETA}
            \label{fig:ecc_sLORETA}
        \end{subfigure}
        \vskip0.1cm
        \begin{subfigure}[T]{7.5cm}
        \centering
            \includegraphics[width=7.0cm]{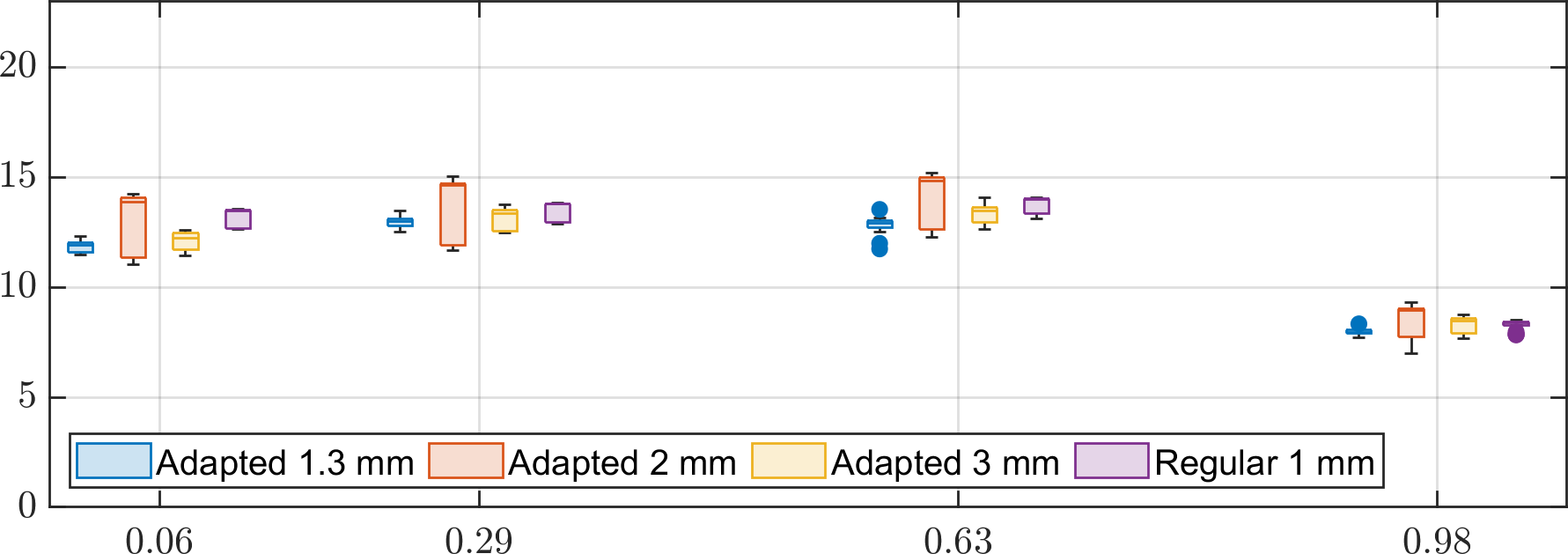}
            \caption{Dipole scan}
            \label{fig:ecc_DS}
        \end{subfigure}
    \end{scriptsize}
\caption{The earth mover's distance (EMD) measuring the amount of work required to transfer a distributional reconstruction into a dipole distribution (Dirac's delta) using Minimun Norm Estimation (MNE) (\ref{fig:ecc_MNE}), Standardized low-resolution brain electromagnetic tomography (sLORETA) (\ref{fig:ecc_sLORETA}), and dipole scan (\ref{fig:ecc_DS}).}
\label{fig:ecc_4methods}
\end{figure}

The source localization estimates for the spherical geometry and four different eccentricities, 0.06, 0.29, 0.63, and 0.98 \%, are shown in Figure (\ref{fig:ecc_4methods}). The relative mutual differences between the adaptive and regular meshing approaches, measured via EMD, were observed to be less prominent in source localization as compared to forward modeling. These EMD differences also depend on the applied inverse method. Matching with the forward simulation results, the adaptive 1.3 mm mesh provided the smallest median in each case, with a maximum of 3.0 mm marginal to the other cases. For MNE and sLORETA, the initial FE mesh resolution (3.0, 2.0, and 1.3 mm) seems to be the governing factor determining the differences in source localization accuracy, while the dipole scan, such a tendency was not observed with the mutual differences between the methods being minor. In each case, the EMDs grow towards the center of the domain (median 8--14 mm in the vicinity of the boundary and 12--15 mm near the center), which is a natural consequence of the ill-posed nature of the source localization problem. 

\subsection{SEP components P14/N22 and P22/N22}
\begin{figure}[t]
\centering
    \begin{scriptsize}
        \begin{minipage}[T]{7.5cm}
            \centering
            \begin{minipage}[T]{1.60cm}
                \centering
                \textbf{Adapted \\ 3.0 mm}
            \end{minipage}
            \begin{minipage}[T]{1.60cm}
                \centering
                \textbf{Adapted \\ 2.0 mm}
            \end{minipage}
            \begin{minipage}[T]{1.60cm}
                \centering
                \textbf{Adapted \\ 1.3 mm}
            \end{minipage}
            \begin{minipage}[T]{1.60cm}
                \centering
                \textbf{Regular \\ 1.0 mm}
            \end{minipage}
        \end{minipage}
        \begin{minipage}[T]{7.5cm}
        \centering
            \begin{subfigure}[T]{1.60cm}
                \centering
                \includegraphics[width=1.5cm]{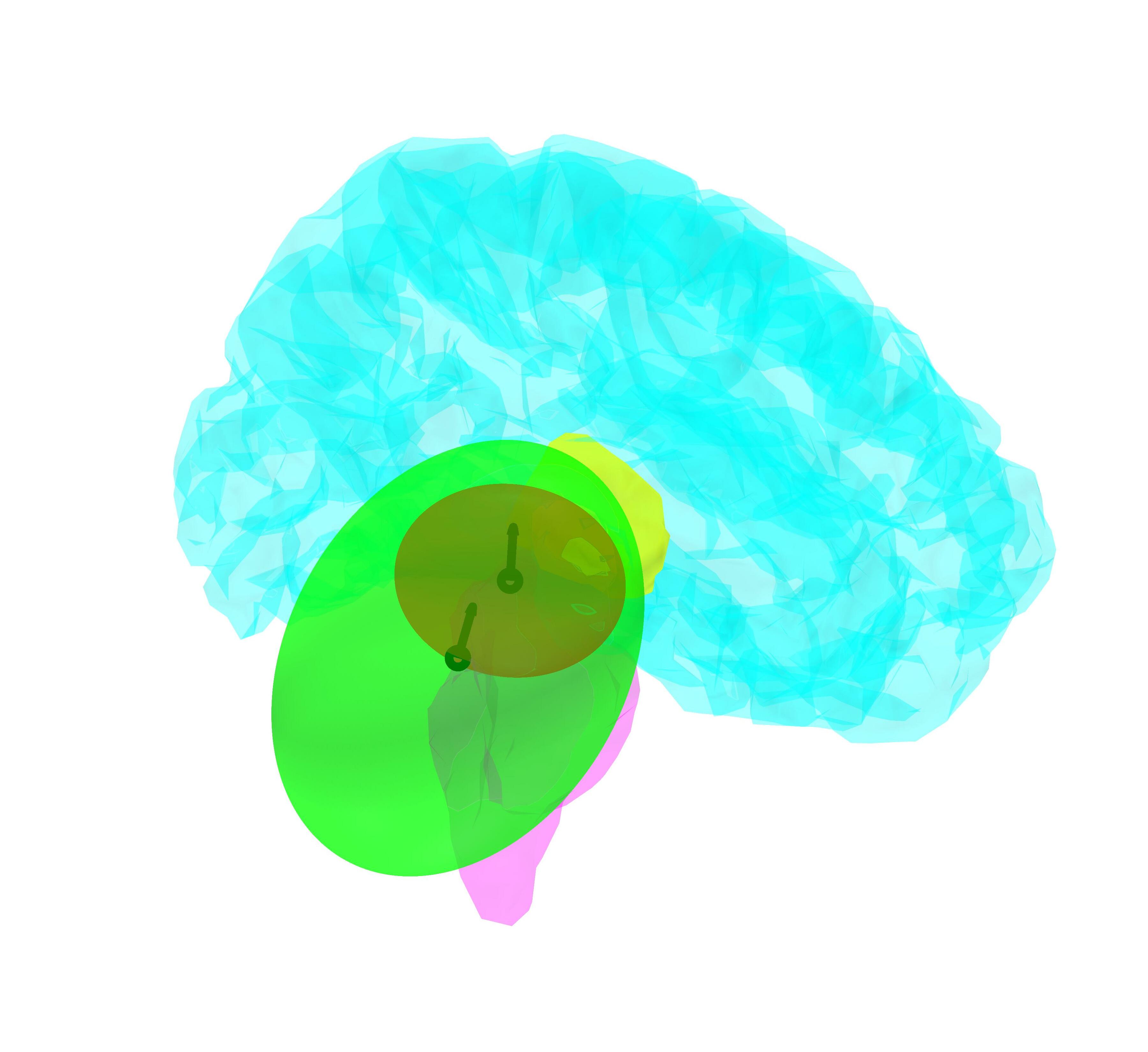}
            \end{subfigure}
            \begin{subfigure}[T]{1.60cm}
                \centering
                \includegraphics[width=1.5cm]{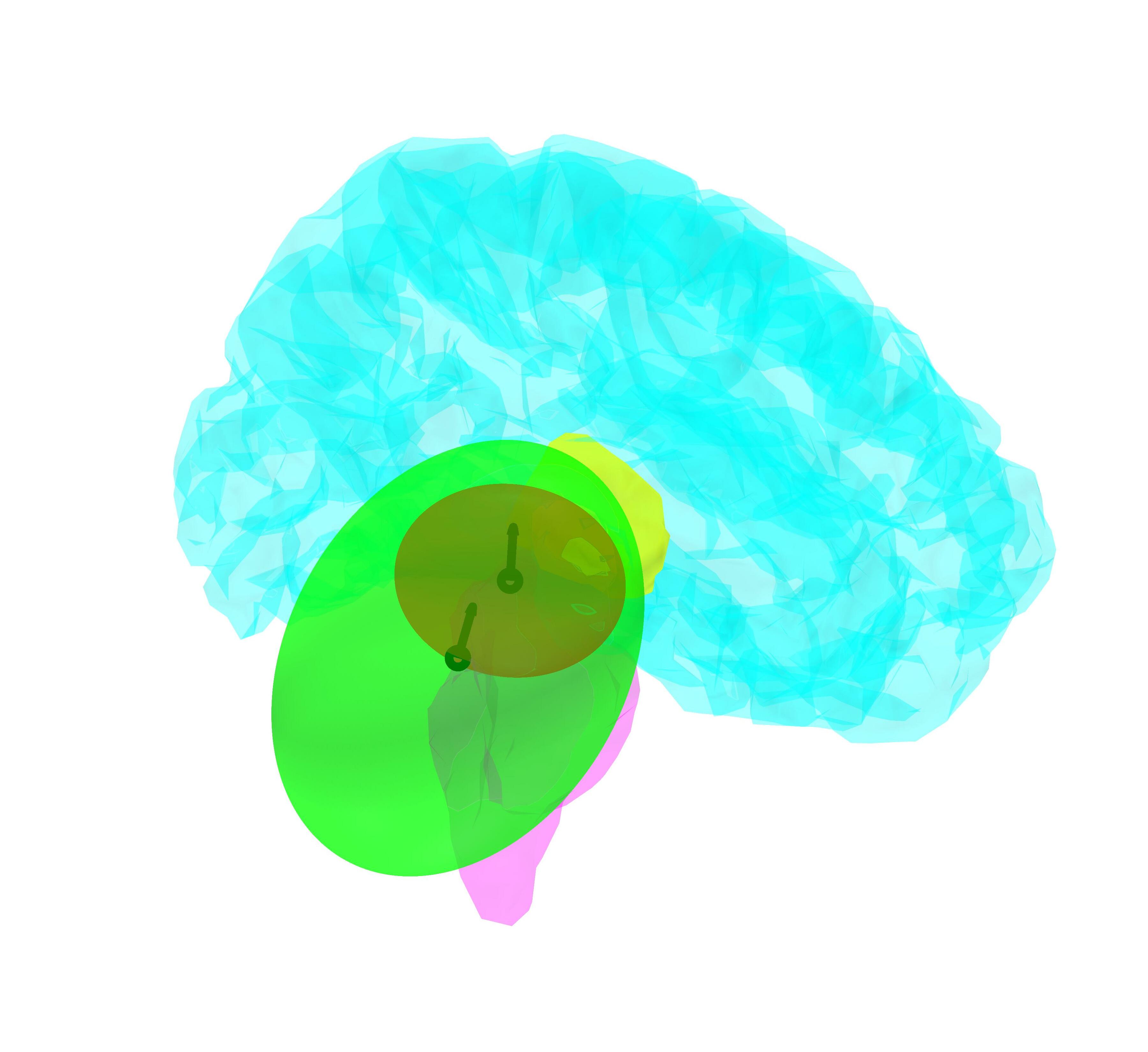}
            \end{subfigure}
            \begin{subfigure}[T]{1.60cm}
                \centering
                \includegraphics[width=1.5cm]{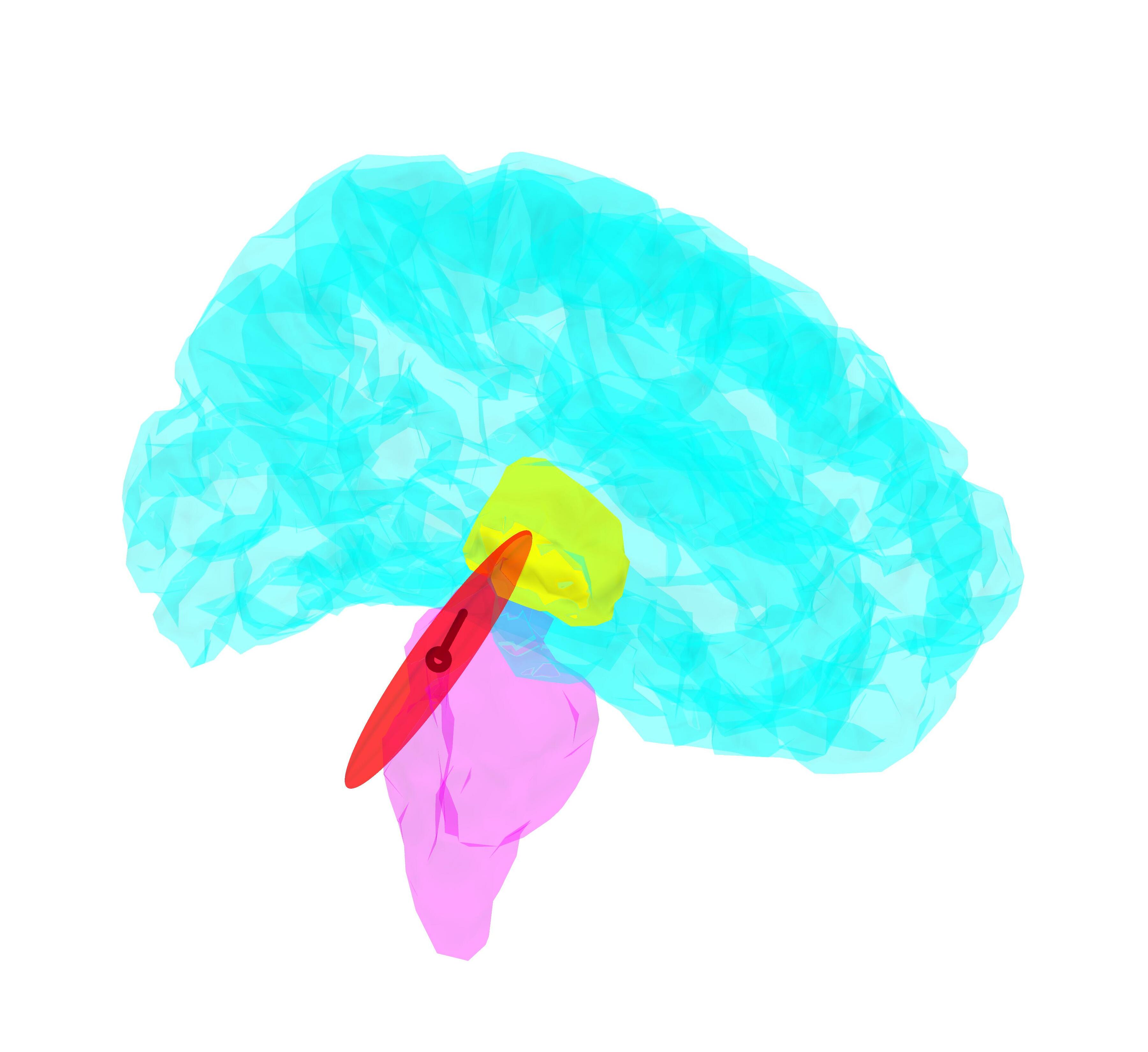}
            \end{subfigure}
            \begin{subfigure}[T]{1.60cm}
                \centering
                \includegraphics[width=1.5cm]{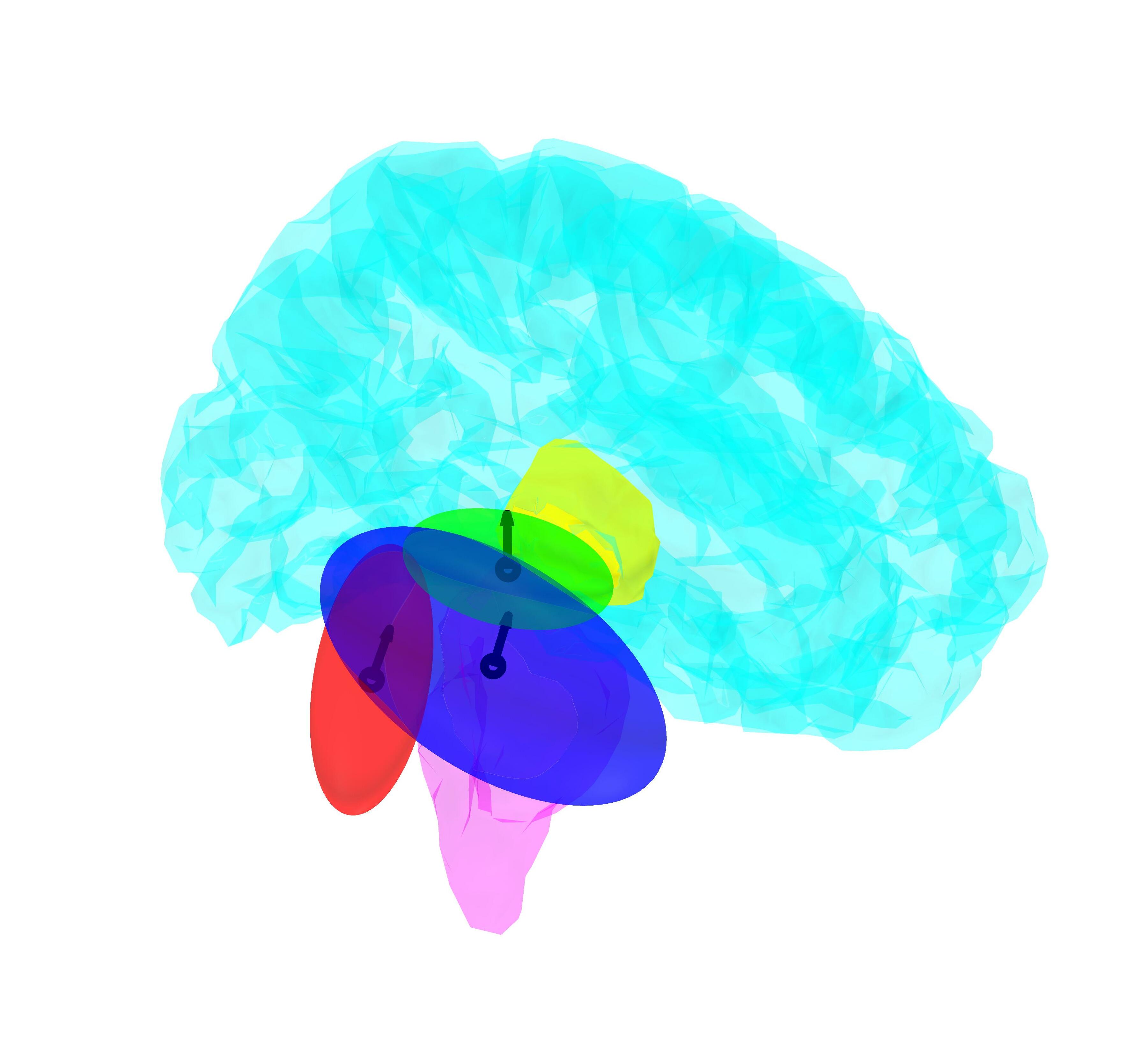}
            \end{subfigure}
            \\ 
            \begin{subfigure}[T]{1.60cm}
                \centering
                \includegraphics[width=1.5cm]{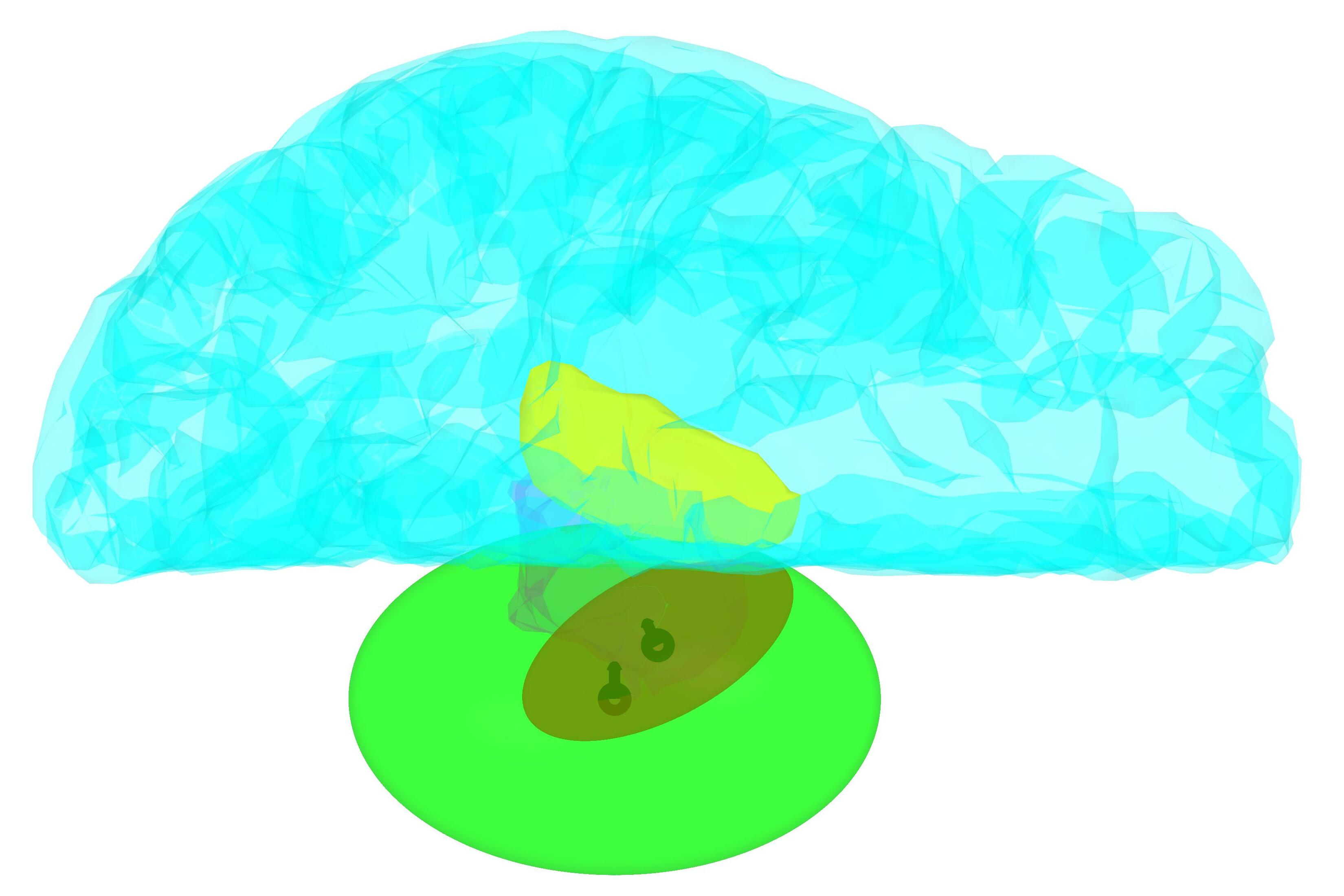}
                \caption{}
                \label{fig:GMM_3mm_p14}
            \end{subfigure}
            \begin{subfigure}[T]{1.60cm}
                \centering
                \includegraphics[width=1.5cm]{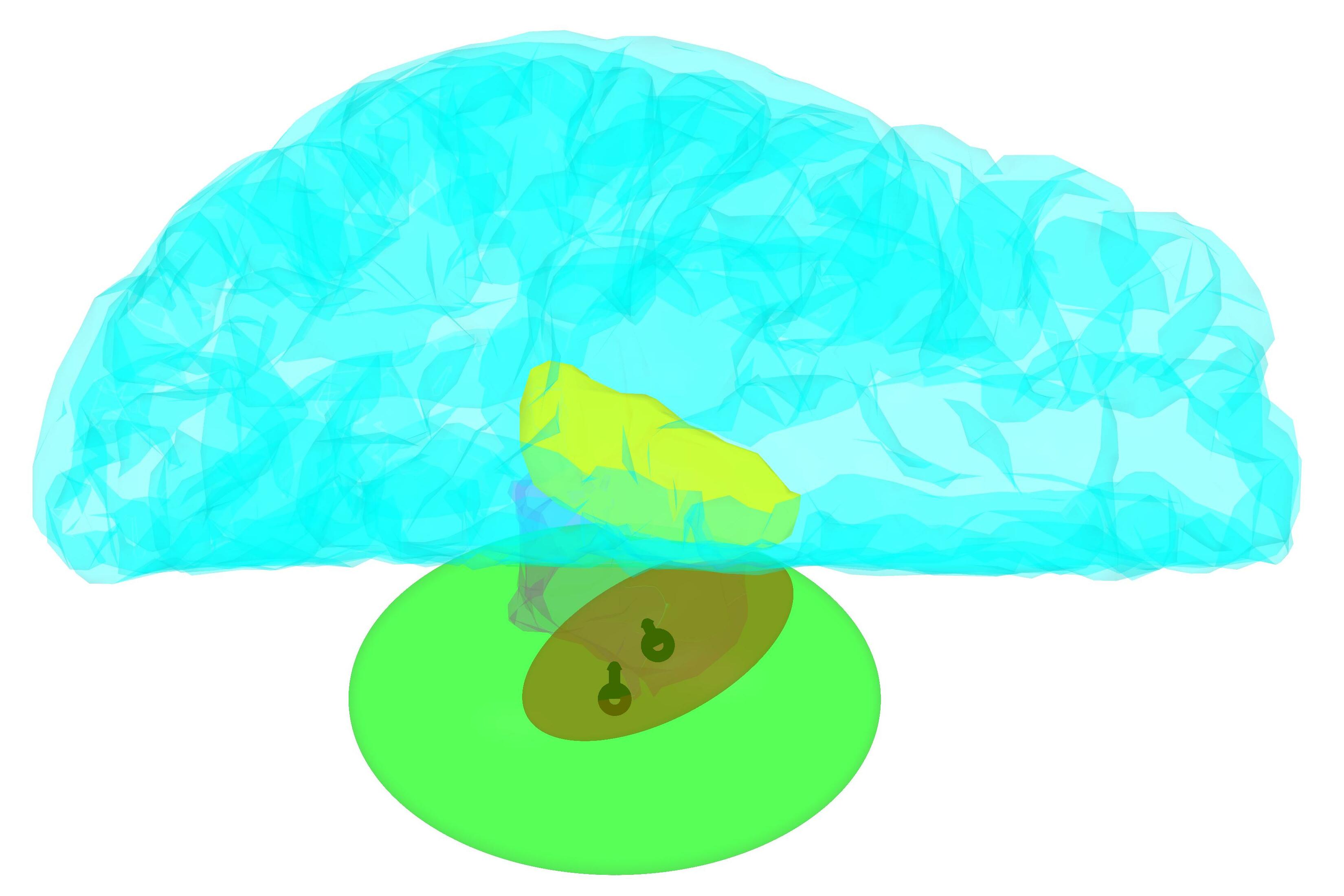}
                \caption{}
                \label{fig:GMM_2mm_p14}
            \end{subfigure}
            \begin{subfigure}[T]{1.60cm}
                \centering
                \includegraphics[width=1.5cm]{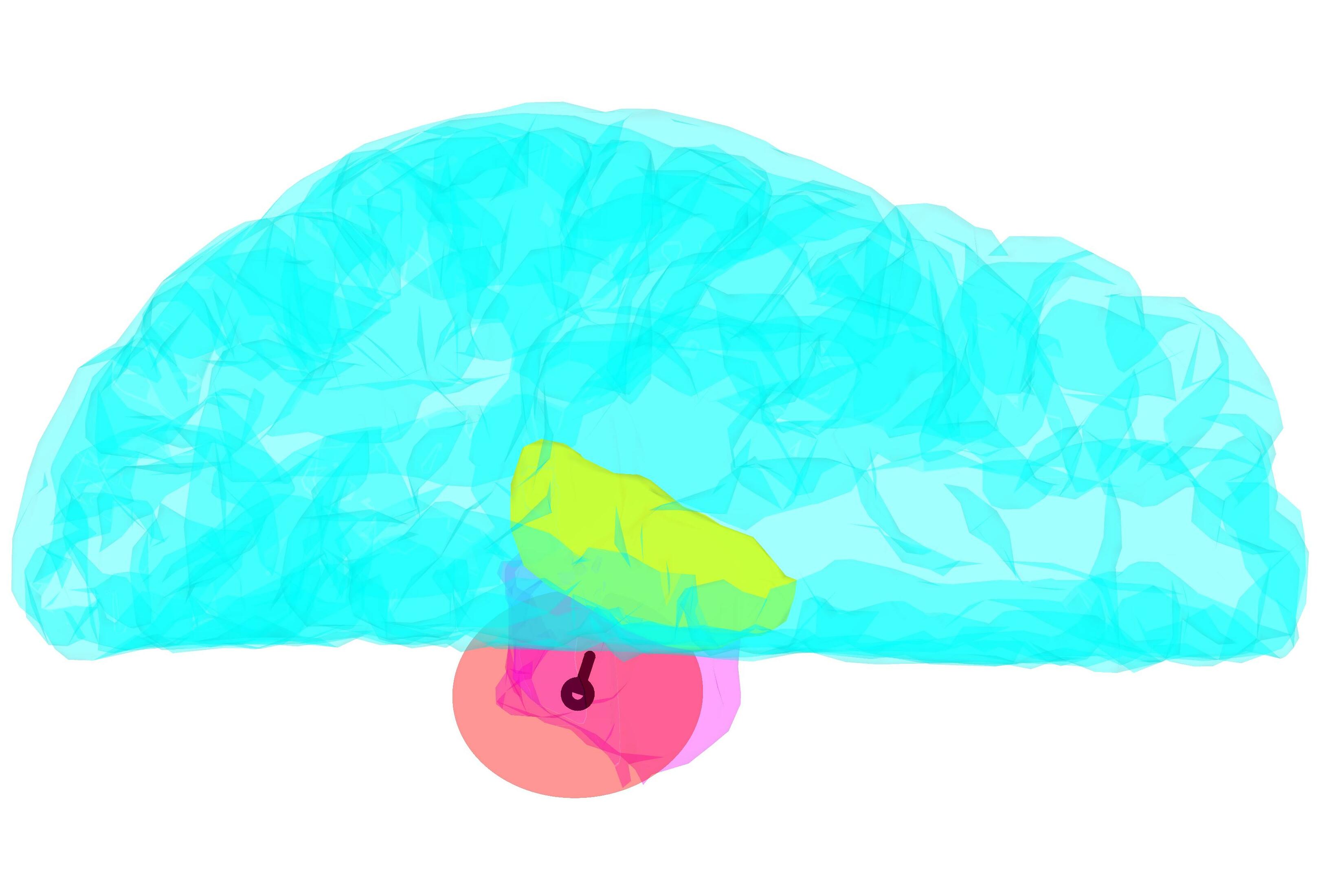}
                \caption{}
                \label{fig:GMM_130mm_p14}
            \end{subfigure}
            \begin{subfigure}[T]{1.60cm}
                \centering
                \includegraphics[width=1.5cm]{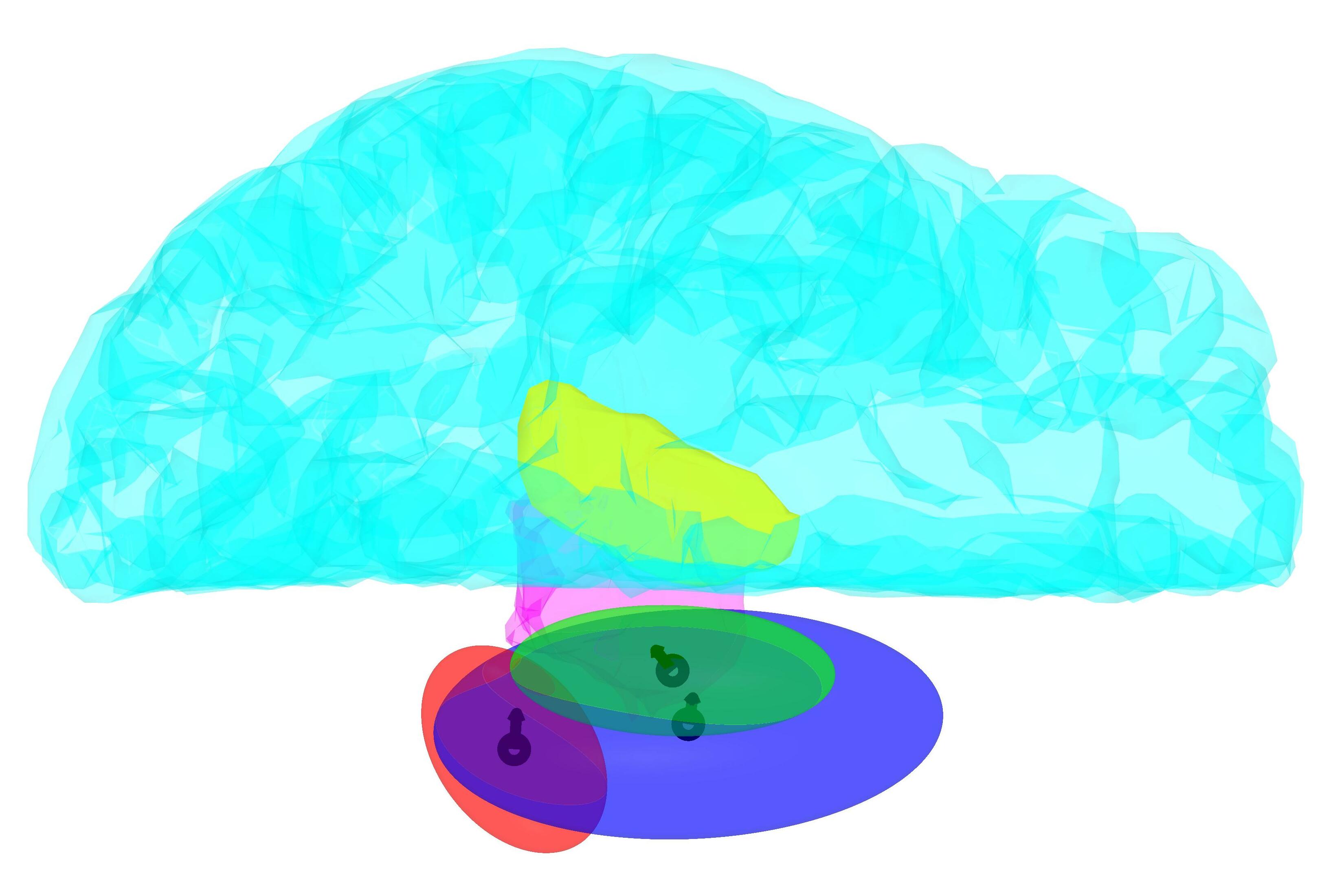}
                \caption{}
                \label{fig:GMM_1mm_p14}
            \end{subfigure}
        \end{minipage}
        \hrule
        \vskip0.4cm
        \begin{minipage}[T]{7.5cm}
            \centering
                \begin{subfigure}[T]{1.60cm}
                \centering
                   \includegraphics[width=1.5cm]{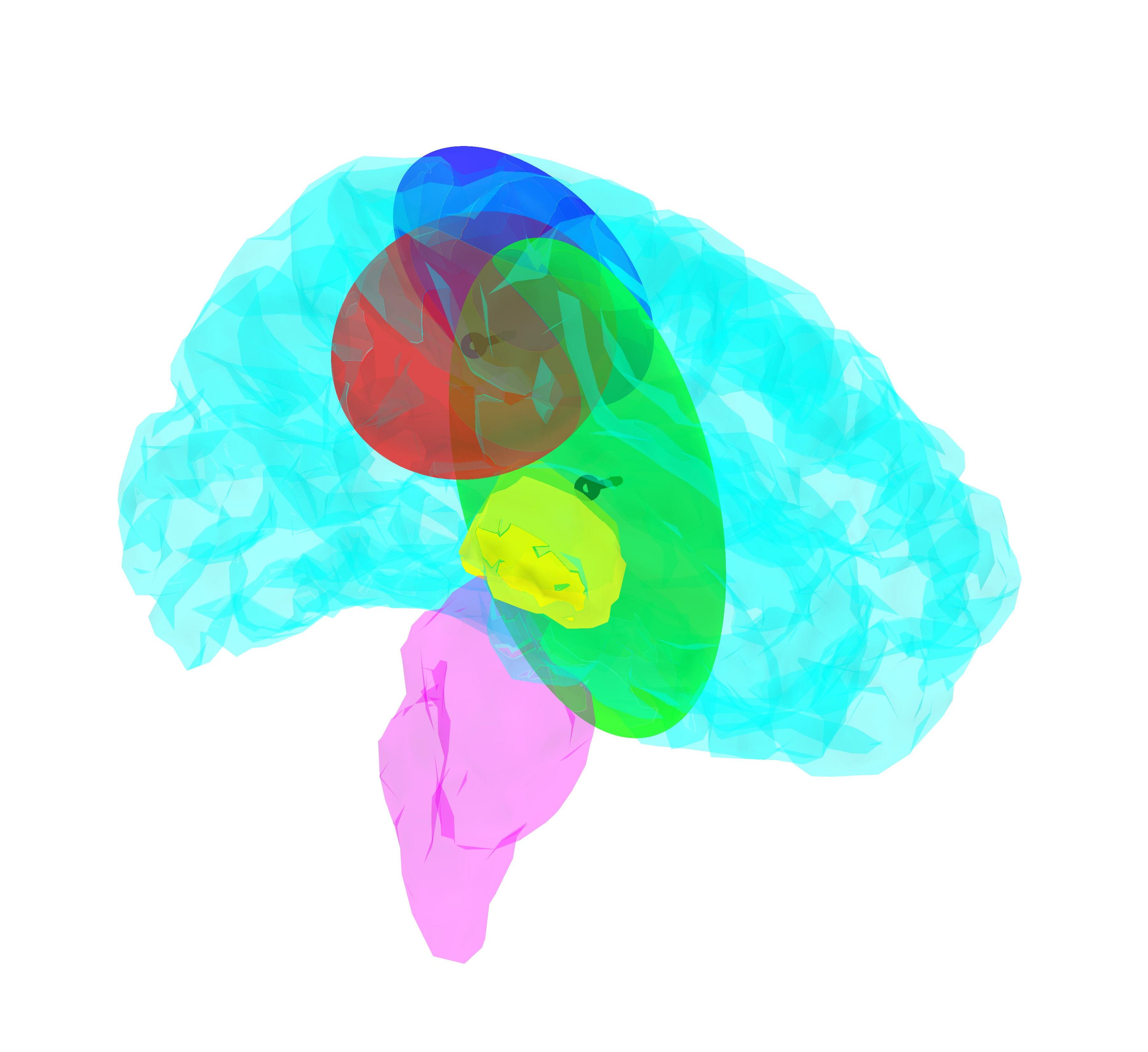}
                \end{subfigure}
                \begin{subfigure}[T]{1.60cm}
                    \centering
                    \includegraphics[width=1.5cm]{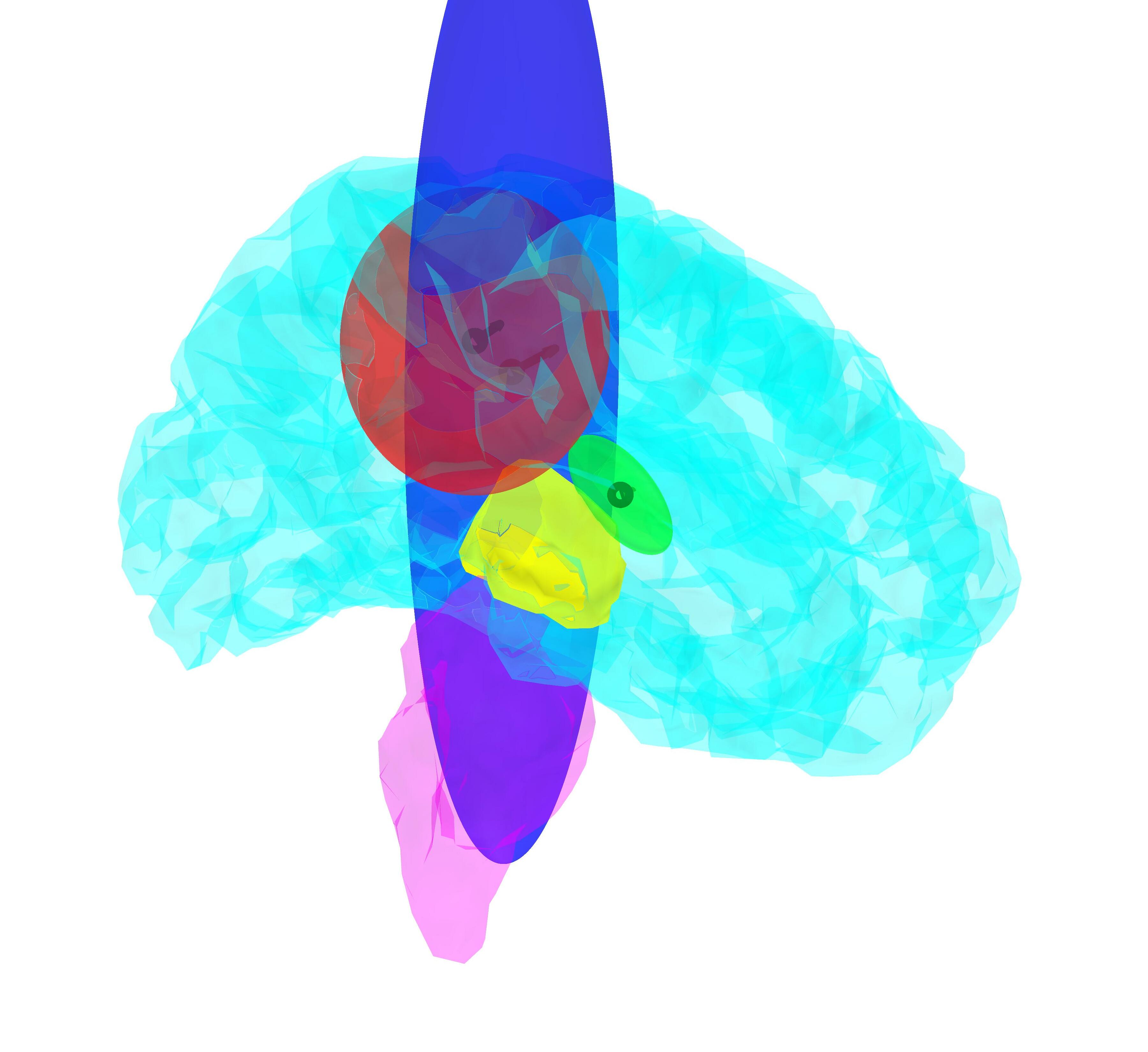}
                \end{subfigure}
                \begin{subfigure}[T]{1.60cm}
                    \centering
                    \includegraphics[width=1.5cm]{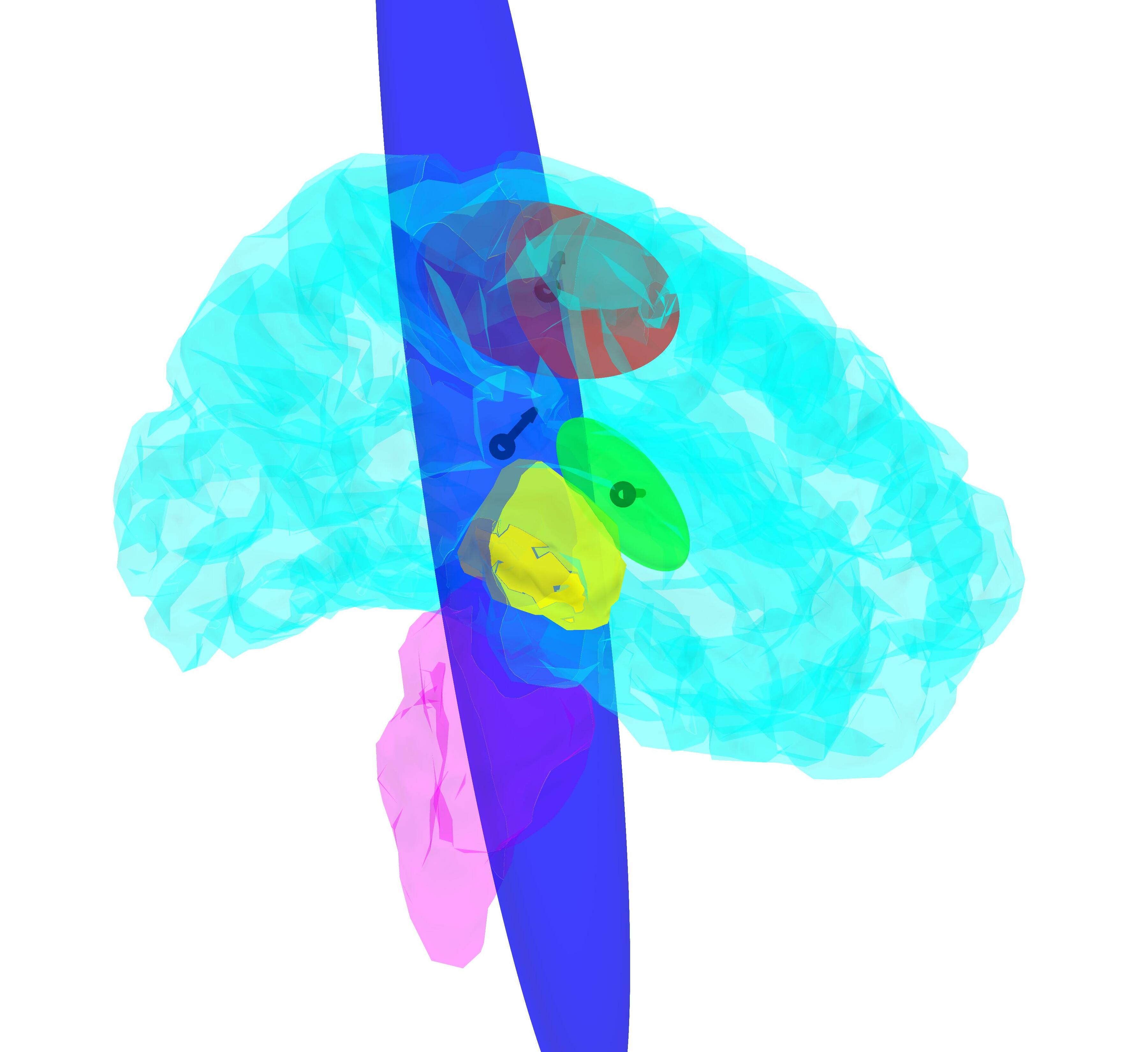}
                \end{subfigure}
                \begin{subfigure}[T]{1.60cm}
                    \centering
                    \includegraphics[width=1.5cm]{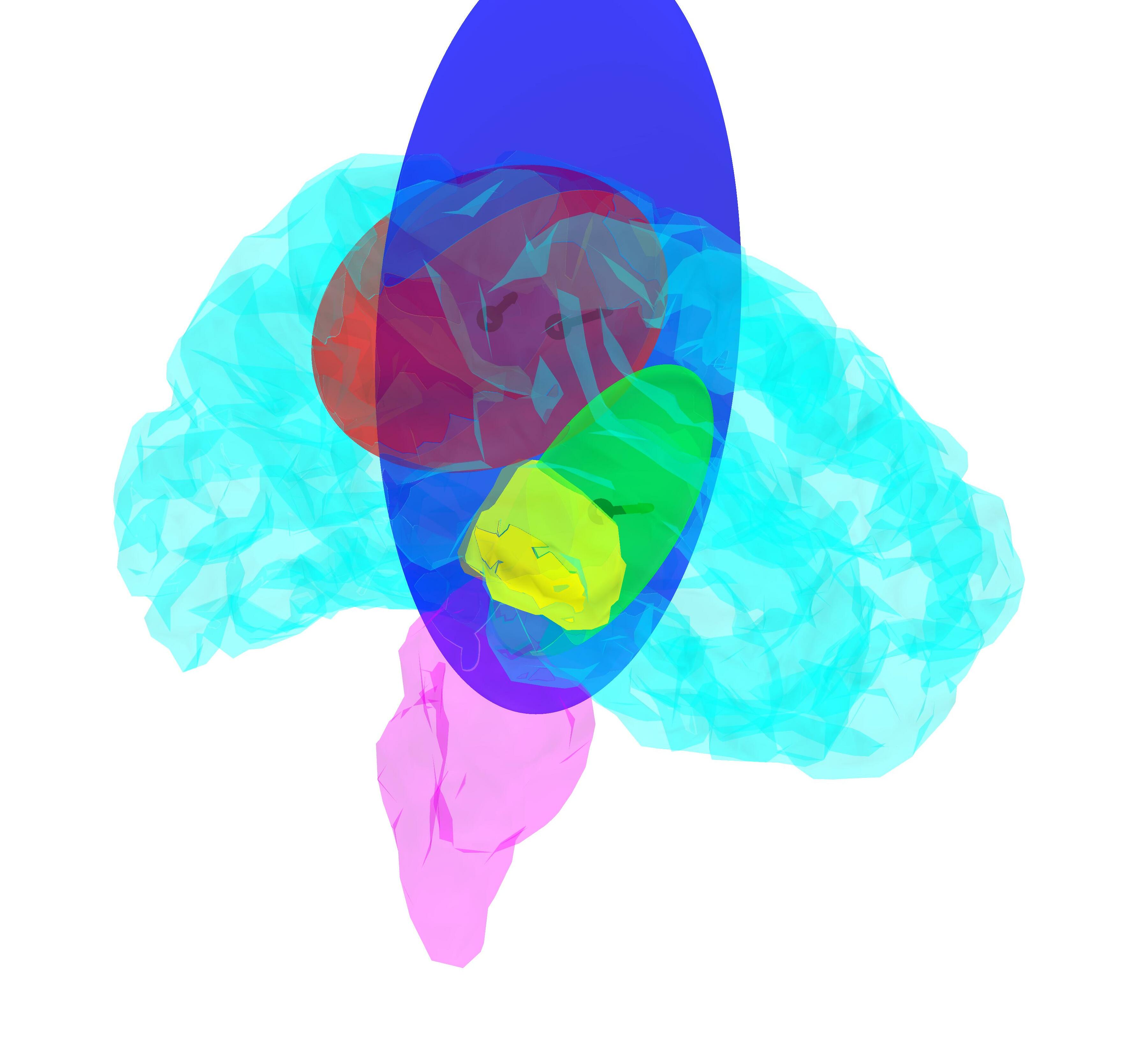}
                \end{subfigure}
                \\ 
                \begin{subfigure}[T]{1.60cm}
                    \centering
                    \includegraphics[width=1.5cm]{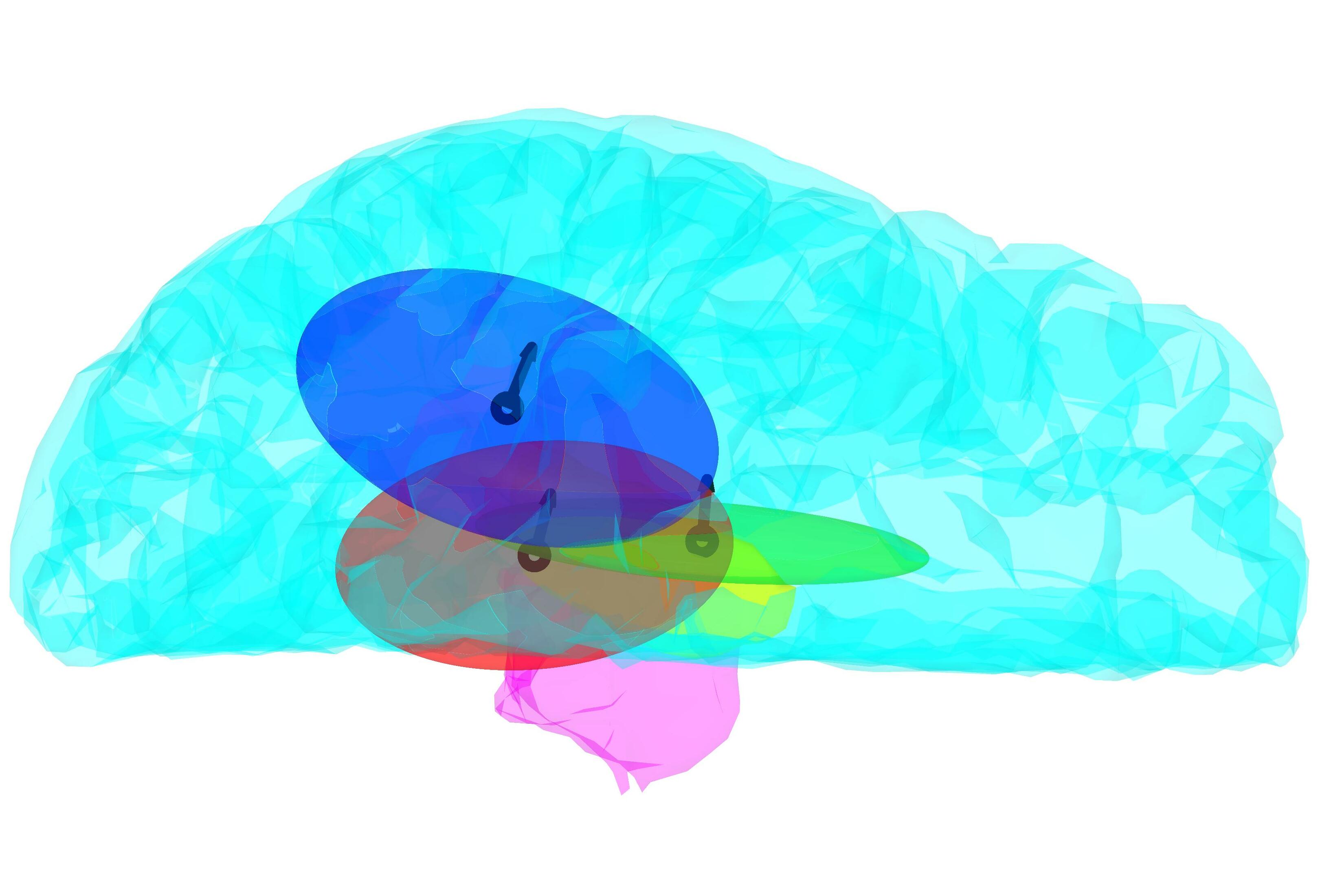}
                    \caption{}
                   \label{fig:GMM_3mm_p22}
                \end{subfigure}
                \begin{subfigure}[T]{1.60cm}
                    \centering
                    \includegraphics[width=1.5cm]{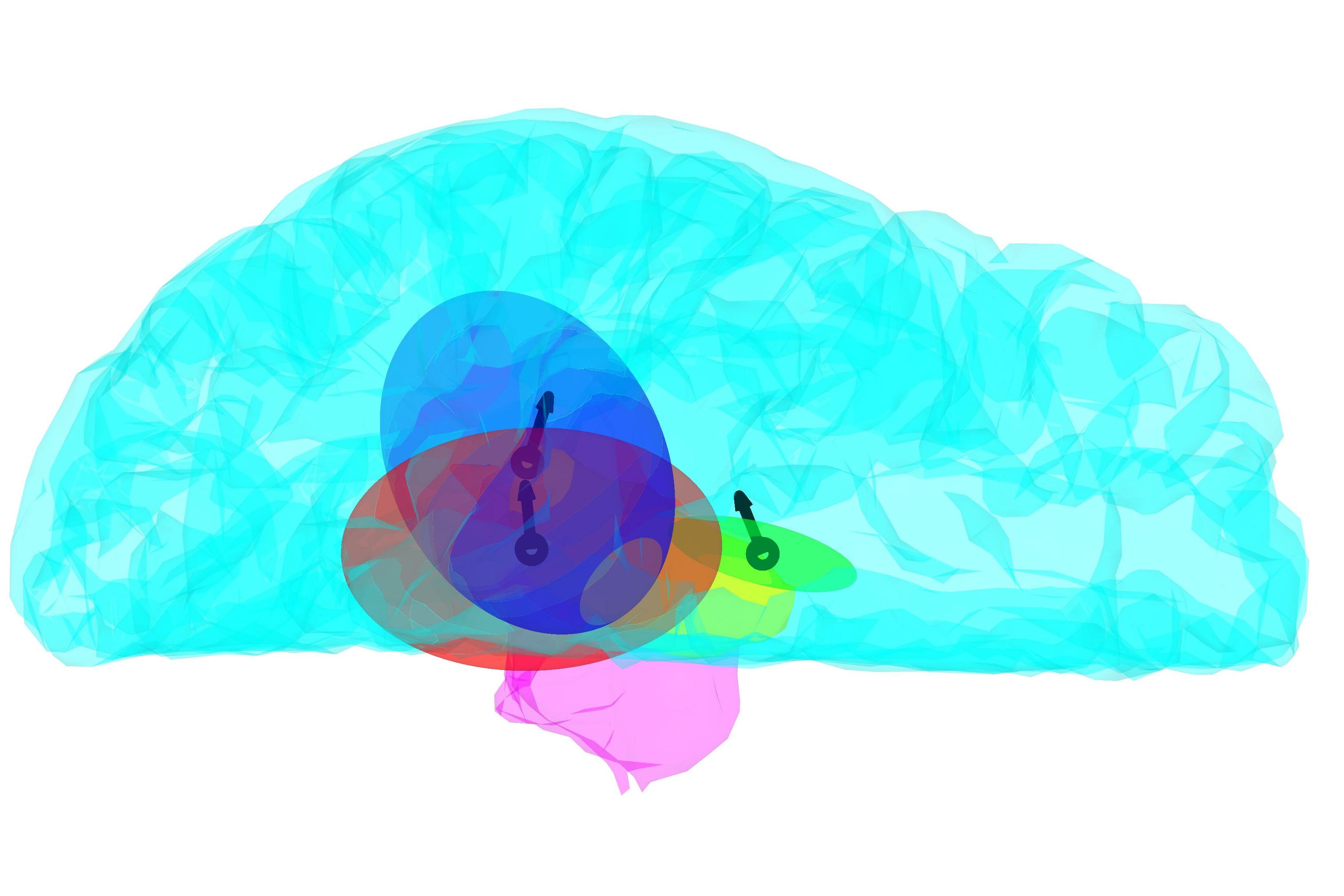}
                    \caption{}
                    \label{fig:GMM_2mm_p22}
                \end{subfigure}
                \begin{subfigure}[T]{1.60cm}
                    \centering
                    \includegraphics[width=1.5cm]{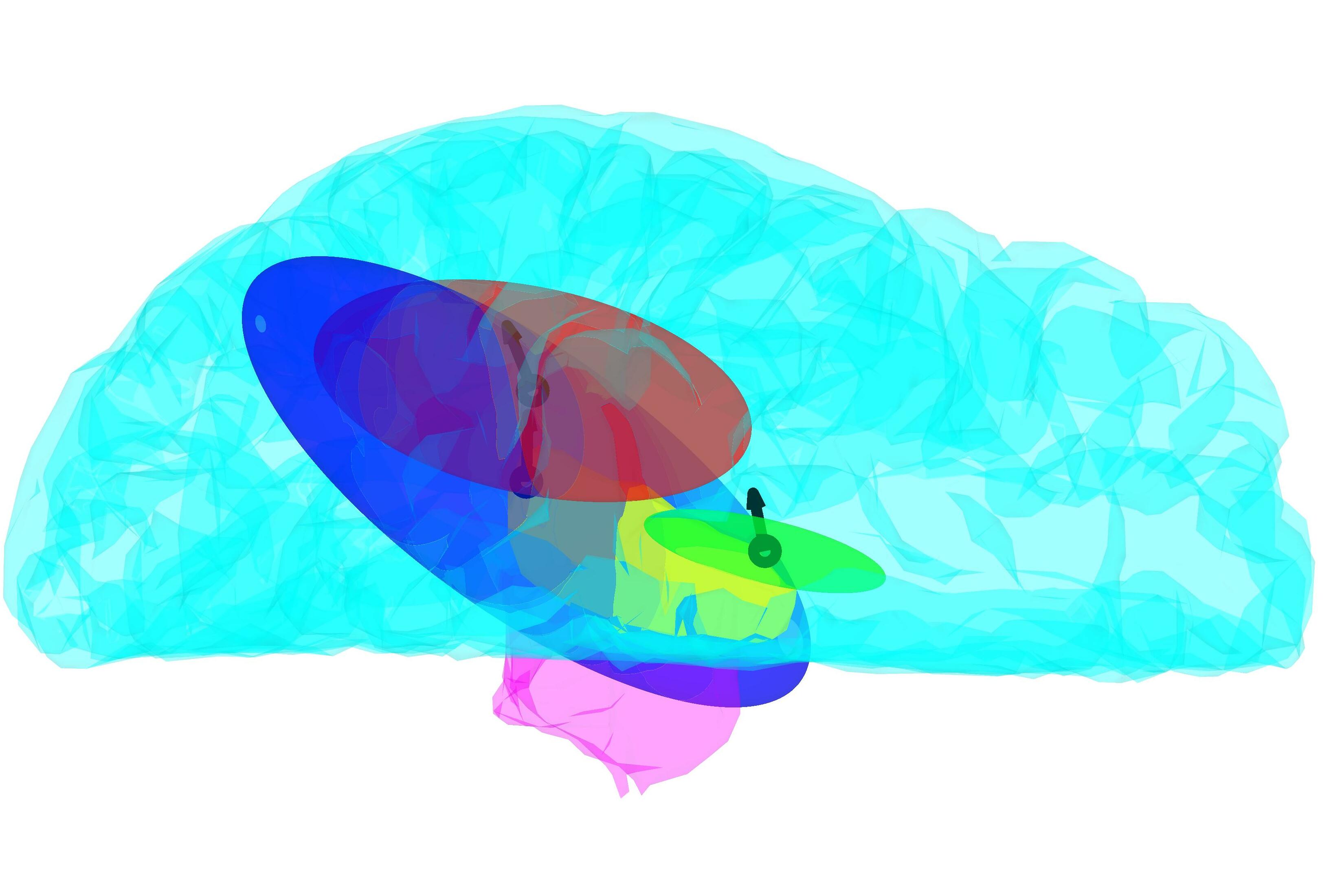}
                    \caption{}
                    \label{fig:GMM_130mm_p2}
                \end{subfigure}
                \begin{subfigure}[T]{1.60cm}
                    \centering
                    \includegraphics[width=1.5cm]{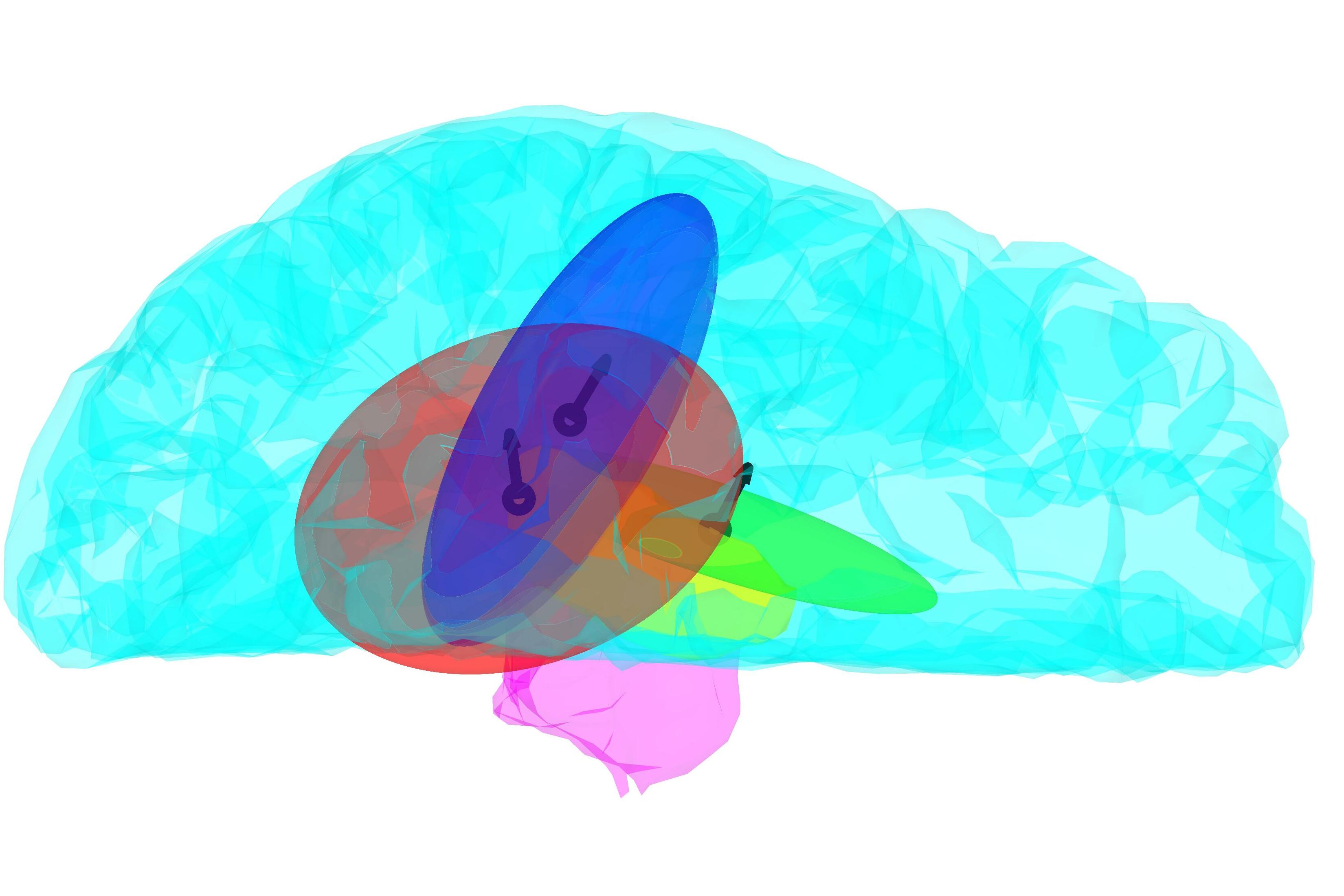}
                    \caption{}
                    \label{fig:GMM_1mm_p22}
                \end{subfigure}
            \end{minipage}
    \end{scriptsize}
\caption{Gaussian mixture modelling (GMM) clusters obtained for the P14/N14 (\ref{fig:GMM_3mm_p14}--\ref{fig:GMM_1mm_p14}) and P22/N22 (\ref{fig:GMM_3mm_p22}--\ref{fig:GMM_1mm_p22}) components of the experimental SEP dataset  \cite{piastra_maria_carla_2020_3888381, rezaei2021reconstructing}. Each cloud visualized shows the 90 \% credibility set of the corresponding cluster. The clouds  are color-labeled according to their measured intensities levels in cubic/millimeter (mm\textsuperscript{3}).}
\label{fig:GMM_p14_p22}
\end{figure}

The GMM clusters corresponding to the SEP components P14/N22 and P22/N22 are shown in Figure \ref{fig:GMM_p14_p22}. Table \ref{table:GMM_Vol} describes the cluster volumes and their best-fitting cluster counts. Clusters have been color-labeled according to their intensities (R = red, G = green, B = blue) in descending order. 

For P14/N14, a single cluster concentrated at the brain stem was found in the case of the adapted 1.3 mm mesh, which is in accordance with the physiological knowledge of the originator (see Section \ref{sec:seps}). Two clusters were found with adapted 2.0 and 3.0 mm meshes, while the regular one resulted in three reconstructed clouds, which is partially contrary to the knowledge of a single originator. The weaker clouds were larger in size, suggesting that they are due to modeling errors, as the actual originator is known to be well-localized in the brain stem.

For P22/N22, three clouds were obtained with each FE mesh. The distinction of the two most intense clouds, R and G, into cortical and thalamic components, was the clearest with 1.3 mm mesh. In the other meshes, the clusters were closer to each other or overlapped; in particular, the centroid of the most intense cluster was located deeper in the brain, suggesting that the reconstructed cortical and sub-cortical activity were partially mixed in the distribution found by sLORETA. The weaker clusters were larger in size, suggesting that they corresponded to modeling errors.

\begin{table*}[h!]
\centering
    \begin{scriptsize}
        \caption{Total computation time (seconds) in CPU (top) and GPU parallelized (bottom) for the FE mesh creation of the spherical three-layer Ary model including post-processing methods, and Matlab's data handling. }
        \label{table:cpu_gpu_sphere}
        \begin{tabular}{lrrrr}
            \toprule
            \multicolumn{1}{c}{} & \multicolumn{3}{c}{\bf Adapted} & \multicolumn{1}{c}{\bf Regular} \\
            \multicolumn{1}{l}{\textbf{Function}} & \textbf{3.0mm} & \textbf{2.0mm} & \textbf{1.3mm} & \textbf{1.0mm} \\
            \hline
            {Surface extraction} & 32   & 114  & 335   & 249   \\ 
            {Labeling}           & 95   & 501  & 4810  & 7450  \\ 
            {Refinement}         & 7    & 23   & 64    & 0    \\ 
            {Smoothing} 	     & 62   & 50   & 151   & 0    \\ 
            {Inflation}          & 23   & 56   & 162   & 0    \\ 
            {Optimization}       & 15   & 223  & 801   & 0    \\ 
            {Data handling}      & 64   & 197  & 579   & 690  \\ 
            \hline
            \textbf{Total}       & 298  & 1164 & 6902  & 8389 \\ 
            \\
            \vspace{0.03cm}
            \\
            {Surface extraction} & 32  & 116  & 354  & 247   \\ 
            {Labeling}           & 17  & 66   & 424  & 597  \\ 
            {Refinement}         & 8   & 24   & 65   & 0    \\ 
            Smoothing 	         & 15   & 51  & 153  & 0    \\ 
            Inflation            & 23   & 59  & 162  & 0    \\ 
            Optimization         & 62   & 243 & 777  & 0    \\ 
            {Data handling}      & 64   & 195 & 581  & 690  \\ 
            \hline
            \textbf{Total}       & 221  & 754 & 2516 & 1534 \\
        \end{tabular}
    \end{scriptsize}
\end{table*}

\begin{table*}[h!]
\centering
    \begin{scriptsize}
        \caption{Total computation time (seconds) in CPU (top) and GPU parallelized (bottom) for the FE mesh creation of the realistic head model including post-processing methods, and Matlab's data handling.}
        \label{table:cpu_gpu_head}
        \begin{tabular}{lrrrr}
            \toprule
            \multicolumn{1}{c}{} & \multicolumn{3}{c}{\bf Adapted} & \multicolumn{1}{c}{\bf Regular} \\
            \multicolumn{1}{c}{\textbf{Function}} & \textbf{3.0mm} & \textbf{2.0mm} & \textbf{1.3mm} & \textbf{1.0mm} \\
            \hline
            {Surface extraction} & 303  & 1124 & 2640   & 1493  \\ 
            {Labeling}           & 546  & 1528 & 11585  & 23947 \\ 
            {Refinement}         & 11   & 39   & 195    & 0     \\ 
            {Smoothing} 	     & 23   & 75   & 261    & 0     \\ 
            {Inflation}          & 90   & 247  & 714    & 0     \\ 
            {Optimization}       & 147  & 528  & 1925   & 0     \\ 
            {Data handling}      & 70   & 222  & 1711   & 642   \\ 
            \hline
            \textbf{Total}       & 1190 & 3763 & 19031  & 26082 \\ 
            \\
            \vspace{0.03cm}
            \\
            {Surface extraction} & 205  & 767  & 2416   & 1158 \\ 
            {Labeling}           & 184  & 451  & 1732   & 1706 \\ 
            {Refinement}         & 10   & 44   & 273    & 0    \\ 
            {Smoothing} 	     & 14   & 44   & 161    & 0    \\ 
            {Inflation}          & 78   & 202  & 553    & 0    \\ 
            {Optimization}       & 223  & 668  & 2084   & 0    \\ 
            {Data handling}      & 77   & 228  & 745    & 636  \\ 
            \hline
            \textbf{Total}       & 791  & 2404 & 7946   & 3500 \\ 
        \end{tabular}
    \end{scriptsize}
\end{table*}

\begin{table}
\centering
    \begin{scriptsize}
    \caption{Volume-value of the up-to-three obtained GMM-based clusters (R=Red, G=Green, B=Blue) ordered in descending order with respect to their intensity, and  measured in cubic millimeters (mm\textsuperscript{3}).}
    \label{table:GMM_Vol}
        \begin{tabular}{rrrr}
        \toprule
        \multicolumn{4}{c}{\bf P14/N14 component } \\
        \multicolumn{1}{c}{\textbf{Mesh size}} & \textbf{R} & \textbf{G} & \textbf{B}\\
        \midrule
        Adapted 3.0 mm   & 889   & 34104  & -    \\
        Adapted 2.0 mm   & 889   & 34104  & -    \\
        Adapted 1.3 mm   & 455   & -      & -    \\
        Regular 1.0 mm   & 393   & 333    & 6537 \\
        \bottomrule
        \multicolumn{4}{c}{\em{(a)} Volume of the P14/N14 clusters. }
        \end{tabular}
        
        \vskip0.4cm
        
        \begin{tabular}{rrrr}
        \toprule
        \multicolumn{4}{c}{\bf P22/N22 component } \\
        \multicolumn{1}{c}{\textbf{Mesh size}} & \textbf{cR} & \textbf{cG} & \textbf{cB}\\
        \midrule
        Adapted 3.0 mm   & 2612   & 4581 & 5697  \\
        Adapted 2.0 mm   & 2606   & 10   & 22317 \\
        Adapted 1.3 mm   & 2353   & 23   & 46303 \\
        Regular 1.0 mm   & 9915   & 824  & 31988 \\
        \bottomrule
        \multicolumn{4}{c}{\em{(b)} Volume of P22/N22 clusters.}
        \end{tabular}
    \end{scriptsize}
\end{table}

\section{Discussion}
\label{sec:discussion}
This study has demonstrated the open-source Matlab-based Zeffiro Interface \cite{he2020zeffiro} capabilities to automatically create a FE mesh \cite{de_munck_wolters_clerc_2012, braess2007finite} for a multi-compartment model of the human head, including both cortical and sub-cortical compartments. The realistic meshes generated were shown to match the given complete head segmentation obtained using FreeSurfer Software Suite \cite{freesurfer}.

Our method is based on nested compartment structures allowing a robust mesh generation for an arbitrary set of segmentation boundaries. One advantage of our method is the capability to generate a mesh regardless of intersecting segmentation boundaries, inhibiting the issues that can follow from intersecting surfaces when using a standard heuristic FE mesh generator \cite{hang2015tetgen, schoberl1997netgen, geuzaine2009gmsh}. 

The resolution of the FE mesh in the vicinity of the tissue boundaries is essential to reduce forward modeling and source localization errors. The recursive labeling approach vastly reduced the computational cost of generating the meshes, thus resulting in high-quality FE simulations with accurate segmentations and compartments. GPU-accelerated recursive labeling further sped up the mesh generation compared to a CPU-only implementation. 

We analyzed the mesh generation and post-processing results for three different initial hexahedral lattice resolutions of 3.0, 2.0, and 1.3 mm. A regular mesh ---a mesh without post-processing effects, i.e., no smoothing, inflation or optimization steps,--- was also generated for comparison. A mesh size of about 1.0 mm is crucial for obtaining a rigorous EEG forward model, e.g., in \cite{rullmann2009eeg}. The regular mesh had the largest spread than the adapted ones, reflected in the measures for EEG modeling accuracy. Due to adaptation, each (non-regular) mesh surpassed the 1.0 mm accuracy considering the distance spread between the given and reconstructed grey matter boundary. Due to increasing computing effort along the initial mesh resolutions, the coarser cases, say, of 3.0 and 2.0 mm, can be considered suitable without GPU acceleration, e.g., with a personal computer or laptop.

In particular, small sub-cortical details, such as the structure of the cingulate cortex, required the finest resolution to be labeled correctly. The finest resolution outperformed the regular 1.0 mm mesh in EEG forward simulation and source localization applications. With this resolution, the median forward accuracy obtained at 98 \% corresponded to 1.1 \% RDM and 2.8 \% MAG. For comparison,  \cite{miinalainen2019realistic, pursiainen2016electroencephalography} obtained approximately 0.3 \% RDM and 0.3 \% MAG with an optimized spherical Stok model \cite{stok1987} using a mesh created with Gmsh. 

The present results have their own advantages and disadvantages, as the Ary model \cite{ary1981location} has been observed to yield up to three times the errors obtained with Stok \cite{pursiainen2011forward}, following from the direct contact between the brain and skull compartment, which in Stok has a CSF layer in between them. Direct comparison between different meshing approaches is also not fully justified; the current strategy sets comparably soft requirements for surface accuracy to avoid overfitting in the realistic case, reflected in the results.

In source localization, the median EMD found (with the finest adapted mesh), 8--10 mm for the superior locations and 12--13 mm for the deeper ones, with 2.0 mm spread in each case, which matches appropriately with the experimentally found estimates for the best possible source localization accuracy in a spherical geometry being 9.2 and 12.8 mm with a standard deviation of 4.4 and 6.2 mm, respectively \cite{cuffin2001realistically}. The forward accuracy was found to be superior with regular 1.0 mm mesh, while in source localization they were slightly less robust and lacked some accuracy with MNE and sLORETA, the deterioration being of maximally 4.0 mm considering median values. This deficiency falls in the standard deviation, suggested for the experimental case \cite{cuffin2001realistically}, and can be interpreted as a lattice or interpolation effect \cite{bauer2015, pursiainen2016electroencephalography}, where sources that are distant from each other than in a denser mesh are used in interpolating a given source. Namely, the differences observed are close to the initial mesh resolution in size.

The GMM clusters, including count, position, and size, obtained using the experimental SEP data \cite{piastra_maria_carla_2020_3888381}, suggest that the resolution of the mesh can have a significant effect on the inverse estimates obtained with a realistic setting considering all the mesh resolutions examined in this study. Reflecting the physiological knowledge of the number, size, and positioning of the GMM clusters found \citep{noel1996origin, buchner1995somatotopy,  mauguiere1983neural, urasaki1990origin, passmore2014origin, buchner1994source, allison1991cortical,  fuchs1998improving, papadelis2011ba3b}, a high resolution and mesh adaptation seem potentially beneficial ways to improve the reconstructions of P14/N14 and P22/N22 components. This observation effect was concerned, especially the adapted 1.3 mm case, and is in line with the other results of this study. Because of the sensitivity to modeling errors, we encountered deep structure-activity of the P14 component to deviate between the lead fields obtained with different meshes. The sensitivity of P22/N22 not only follows its depth contribution but also from the simultaneous cortical originator in the pre-or post-central gyrus. The GMM deviations induced by the different FE meshes highlight the importance of accurate forward modeling that even a few percent differences in RDM or MAG can lead to significant effects when reconstructing weakly distinguishable activity. Previously, the detectability of the P14/N22 and P22/N22 originators has been suggested and analyzed in \cite{rezaei2021reconstructing}. The distinguishability of the deep sub-cortical components based on non-invasive experimental EEG data suggested recently \cite{seeber2019subcortical}.

Our present approach constitutes an independent, time-effective tetrahedral mesh generator, free of supplementary pre-processing aspects. Other advanced tetrahedral mesh generators for the brain have been created recently by coupling multiple pre-processing tools with the well-known mesh generators, possibly most prominently, by recurring Matlab's toolbox iso2mesh \cite{fang2009tetrahedral}. In \cite{piastra2021comprehensive}, a tetrahedral discretization of FreeSurfer's Aseg atlas has been obtained using a combination of Computational Geometry Algorithms Library (CGAL) \cite{fabri2009cgal} and the iso2mesh toolbox, resulting in a FE mesh of 0.8 M nodes and 5.3 M tetrahedrons. In \cite{ROAST2018}, the iso2mesh-based mesh generator has been introduced for transcranial electrical stimulation (tES). Compared to SimBio-Vgrid driven studies \cite{dannhauer2011modeling, antonakakis2019effect}, we used a higher number of compartments, including a complete set of sub-cortical ones. Based on our results, we are confident that a finer mesh resolution might be necessary to completely distinguish the deeper structures as they appear in the Aseg atlas, justifying the examination of denser meshes and GPU acceleration as a way to speed up the labeling process. 

Our future work will include further tests with experimental data as well as methodological and computational considerations, e.g., possibilities to speed up the surface extraction process, which is a major contributor to the total meshing time and involves complex indexing operations. The present meshing approach will be applied in the further development of advanced FEM forward simulation tools such as Duneuro \cite{schrader2021duneuro} or packages for analyzing brain activity, e.g., the current implementation platform ZI \cite{he2020zeffiro}. It can potentially provide future directions for developing similar platforms such as the well-known Brainstorm \cite{tadel2011brainstorm} whose forward model was originally built upon the boundary element method.

\backmatter

\bmhead{Acknowledgments}
{\em Funding.} 
FGP, JL, MS, and SP were supported by the Academy of Finland through the Center of Excellence in Inverse Modelling and Imaging 2018--2025 (336792), DAAD project (334465), and the ERA "Personalised diagnosis and treatment for refractory focal pediatric and adult epilepsy" (PerEpi) project (344712). JL is supported by Väisälä Fund’s (Finnish Academy of Science and Letters) one-year young researcher grant admitted in 2021.

{\em Conflicts of interest.} 
The authors certify that this study is a result of purely academic, open, and independent research. They have no affiliations with or involvement in any organization or entity with a financial interest or non-financial interest such as personal or professional relationships, affiliations, knowledge, or beliefs in the subject matter or materials discussed in this manuscript.

{\em Data availability.}
The datasets used in this manuscript are open access under the license Open Data Commons Attribution License ODC-BY 1.0.

{\em Code availability.}
All scripts used in this manuscript were designed and implemented in Zeffiro Interface (ZI), @ 2018- Sampsa Pursiainen \& ZI Development Team, an open source code package constituting an accessible tool for multidisciplinary finite element (FE) based forward and inverse simulations in complex geometries.

{\em Authors’ contributions.}
FGP wrote the main manuscript text and prepared all figures, JL, MS, and SP contributed significantly to the numerical experiments. All authors reviewed the manuscript.


\bibliography{cas-refs}

\end{document}